\documentclass[11pt]{amsart}
\usepackage{pdfpages}
\usepackage[english]{babel}
\usepackage[utf8]{inputenc} 
\usepackage[T1,OT1]{fontenc} 

\usepackage{amssymb} 
\usepackage{amsmath} 
\usepackage{amsthm} 
\usepackage{graphicx,subfigure}	
\usepackage{fullpage} 
\usepackage{epsfig,color}
\usepackage{enumerate}
\usepackage{natbib}
\usepackage[final]{showlabels}
\usepackage{enumitem}

\usepackage{float} 

\usepackage{tikz}
\usetikzlibrary{shapes}

\usepackage{pgfplots}
\usepackage{todonotes}

\definecolor{green_mermoz}{rgb}{0,0.5,0}
\usepackage{hyperref}

\hypersetup{colorlinks = true, linkcolor=green_mermoz, citecolor=orange}

\usepackage[noabbrev]{cleveref}
\creflabelformat{equation}{#2(#1)#3}
\crefformat{equation}{#2(#1)#3}
\crefformat{equations}{#2(#1)#3}
\crefname{thm}{Theorem}{Theorems}
\crefname{ass}{Assumption}{Assumptions}
\crefname{problem}{Problem}{Problems}
\crefname{prop}{Proposition}{Propositions}
\crefname{lem}{Lemma}{Lemmas}
\creflabelformat{problem}{#2(#1)#3}

\definecolor{k4}{rgb}{0.8,0.8,0.8}
\definecolor{k3}{rgb}{0.6,0.6,0.6}
\definecolor{k2}{rgb}{0.4,0.4,0.4}
\definecolor{k1}{rgb}{0.2,0.2,0.2}

\newcommand{\N}{\mathbb{N}}

\newcommand{\R}{\mathbb{R}}

\newcommand{\eps}{\varepsilon}

\newcommand{\p}{\partial}

\newcommand{\abs}[1]{\left\lvert#1\right\rvert}
\newcommand{\norm}[1]{\left\lVert#1\right\rVert}

\newcommand{\normf}[1]{\left\lVert#1\right\rVert_\F}

\newcommand{\nast}[1]{\left\lVert#1\right\rVert_\ast}

\newcommand{\bz}{\overline{z}}

\newcommand{\E}{\mathcal{E}}

\newcommand{\Va}{V^\ast}
\renewcommand{\epsilon}{\varepsilon}

\DeclareMathOperator{\I}{I}

\DeclareMathOperator{\zs}{z_\ast}
\DeclareMathOperator{\zds}{\dot{z}_\ast}
\DeclareMathOperator{\zopt}{z_\ast^s}

\newcommand{\D}{\mathcal{D}_\e}
\renewcommand{\I}{\mathcal{I}}

\renewcommand{\leq}{\leqslant}
\renewcommand{\geq}{\geqslant}

\renewcommand{\eps}{\varepsilon}

\renewcommand{\tilde}{\widetilde}

\def\F{\mathcal{F}}

\newcommand{\sign}{\mathrm{sign}}

\newcommand{\e}{\varepsilon}

\newcommand{\vphia}{\varphi_\alpha}

\newcommand{\Oun}{O^\ast(1)+O(\e)}

\newcommand{\Vea}{V_\e^\ast}
\newcommand{\qea}{q_\e^\ast}

\newcommand{\qdea}{\dot{q}_\e^\ast}

\newcommand{\pde}{\dot{p}_\e}
\newcommand{\qde}{\dot{q}_\e}

\newcommand{\kde}{\dot{\kappa}_\e}

\newcommand{\Iea}{\mathcal{I}_\e^\ast}

\newcommand{\East}{\mathcal{E}^\ast}

\newcommand{\qa}{q^\ast}
\newcommand{\pa}{p^\ast}

\newcommand{\qda}{\dot{q}^\ast}
\newcommand{\pda}{\dot{p}^\ast}

\renewcommand{\Vea}{\Va}
\renewcommand{\qea}{\qa}
\renewcommand{\qdea}{\qda}

\newtheorem{thm}{Theorem} 
\newtheorem{lem}[thm]{Lemma} 
\newtheorem{prop}[thm]{Proposition} 
\newtheorem{remark}[thm]{Remark} 
\newtheorem{cor}[thm]{Corollary}
\newtheorem{definition}[thm]{Definition}
\newtheorem{ass}[thm]{Assumption} 
\numberwithin{equation}{section}
\numberwithin{thm}{section}

\def\ds{\displaystyle}
\usepackage{epstopdf}

\title{The Cauchy problem for the infinitesimal model in the regime of small variance}

\author{Florian Patout} 
\address{UMPA, UMR 5669 CNRS \& Ecole Normale Sup\'erieure de Lyon, Lyon, France} 
\email{florian.patout@ens-lyon.fr}
\date{\today}

\begin{document}

\begin{abstract}
	We study the asymptotic behavior of solutions of the Cauchy problem associated to a  quantitative genetics model with a  sexual mode of reproduction. It combines  trait-dependent mortality and a nonlinear integral reproduction operator "the infinitesimal model" with a parameter describing the standard deviation between the offspring and the mean parental traits. We show that under mild assumptions upon the mortality rate $m$, when the deviations are small, the	 solutions stay close to a Gaussian profile with small variance, uniformly in time.  Moreover we characterize accurately the dynamics of the mean trait in the population. Our study extends previous results on the existence and uniqueness of stationary solutions for the model. It relies on perturbative analysis techniques together with a sharp description of the correction measuring the departure from the Gaussian profile.
\end{abstract}
\maketitle
\section{Introduction}\label{sec intro}
We investigate solutions $f_\e \in L^1(\R_+\times \R)$ of  the following Cauchy problem: 
\begin{equation}\label[problem]{Ptfeps}\tag{$P_tf_\e$}
\left\{\begin{array}{ll}
\e^2 \p_t	 f_\e(t,z) + m(z) f_\e(t,z) = \mathcal{B}_\e(f_\e)(t,z), \quad t>0, \ z \in \R,  \\

f_\e(0,z) = f_\e^0(z).
\end{array} \right.
\end{equation}
where $\mathcal{B}_\e (f)$ is the  following nonlinear, homogeneous mixing operator associated with the infinitesimal model \cite{Fisher1918}, see also \cite{weberinfinitesimal} for a modern perspective :
\begin{equation}
\mathcal{B}_\e (f)(z) := \dfrac{1}{\eps \sqrt{\pi}}   \iint_{\mathbb{R}^{2}}  \exp\left[ -  \dfrac{1}{\eps^2} \left( z - \dfrac{ z_1 +   z_2}2 \right)^2 \right ]  f( z_1)\dfrac{  f(  z_2)}{ \int_{\R}   f(  z_2')\, d  z_2'}\, d  z_1 d  z_2.\label{eq:Beps_cauchy}
\end{equation}
This problem originates from quantitative genetics in the context of evolutionary biology. The variable $z$ denotes a phenotypic trait, $f_\e$ is the distribution of the population with respect to $z$ and $m$ is the trait-dependent mortality rate.

The mixing operator $\mathcal{B}_\e$ models  the inheritance of quantitative traits in the population, under the assumption of a sexual mode of reproduction. As formulated in \cref{eq:Beps_cauchy}, it is assumed that offspring traits are distributed normally around the mean of the parental traits $(z_1 + z_2)/2$, with a constant variance, here $\eps^2/2$.  We are interested in the evolutionary dynamics resulting in the selection of well fitted (low mortality) individuals \textit{i.e.} the concentration of the distribution around some dominant traits with standing variance.

In theoretical evolutionary biology, a broad literature deals with this model to describe sexual reproduction, see {\em e.g.} \cite{slatkin1970selection,
	roughgarden1972evolution,
	slatkin_niche_1976,
	bulmerbook,menormal,
	tufto_quantitative_2000,	
	barfield_evolution_2011, huisman2012comparison,
	cotto-ronce,
	weberinfinitesimal,turelli_comm}.
 
We are interested in the asymptotic behavior of the trait distribution $f_\e$ as $\eps^2$ vanishes. It is expected that the profile concentrates around some particular traits under the influence of selection.

The asymptotic description of concentration around some particular trait(s) has been extensively investigated for various linear operators $\mathcal{B}_\e$ associated with asexual reproduction such as, for instance, the diffusion operator $f_\e(t,z) + \e^2 \Delta f_\e(t,z)$, or the convolution operator $\frac1\e K(\frac z\e)*f_\e(t,z)$ where $K$ is a probability kernel with unit variance, see \cite{diekmann2005,perthame-book,barlesperthamecontemporary,
	barles2009,
	lorz2011dirac}  for the earliest investigations, and  \cite{meleard_singular_2015,
	mirrahimi2018singular,bouin_thin_2018} for the case of  fat-tailed kernel $K$.  
In those linear cases, the asymptotic analysis usually leads to a Hamilton-Jacobi equation after performing the Hopf-Cole transform $u_\e = -\e \log f_\e$. Those problems require a careful well-posedness analysis for uniqueness and convergence as $\e \to 0$ see: \cite{barles2009,mirrahimi_class_2015,
	calvez_uniqueness_2018}.


Much less is known about the operator $\mathcal{B}_\e$ defined by \cref{eq:Beps_cauchy}. From a mathematical viewpoint, in the field of probability theory \cite{weberinfinitesimal} derived the model from a microscopic framework.  In \cite{raoulmirrahimiinfinitesimal,Raoul}, the authors deal with a different scaling than the current  small variance assumption $\e^2 \ll 1$ , and add a spatial structure in order to derive the celebrated Kirkpatrick and Barton system \cite{KirBar97}. \vspace{1cm}

\begin{figure}
	\begin{center}
		\includegraphics[scale=0.83,page=1]{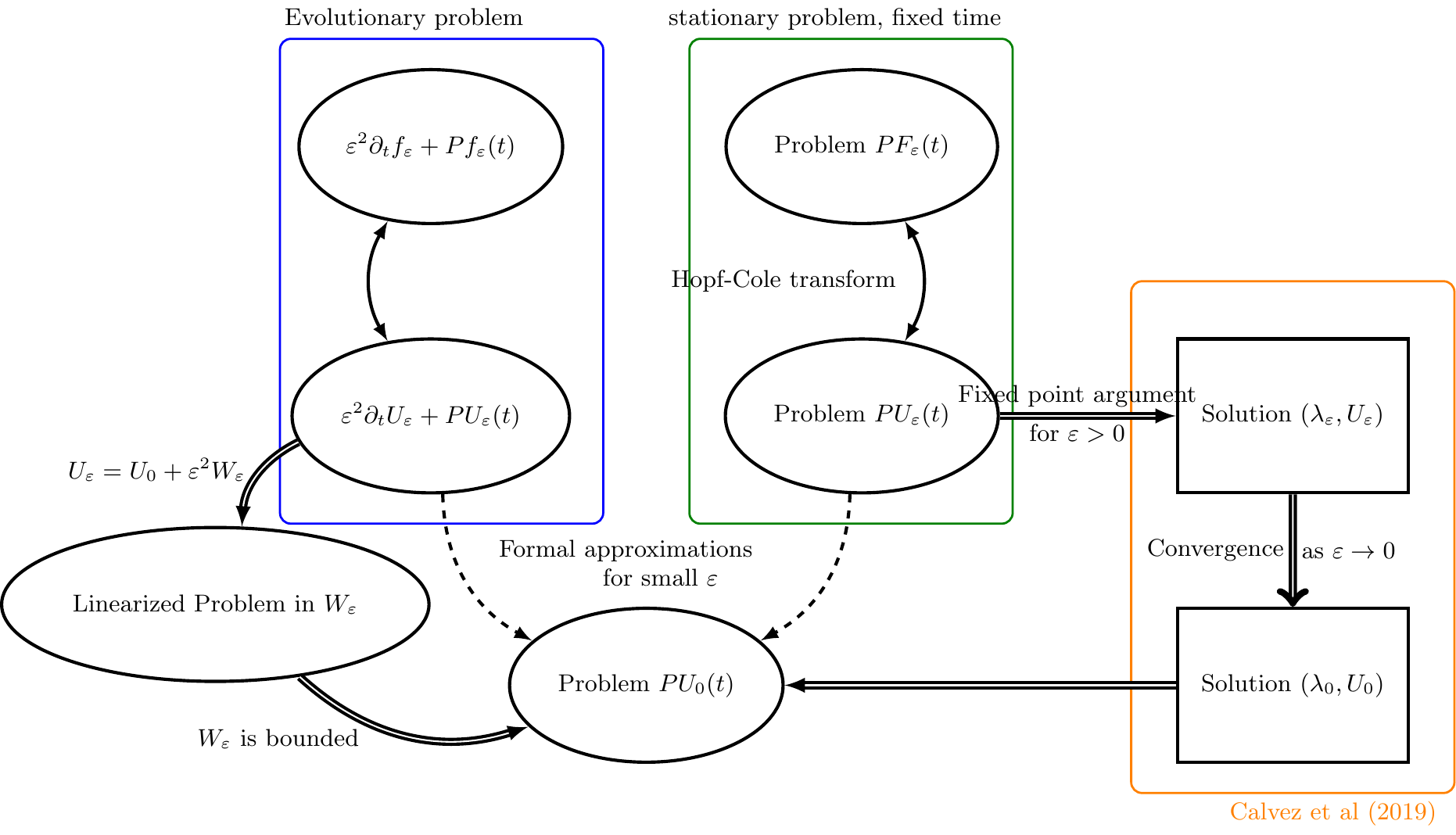}
	\end{center}
	\caption{Scope of our paper compared to precedent work	}
	\label{fig scheme}
\end{figure}

Gaussian distributions will play a pivotal role in our analysis as they are left invariant by the infinitesimal operator $\mathcal{B}_\e$, see \cite{menormal,raoulmirrahimiinfinitesimal}. In \cite{spectralsex}, the authors studied  special "stationary" solutions, having the form: 
\begin{align*}
\exp \left(\dfrac{\lambda_\e  \, t }{\e^2}  \right) F_\e(z),  \quad \text{ with } F_\e(z)= \frac{1}{\e \sqrt{ 2\pi}}\exp \left( -\frac{(z-\zs)^2}{2 \e^2} - U_\e^s(z) \right).
\end{align*} 
In this paper we tackle the Cauchy \cref{Ptfeps}, and we hereby look for solutions that are close to Gaussian distributions uniformly in time:
\begin{align}\label{hopfcole}
f_\e(t,z) = \frac{1}{\e \sqrt{ 2\pi}}\exp \left( \dfrac{\lambda(t)}{\e^2} - \dfrac{(z - \zs(t) )^2}{2\epsilon^2}  - U_\epsilon(t,z) \right). 
\end{align}
The scalar function $\lambda$  measures the growth (or decay according to its sign) of the population.  The mean of the Gaussian density, $\zs$ is also the trait at which the population concentrates when $\e \to 0$. The pair $(\lambda,\zs)$ will be determined by the analysis at all times. It is somehow related to invariant properties of the operator $\mathcal{B}_\e$. The function $U_\e$ measures the deviation from the Gaussian profile induced by the selection function $m$. It is a cornerstone of our analysis that $U_\e$ is Lipschitz continuous with respect to $z$, uniformly in $t$ and $\e$.  
Plugging the transformation \cref{hopfcole} into \cref{Ptfeps} yields the following equivalent one: 
\begin{multline}\label[problem]{PtUeps}\tag{$P_tU_\e$}
-\e^2 \p_t U_\e(t,z) + \dot{\lambda} (t) + (z-\zs(t) ) \zds(t) + m(z) = \\ I_\e(U_\e)(t,z) \exp \Big( U_\e(t,z)-2U_\e\left(t, \bar{z}(t) \right) + U_\e(t,\zs(t) ) \Big), 
\end{multline}
where $\bar{z}(t)$ is the midpoint between $z$ and $\zs(t)$ :
\begin{align*}
\bar{z}(t) = \dfrac{z + \zs(t)}{2},
\end{align*}
and the functional $I_\e$ is defined by
\begin{multline}\label{def_Ieps}
I_\e(U_\e)(t,z) = \\  \dfrac{\ds  \iint_{\R^2} \exp \left[-\dfrac{1}{2}y_1 y_2 - \dfrac{3}{4}\left (y_1^2+y_2^2\right ) +2 U_\e\left(t, \bz \right) - U_\e\left(t, \bz+\epsilon y_1 \right) - U_\e\left(t, \bz +\epsilon y_2\right)\right] d y_1 d y_2  }
{\ds   \sqrt{\pi} \int_\R \exp\left[ - \frac12 y^2 +  U_\e(t,\zs)- U_\e(t,\zs+\epsilon y) \right] d y}. 
\end{multline}
This functional is the residual shape of the infinitesimal operator \cref{eq:Beps_cauchy} after suitable transformations. It was first introduced in the formal analysis of \cite{main} and in the study  of the corresponding stationary problem in \cite{spectralsex}. The Lipschitz continuity of $U_\e$ is pivotal here as it ensures that $I_\e(U_\e) \to 1$ when $ \e \to 0$. Thus for small $\e$, we expect that the \cref{Ptfeps}  is well approximated by the following one : 
\begin{align}\label{PtU0int}
\dot{\lambda}(t) + (z-\zs(t))  \zds(t) + m(z) =  \exp \Big( U_0(t,z)-2U_0\left(t, \bar{z}(t) \right) + U_0(t, \zs(t) ) \Big).
\end{align}
Interestingly,  this characterizes the dynamics of $(\lambda(t),\zs(t))$.  By differentiating \cref{PtU0int} and evaluating at the point  $z = \zs(t)$, then simply evaluating \cref{PtU0int} at $z=\zs(t)$, we find the following pair of relationships:
\begin{align}
\label{eqdef zs} \zds(t) + m' (\zs(t)) = 0,\\
\label{eqdef lambda} \dot{\lambda}(t) + m(\zs(t)) =  1.
\end{align}
Then, a more compact way to write the limit problem for $\e=0$ is 
\begin{align}\label[problem]{PtU0}\tag{$P_tU_0$}
M(t,z) =  \exp \Big( U_0(t,z)-2U_0\left(t, \bar{z}(t) \right) + U_0(t, \zs(t) ) \Big), 
\end{align}
with the notation
\begin{align}\label{def Gamma}
M(t,z):= 1+m(z)-m(\zs(t))-m'(\zs(t))(z-\zs(t)).
\end{align}
It verifies from \cref{eqdef lambda,eqdef zs}:
\begin{equation}\label{rel Mzs}
M(t,\zs(t)) =1, \quad \p_z M(t,\zs(t))= 0.
\end{equation}
An explicit solution of \cref{PtU0} exists under the form of an infinite series:
\begin{equation}\label{def Va}
\Va(t,z):= \sum_{k\geq 0} 2^k \log \Big(  M\left(t,\zs(t) +2^{-k}(z-\zs(t)) \right)\Big).
\end{equation}
Interestingly this series is convergent thanks to the relationships of \cref{rel Mzs}. The function $\Va$ is a solution of \cref{PtU0}, but not the only one. There are two degrees of freedom when solving \cref{PtU0}, since adding any affine function to $U_0$ leaves the right hand side unchanged. Therefore, a  general expression of solutions is the following, where the scalar functions $p_0$ and $q_0$ are arbitrary:
\begin{equation}\label{U0}
U_0(t,z) = p_0(t) + q_0(t)  (z-\zs(t)) + \Va(t,z). 
\end{equation} 
We have foreseen that the Lipschitz regularity of $U_\e$  was the way to guarantee that $I_\e(U_\e) \to 1$ as $\e \to 0$. As a matter of fact, an important part of \cite{spectralsex} is dedicated to prove such regularity for $U_\e^s$ the solution of the stationary problem : 
\begin{align}\label[problem]{PUe stat}\tag{$PU_\e$ stat}
\lambda_\e^s + m(z) & = I_\e(U_\e^s)(z) \exp \left( U_\e^s(z)-2U_\e^s\left(\frac{z+\zopt}{2} \right) + U_\e^s(\zopt) \right),  \quad z \in \R.
\end{align}
The authors  introduced an appropriate functional space controlling Lipschitz bound. They were then able to show the existence of  $U_\e^s$ and its (local) uniqueness in that space. They also proved that $U_\e^s$ was converging when $\e \to 0$ towards solutions of \cref{PtU0}, see \Cref{fig scheme} for a schematic comparison of the scope of the present article article compared to previous work. 

Here, to tackle the non stationary \cref{PtUeps}, we make the following assumptions of asymptotic growth on the selection function $m$,  when $\abs{z} \to \infty$.
\begin{ass}\label{def m cauchy}
	$ $\\
	We suppose that the function $m$ is a $\mathcal{C}^5(\R)$ function, bounded below. We define the scalar function $\zs$ as the following gradient flow:
	\begin{align}\label{gradient flow}
	\zds(t) = - m'(\zs(t)), \quad  t>0,
	\end{align}
	associated to an initial data $\zs(0)$ prescribed. Next, we make the following assumptions :
	\begin{itemize}[label=$\triangleright$]
		\item We suppose that $\zs(0)$ lies next to a non-degenerate \emph{local} minimum of $m$, $\zopt$ such that		\begin{align}\label{zs infty}
		\zs(t)\xrightarrow[t \to \infty]{} \zopt.
		\end{align}
		\item We also require that there exists  a uniform positive lower bound on $M$ : 
		\begin{align}\label{cond Gamma}
	\inf_{(t,z) \in \R_+\times \R} M(t,z) >0.
		\end{align}
		\item We make growth assumptions on $M$ in the following way:
		\begin{align}\label{decay Gamma}
		\text{for } \,  k = 1,2,3,4,5 : \quad 	(1+\abs{z-\zs})^\alpha \dfrac{\p_z^k M(t,z)}{ M(t,z)} \in L^\infty(\R_+  \times \R)\, ,
		\end{align}
		for some $0<\alpha<1$, the same than in \cref{def F}. 
		\item We make a final assumption upon the behavior of $m$ at infinity, that is roughly that it has superlinear growth, uniformly in time:
		\begin{align}\label{cond m infty}
		\limsup_{z \to \infty} \abs{\frac{ M(t, \bz)}{M(t,z)}} :=a < \frac12, \quad \limsup_{z \to \infty} \abs{\frac{ \p_z M(t, \bz)}{\p_z M(t,z)}} < \infty .
		\end{align}
	\end{itemize}
\end{ass}
The first assumption on $m$ and $\zs$ guarantee the following local convexity property, at least for times $t$ large enough :
\begin{align}\label{sign pzzm}
\exists \, \mu_0 >0, \ \exists t_0>0, \text{ such that } \forall t\geq t_0, \,		m''(\zs(t)) \geq \mu_0.
\end{align}
\begin{remark}
	Based on the formulation of \cref{PtU0}, the function $M$ must be positive. We require a uniform bound in \cref{cond Gamma} for technical reasons. It corresponds to a global assumption on the behavior of $\zs$ and $m$, that further reduces the choice of $\zs(0)$. This condition holds true for globally convex functions $m$. However we do not want to restrict our analysis to that case, so we suppose more generally that \cref{cond Gamma} is verified. A more detailed discussion on the behavior of the solution whether this condition is verified or not is carried out in \cref{sec num} with some numerical simulations displayed. 
	Moreover, the decay assumptions \cref{decay Gamma} and \cref{cond m infty} hold true if $m$ behaves like a polynomial function at least quadratic as $\abs{z}\to + \infty$.
\end{remark}
The purpose of this work is to rigorously prove the convergence of the solution of \cref{PtUeps} towards a particular solution of \cref{PtU0}.   
Given the general shape of $U_0$, see \cref{U0}, it is natural to decompose $U_\e$ by separating the affine part from the rest: 
\begin{align}\label{decomp Ue}
U_\e(t,z) = p_\e(t) + q_\e(t)(z-\zs(t)) + V_\e(t,z).
\end{align}
We require accordingly that at all times $t>0$, \begin{align*}
V_\e(t,\zs)= \p_z V_\e(t,\zs)= 0,
\end{align*}
which is another way of saying that the pair $(p_\e,q_\e)$ tune the affine part of $U_\e$. The pair $(q_\e,V_\e)$ is  the main unknown of this problem. It is expected that $V_\e$ converges to $\Va$ when $\e \to 0$. Our analysis will be able to determine the limit of $q_\e$ even if it cannot be identified by the problem at $\e=0$. Indeed in \cref{PtU0}, the linear part $q_0$ can be any constant. Our limit candidate for $q_\e$ is $\qa$, that we define as the solution  of the following differential equation
\begin{align} 
\label{def qda}   	\qda(t)   = -m''(\zs(t)) \qa(t) + \frac{ m^{(3)}(\zs(t))}{2}-2 m''(\zs(t))m'(\zs(t)),
\end{align}
corresponding to an initial value of $\qa(0)$.
Moreover we define $\pa$ as the function which verifies for a given $\pa(0)$,
\begin{align}
\label{def pda} \pda(t) =-  m'(\zs(t)) \qa(t) + m''(\zs(t)).
\end{align}
Finally, the function
\begin{equation}\label{def Uast}
U^\ast(t,z):= \pa(t) + \qa(t)(z-\zs(t)) +\Va(t,z) 
\end{equation}
will be our candidate for the limit of $U_\e$ when $\e \to 0$.  The problem for $V_\e$ equivalent to \cref{PtUeps}, using \cref{decomp Ue}, is: 
\begin{multline}\label[problem]{PtVeps}\tag{$P_tV_\e$}
M(t,z)  - \e^2 \Big( \pde(t) + \qde(t)(z- \zs(t)) + m'(\zs(t))q_\e(t)\Big) -\e^2  \p_t V_\e(t,z)  \\ = \I_\e(q_\e,V_\e)(t,z)   \exp \Big( V_\e(t,z)-2 V_\e\left(t, \bar{z}(t) \right) + V_\e(t, \zs(t) ) \Big).
\end{multline}
One can notice that thanks to cancellations the functional $I_\e(U_\e)$ does not depend on $p_\e$, which explains for the most part why we focus upon $(q_\e, V_\e)$. We choose to write $\I_\e(q_\e,V_\e)(t,z) := I_\e(U_\e)(t,z)$  as a functional of both unknowns because we will study variations in both directions.
One of the main difficulties to prove the link between \cref{PtVeps,PtU0} is that formally, the terms with the time derivatives in $q_\e$ and $V_\e$ vanish when $\e \to 0$. This makes our study  belong to the class of singular limit problems.

Before stating our main result we need to define appropriate functional spaces. We first define a reference space $\E$, similar to the one introduced in \cite{spectralsex} for the study of the stationary \cref{PUe stat}. However, compared to that case we will need more precise controls, which is why we introduce a subspace $\F$ with more stringent conditions. 
\begin{definition}[Functional spaces]\label{def F}$ $\\
We define $
\alpha < 2 - \frac{\ln 3}{\ln 2},
$
such that $\alpha \in \left(0,1\right)$ and the corresponding functional space \begin{multline*}
	\E =\Big\{  v \in \mathcal{C}^3(\R_+ \times \R) \, | \, \forall t>0 , \, v(t, \zs(t)) = \p_z v(t, \zs(t))= 0 \Big\}  \\ 
\bigcap	 \left\{v \in \mathcal C^3(\R_+ \times \R)  \left|
	\begin{array}{rrr} \ds  & \abs{ \p_z v(t,z)}  \\
	\ds  & \Big( 1+ \abs{z-\zs(t)} \Big)^\alpha  \abs{\p_z^2 v(t,z) } \\
	\ds &  \Big( 1+ \abs{z-\zs(t)} \Big)^\alpha  \abs{  \p_z^3 v(t,z) } 
	\end{array} \right.  \in L^\infty(\R_+  \times \R) \right\}	\end{multline*}
	equipped with the norm 
	\begin{multline*}
	\norm{v}_\E = \max \left(   \sup_{ (t,z) \in\R_+ \times \R} \abs{ \p_z v(t,z) } , \  \sup_{ (t,z) \in\R_+\times \R} \Big(1+ \abs{z-\zs(t)}\Big)^\alpha  \abs{ \p_z^2  v(t,z) } \right., \\ \left. \sup_{ (t,z) \in\R_+\times \R}  \Big(1+ \abs{z-\zs(t)}\Big)^\alpha \abs{ \p_z^3  v(t,z) }\right).
	\end{multline*}
	We also define the subspace : 
	\begin{align*}
	\F := \E  \cap \left\{v \in \mathcal C^1(\R_+ \times \R)    \left|
	\begin{array}{rrr} \ds  & \abs{ 2 v(t,\bz(t)) -v(t,z) } \\
	\ds  &   \Big(1+ \abs{z-\zs(t)}\Big)^\alpha \abs{ \p_z v(t,\bz(t)) -\p_z v(t,z) }
	\end{array} \right.  \in L^\infty(\R_+  \times \R) \right\}
	\end{align*}
	and we associate the corresponding norm :
	\begin{equation*}
	\normf{v} = \max \left(  \norm{v}_\E, \sup_{ (t,z) \in\R_+\times \R} \abs{ 2 v(t,\bz(t)) - v(t,z) },\sup_{ (t,z) \in\R_+\times \R} \Big(1+ \abs{z-\zs(t)}\Big)^\alpha  \abs{ \p_z v(t,\bz(t)) - \p_z v(t,z) } \right). 
	\end{equation*}
We will use the notational shortcut $\vphia$ for the weight function :
\begin{align*}
\vphia (t,z):= \Big( 1+ \abs{z-\zs(t)} \Big)^\alpha.
\end{align*}
\end{definition}
Since most of this paper is focused around the pair $(q_\e,V_\e) \in \R \times \F$, we will use the convenient notation $\norm{(q,V)} := \max \Big( \abs{q}, \normf{V} \Big)$. Our main theorem is the following convergence result :
\begin{thm}[Convergence]\label{main theo}$ $\\
	There exist $K_0$, $K'_0$ and $\e_0>0$ such that if we make the following assumptions on the initial condition, for all $ \e \leq \e_0 $ :
	\begin{align*}
	\quad \normf{V_\e(0,\cdot)-\Va(0,\cdot)} & \leq \e^2 K_0, \\
	\abs{q_\e(0)-\qa(0)} & \leq \e^2 K_0,  \text{ and }  \\\quad  \abs{p_\e(0)-\pa(0)} & \leq \e^2 K_0,
	\end{align*}
	then we have uniform estimates of the solutions of the Cauchy problem:
	\begin{align*}
\sup_{t>0}	\  \normf{V_\e-V^\ast}   & \leq \e^2 K'_0, \\
	\sup_{t>0} \abs{q_\e(t)-\qa(t)} &  \leq \e^2 K'_0\\
  	\sup_{t>0} \abs{p_\e(t)-\pa(t)} &  \leq \e^2 K'_0,
	\end{align*}
	where $\qa$ is the solution of \cref{def qda} associated to $\qa(0)$ and $\pa$ is the solution of \cref{def pda} associated to $\pa(0)$. The function $\Va$ is defined  in \cref{def Va}.
\end{thm}
Therefore, as predicted, the limit of $U_\e$ when $\e \to 0$ is the function $ \pa(t) + \qa(t)(z-\zs(t)) +\Va(t,z)$.  \cref{main theo} establish the  stability, with respect to $\e$ and uniformly in time, of Gaussian distributions around the dynamics of the dominant trait driven by a gradient flow differential equation. 

In \cite{spectralsex}  a fixed point argument was used to build solutions of the stationary problem when $\e \ll 1$. However, this method can no longer be applied in this case since the derivative in time breaks the structure that made the stationary problem equivalent to a fixed point mapping. In fact, in the present article, \cref{PtU0int} and \cref{PtUeps} are of different nature due to the fast time relaxation dynamics. This is one of the main difficulty of this work compared to \cite{spectralsex}.  
For this reason we replace the fixed  point argument by a perturbative analysis. We introduce the following corrector terms, $\kappa_\e,W_\e$, our aim is to bound them uniformly :
\begin{align}
\label{perturb V}V_\e(t,z) & = \Va(t,z) +   \e^2 W_\e(t,z),	\\ 
\label{perturb kappa} q_\e(t) &  = \qa(t) + \e^2 \kappa_\e(t).
\end{align} 
The scalar $\qa$ , perturbed by $\e^2 \kappa_\e$, will tune further the affine part of the solution. The function $W_\e$ measures the error made when approximating \cref{PtUeps} by  \cref{PtU0}. 
We choose not to perturb $p_\e$ because we will see in \cref{equation pe} that it can be straightforwardly deduced from the analysis.

 This decomposition highlights a crucial part of our analysis, coming back to the initial \cref{Ptfeps}. Strikingly, the main contribution (in $\e$) to the solution is quadratic, see \cref{hopfcole}, and therefore it does not belong to the space of the corrective term $V_\e$. 
The order of precision  is quite high since we are investigating the error made when approximating $f_\e$ by almost Gaussian distributions : $W_\e$ is of order $\e^2$, while $U_\e$ is of order $1$ in $\e$.   The objective of this article is to show that $\kappa_\e$ and $W_\e$ are uniformly bounded with respect to  time and $\e$. 

\section*{Acknowledgment}
This project has received funding from the European Research Council (ERC) under the European Union’s Horizon 2020 research and innovation programme (grant agreement No 639638), it was also supported by the LABEX MILYON (ANR-10-LABX-0070) of Université de Lyon, within the program "Investissements d'Avenir" (ANR-11-IDEX- 0007) operated by the French National Research Agency (ANR).

\section{Heuristics and method of proof}\label{sec heuristics}
For this section only, we focus on the function $U_\e$ instead of $V_\e$ to get a heuristic argument in favor of  the decomposition  \cref{decomp Ue} and some elements supporting \cref{main theo}. We will denote $R_\e$ the perturbation such that we look for solutions of \cref{PtUeps} under the following form:
\begin{align*}
U_\e(t,z) = U^\ast(t,z) + \e^2 R_\e(t,z).
\end{align*}
The function $U^\ast$, defined in \cref{def Uast} also solves \cref{PtU0}. Plugging this perturbation into  \cref{PtUeps} yields the following perturbed equation for $R_\e$:
\begin{multline*}
M(t,z)-\e^2 \p_t U^\ast(t,z) -\e^4 \p_t R_\e(t,z)  =  I_\e( U^\ast + \e^2 R_\e)(t,z) \times \\ 
\exp \Big( U^\ast(t,z)-2 U^\ast \left(t, \bar{z}(t) \right) + U^\ast(t,\zs(t))\Big) \exp \left( \e^2\Big( R_\e(t,z)-2 R_\e \left(t, \bar{z}(t) \right) +R_\e(t,\zs(t))   \Big) \right)	. 
\end{multline*}
By using \cref{PtU0}, one gets that $R_\e$ solves the following :
\begin{multline*}
M(t,z)-\e^2 \p_t U^\ast(t,z) -\e^4 \p_t R_\e(t,z) = \\
I_\e( U^\ast + \e^2 R_\e)(t,z) M(t,z) \exp \left( \e^2\Big( R_\e(t,z)-2 R_\e \left(t, \bar{z}(t) \right) +R_\e(t,\zs(t))  \Big) \right)	. 
\end{multline*}
To prove the boundedness of $R_\e$, solution to this nonlinear equation, we shall linearize it and show a stability result on the linearized problem (see \cref{stab We}). We explain here the heuristics about the linearization.
We have already said that $I_\e$ is expected to converge to 1. Therefore by linearizing the exponential, a natural linearized equation appears to be :
\begin{align}\label{lin We heur}
\e^2   \p_t \tilde{R_\e}(t,z) = M(t,z)  \Big( -\tilde{R_\e}(t,z) + 2 \tilde{R_\e} \left(t, \bar{z}(t) \right)-\tilde{R_\e}(t,\zs(t))  \Big) , 
\end{align}
For clarity we denote by $T$ the linear operator:
\begin{align*}
T(R)(t,z):= M (t,z) \Big( 2 R(t, \bz(t))-R(t,z)+ R(t,\zs(t)) \Big).
\end{align*}
We know precisely what are the eigen-elements of this linear operator.  The eigenvalue $0$ has multiplicity two, the eigenspace consisting of affine functions.  More generally one can get every eigenvalue by differentiating iteratively the operator and evaluating at $z = \zs$. This corresponds to the following table:
\begin{center}
\begin{table}[H]	\begin{tabular}{|l||c|c|c|c|c|c|c|}
		\hline 
		&&&&&\\[-10pt]
		Eigenvalue  :  & $0$ & $0$ & $-\frac12$ & $-\frac34$ & ...\\ &&&&&\\[-10pt]
		Dual eigenvector :  & $\delta_{\zs(t)}$ & $\delta'_{\zs(t)}$ & $\delta''_{\zs(t)}$ &   $\delta^{(3)}_{\zs(t)}$  & ...\\[-10pt] &&&&&\\
		\hline	
	\end{tabular}\caption{Eigen-elements of $T$.}\label{tab}
	\vspace{5pt}
	\end{table}
\end{center}
This explains why $R_\e$ should be decomposed between affine parts and the rest, and as a consequence, also the solution $U_\e$ we are investigating. The scalars $p_\e$ and $q_\e$ of the decomposition \cref{decomp Ue} correspond to the projection of $U_\e$ upon the eigenspace associated to the (double) eigenvalue $0$. On the other hand the rest is expected to remain uniformly bounded since the corresponding eigenvalues are negative, below $-\frac12$. 

Beyond the heuristics about the stability, this linear analysis also illustrates the discrepancy between $V_\e$ and $q_\e$ in \cref{main theo}. While $V_\e$ is expected to relax to a n explicit bounded value arbitrary quickly as $\e \to 0$ (fast dynamics), this is not true for $\qa$ which solves a differential equation (slow dynamics): 
\begin{align*}
& 	\qda(t)  = -m''(\zs(t)) \qa(t) + \frac{ m^{(3)}(\zs(t))}{2}-2 m''(\zs(t))m'(\zs(t)).
\end{align*}
We can infer that the  second eigenvalue $0$ in \Cref{tab} is in fact inherited from $-\e^2 m''(\zs(t))$ at $\e>0$, which explains that we can read $\qa$  at this order.

The technique we will use in the following sections to bound $W_\e$ in $\F$ will seem more natural in the light of this formal analysis. The first step will be to work around $\zs$, the base point of the dual eigen-elements in \Cref{tab}. We will derive uniform bounds up to  the third derivative to estimate $W_\e$, see \cref{stab We}. 

By plugging the expansions of \cref{perturb V} and \cref{perturb kappa} associated to the decomposition \cref{decomp Ue} and the logarithmic  transform \cref{hopfcole} into our original model \cref{Ptfeps}, we obtain the following main reference equation that we will study in the rest of this article :
\begin{multline}\label{perturb Weps}
M(t,z)-\e^2\Big(  \pde(t) + \qdea(t)(z-\zs) + m'(\zs)\qa(t)   + \p_t \Vea(t,z) \Big) \\ -\e^4 \Big(  \kde(t)(z-\zs) + m'(\zs) \kappa_\e(t) + \p_t W_\e(t,z) \Big)   = M(t,z) \I_\e( \qea+\e^2 \kappa_\e, \Vea + \e^2 W_\e) \times \\  \exp \left( \e^2\Big( W_\e(t,z)-2 W_\e \left(t, \bar{z}(t) \right) + W_\e ( \zs(t) ) \Big) \right)	.
\end{multline}
Our main objective will be to linearize \cref{perturb Weps}, in order to deduce the boundedness of the unknowns, $(\kappa_\e\,W_\e)$, by working on the linear part of the equations. We will need to investigate different scales (in $\e$) to capture the different behaviors of each contribution.

We will pay attention to  the remaining terms. We will use the classical notation $O(1)$ and $O(\e)$, and we will write $ \norm{(\kappa_\e,W_\e)}O(\e)$ to illustrate when the constants of $O(\e)$ depend on $(\kappa_\e,W_\e)$. We also define a refinement of the classical notation $O(\e)$ : 
\begin{definition}[$O^\ast(\e^\alpha)$]\label{def Oast}$ $\\
	For $\alpha \in \N$, we say that a function $g(\e,t,z)$	is such that $g(\e,t,z)=O^\ast(\e^\alpha)$ if there exists $\e^\ast$ such that for all $\e \leq \e^\ast$ it verifies :
	\begin{align*}
	\abs{g(\e,t,z)} \leq C^\ast \e^\alpha, 
	\end{align*}
	and the constants $\e^\ast, C^\ast$ depend only on the pair $(\qa,\Va)$. 
\end{definition}
More generally, when we write $O(\e)$, the constants involved may \textit{a priori} depend upon the pair $(\kappa_\e,W_\e)$. Our intent is to make the dependency of the constants clear when we linearize. This will prove to be a crucial point when we will go back to the non-linear problem \cref{perturb Weps}. We will see that all the terms that do not have a sufficient order in $\e$ to be negligible will be  $O^\ast(1)$, and therefore uniformly bounded independently of $(\kappa_\e,W_\e)$. A key point of our analysis is to segregate those terms when doing the  linearization. 

The rest of the paper is organized as follows :
\begin{itemize}[label=$\triangleright$]
	\item First we prove some properties upon the reference  pair $(\qa,\Va)$ around which all the terms of \cref{perturb Weps} are linearized. 
	\item A key part of our perturbative analysis is to be able to linearize $\I_\e$, which we do in \cref{sec estim Ie} thanks to cautious  estimates upon the directional derivatives.
	\item 
	We  derive an equation on $\kappa_\e$ in \cref{sec kappa}, and later a linear approximated equation for $W_\e$, and more importantly for all of its derivatives in \cref{sec lin}, while controlling precisely the error terms.
	\item We show the boundedness of the solutions of the linear problem in the space $\F$, see \cref{sec lin stab}, mainly through maximum principles and a dyadic division of the space to take into account the non local behavior of the infinitesimal operator. This is the content of \cref{stab We}.
	\item Finally, we tackle the proof of \cref{main theo} in the \cref{sec proof}.
	\item  To conclude, in \cref{sec num} we discuss some of our  assumptions made in \cref{def m cauchy}, illustrated by some numerical simulations. 
\end{itemize}
\newpage
\section*{Index of Notations}
 \begin{tabular}{@{} l l}
$U_\e$	& 
Perturbation of Gaussian distribution to solve the Cauchy problem \cref{Ptfeps}, see \cref{hopfcole}\\
$m$ & Selection function \\
$I_\e$	& Residual shape of the infinitesimal operator after transformation, defined in  \cref{def_Ieps} \\
$\zs$ & Dominant trait in the population, solves a gradient flow ODE: \cref{gradient flow}  \\
$\bar z $	& $\frac{z+\zs}{2}$ \\
$M$	& $1+m(z)-m(\zs(t))-m'(\zs(t))(z-\zs(t))$ \\
$p_\e,q_\e,V_\e$	& $U_\e(t,z) = p_\e(t) + q_\e(t)(z-\zs(t)) + V_\e(t,z)$ \\
$\I_\e$ & Same thing as $I_\e$ but as a function of \emph{two} variables : $\I_\e(q_\e,V_\e) = I_\e(U_\e)$ \\
$U^\ast$	& Limit of  $U_\e$ when $\e \to 0$ \\
$\pa, \qa,\Va$ & $U^\ast(t,z) = \pa(t) + \qa(t)(z-\zs(t)) + \Va(t,z)$ \\
$ (\kappa_\e, W_\e) $ & $q_\e = \kappa_\e + \e^2 \kappa_\e, V_\e = \Va+\e^2 W_\e$ \\
$\Xi_\e(t,z)$ & $W_\e(t,z)-2 W_\e(t,\bz)$\\
	$\E$    &  Functional space to measure $V_\e$ \\
	$\vphia$ & Weight function, $\vphia(t,z) = \Big( 1 + \abs{z-\zs(t)}\Big)^\alpha$ \\
	$\F$    &  Functional space to measure $W_\e$ \\
	$\norm{(g,W)}$ & $\max \Big(\abs{g}, \norm{W}_\F \Big)$, or $\max \Big(\abs{g}, \norm{W}_\E\Big)$ depending on the context \\
		$\East$    &  Functional space to measure $\Va$\\
		$K^\ast$ & Uniform bound of $\Va $ in $\East$
\\
$O^\ast(\e)$   & Special negligible term $O(\e)$ where the constants depend only on $K^\ast$ \\
$\Iea$	& $\I_\e(\qa,\Va)$ \\
$Q(y_1,y_2)$ & $\frac 12 y_1y_2 +\frac 34 (y_1^2+y_2^2)$\\
$Y$ & $(y_1,y_2)$ \\
$\D(V)(Y,t,z)$  & $V(t, \bz) - \frac12 V(t,\bz + \e y_1)- \frac12 V(t,\bz + \e y_2)$ \\
$\D^\ast (V) (y,t)$ &  $V(t, \zs) -  V(t,\zs +\epsilon y)$ \\	$\norm{W}_\infty$   &  $\ds \sup_{ (t,z) \in\R_+ \times \R} \abs{W(t,z)} $ \\
$dG_\e^\ast, dN_\e^\ast$	&  Probability densities that simplify some notations in \cref{sec estim kappa}, defined in \cref{def_dG} and \cref{def_dN}. \\
$O^\ast_0 (1)$  & $ \max\Big( O^\ast(1), O^\ast(1)\norm{  W_\e(0,\cdot)}_\F \Big)$ (where we slightly abuse notation) \\
$B_0$ & Ball that contains $\zs$, see \cref{fig rings} \\
$D_n$   & The nth dyadic ring defined in \cref{def Dn}, see \cref{fig rings} \\
$\norm{W}_\infty^0$, \, $\norm{W}_\infty^n$    &  $\ds \sup_{ (t,z) \in\R_+ \times B_0} \abs{W(t,z)} $ , \quad  $\ds \sup_{ (t,z) \in\R_+ \times D_n} \abs{W(t,z)} $ \\
  & \\

	&\\
\end{tabular}\newpage
\section{Preliminary results : estimates of $\Iea$ and $\Vea$ }\label{sec estim Iea}
\subsection{Control of $(\qa,\Va)$}$ $\\
Before tackling the main difficulties of this article, we first state some controls on the function $\Va$, solution of \cref{PtU0}.  Most of them use the explicit expression of \cref{def Va} and were proved in \cite{spectralsex}. To be able to measure this function we introduce another functional space, with more constraints. 
\begin{definition}[Subspace $\East$]\label{def East}$ $\\
	We define $\East$ as the following subspace of $\E$ :
\begin{align*}
\East := \E  \cap \left\{v \in \mathcal C^5(\R_+ \times \R)    \left|
\begin{array}{rrr} \ds  & \vphia(t,z) \abs{\p_z^4 v(t,z)} \\
\ds  &  \vphia(t,z) \abs{\p_z^5 v(t,z)}
\end{array} \right.  \in L^\infty(\R_+  \times \R) \right\}
\end{align*}
We equip it with the norm  $\nast{\cdot}$ :
	\begin{equation*}
	\nast{v} = \max \left(  \norm{v}_\E, \sup_{ (t,z) \in\R_+\times \R} \vphia(t,z) \abs{\p_z^4 v(t,z)},\sup_{ (t,z) \in\R_+\times \R} \vphia(t,z) \abs{\p_z^5 v(t,z)} \right). 
	\end{equation*}
\end{definition}
Our intention with the successive definitions of the functional spaces is to be able to measure each term of the decomposition made in \cref{perturb V} as follows:
\begin{align*}
\underbrace{V_\e}_{ \E} = \underbrace{V^\ast_{}}_{  \East}  + \e^2  \underbrace{W_\e}_{ \F}.
\end{align*}
The fact that $\Va \in \East$ is part of the claim of the following lemma.
\begin{lem}[Properties of $\Va$]\label{prop Va}$ $\\
	The function $\Va$ belongs to the space $\East$. Moreover,
	\begin{align}\label{link pVea pm}
	\p_z^2 \Vea(t,\zs) =  2 m'' (\zs),  \quad \quad 
	 \p_z^3 \Vea(t,\zs) =\frac43 m^{(3)}(\zs).
	\end{align}
\end{lem}
\begin{proof}$ $\\
	Precise estimates of the summation operator that defines $\Va$ in \cref{def Va} are studied in \cite{spectralsex}. They can be applied there thanks to the decay assumptions about  $M$, \cref{decay Gamma}. The only difference here is that an uniform bound for the fourth and fifth derivative are required. The proofs of those bounds rely solely upon the assumption made in \cref{decay Gamma}, for the fourth and fifth  derivative of $M$.  This shows that $\Va \in \East$.  Explicit computations based on the formula \cref{def Va} prove the relationships \cref{link pVea pm}.
\end{proof}
A consequence of \cref{prop Va} is that since $ m'' (\zs(t))>0$ for $t>t_0$, thanks to \cref{sign pzzm}, which implies that $\Va$ is locally convex around $\zs(t)$. However we need more information upon $\Va$ than the space it belongs to. We will bound $(\qa,\Va)$ independently of time. This is the content of the following result: 
\begin{prop}[Uniform bound on $(\qa,\Vea)$]\label{prop Vea bound}$ $\\
	There exists  a constant $K^\ast$ such that for $j=0,1, 2$ and $3$, we have 
	\begin{align*}
	\max \Big( \nast{\Vea}, \norm{\qa}_{L^\infty(\R_+)}, \norm{\p_t \, \p_z^j \Vea}_{L^\infty(\R_+ \times \R)} \Big) & \leq K^\ast.
	\end{align*}
\end{prop}
\begin{proof}[Proof of \cref{prop Vea bound}]$ $\\
	For the estimates upon $\Va$ and $\p_t \Va$, it is a direct consequence of the definition of $\East$ and the explicit formula \cref{def Va}. The technique to bound the sums is to distinguish between the small and large indices, it was detailed in \cite{spectralsex}. 
	
For $\qa$, one must look to \cref{def qda}. The boundedness of $\qa$ is a straightforward consequence of the convexity of $m$ at $\zs(t)$ for large times, see \cref{sign pzzm} and the convergence of $\zs$ to bound the other terms.
\end{proof}
\subsection{Estimates of $\Iea$ and its derivatives}$ $\\
We next define a notational shortcut for the functional $\I_\e$ introduced in \cref{def_Ieps}, when it is evaluated at the reference pair $(\qa,\Va)$:
\begin{align*}
\Iea := \I_\e(\qa,\Va).
\end{align*}
This section is devoted to get precise estimates of this function. This will be crucial for the  linearization of $\I_\e(\qa + \e^2 \kappa_\e, \Va+\e^2W_\e)$ as can be seen on the full \cref{perturb Weps}. 
\begin{prop}[Estimation  of $\Iea$]\label{estim Iea}$ $\\
	\begin{align*}
	\Iea(t,z)= 1 + O^\ast(\e^2),
	\end{align*}	
	where the constants of $O^\ast(\e^2)$ depend only on $K^\ast$, introduced in \cref{prop Vea bound}, as defined by \cref{def Oast}. 
\end{prop}
The proof  consists in exact Taylor expansion in $\e$. Very similar expansions were performed in \cite[Lemma 3.1]{spectralsex}, we adapt the method of proof here, since it will be used extensively throughout this article.
\begin{proof}[Proof of \cref{estim Iea}]$ $\\
	We recall that by \cref{prop Vea bound}, $\max(\abs{\qa},\nast{\Va}) \leq K^\ast$ , and, by definition 
	\begin{multline*}
	\Iea(t,z)  =  \\ \frac{\ds \iint_{\R^2}  e^{-Q(y_1,y_2)} \exp \Big( -\e\qea(t)(y_1+y_2) +2 \Vea(t, \bz) - \Vea(t, \bz +\e y_1)-\Vea(t,\bz +\e y_2) \Big) dy_1 dy_2 }{\ds \sqrt{\pi} \int_{\R} e^{ -\frac12 \abs{y}^2} \exp \Big(-\e\qea(t) y  +  \Vea(t,\zs)-   \Vea(t,  \zs + \e y) \Big) dy} \\
	:= \frac{N(t,z)}{D(t)},
	\end{multline*}
where $Q$ the quadratic form appearing after the rescaling of the infinitesimal operator in \cref{def_Ieps}: 
\begin{align*}
Q(y_1,y_2):= \frac12 y_1y_2 + \frac34(y_1^2 + y_2^2).
\end{align*}
This quadratic form will appear very frequently in what follows, mostly, as here, through the  bi-variate Gaussian distribution it defines. Once and for all we state that a correct normalization of this Gaussian distribution is:
\begin{align*}
\frac 1{\sqrt{2} \pi} \iint_{ \R^2 } e^{-Q(y_1,y_2)} dy_1dy_2 = 1.
\end{align*}
We start the estimates with the more complicated term, the numerator $N$. With an exact Taylor expansion inside the exponential, there exists generic $0 <\xi_i<1$, for $j=1,2$, such that
	\begin{multline*}
	N(t,z) = \ds \dfrac{1}{\sqrt 2 \pi} \iint_{\R^2} e^{-Q(y_1,y_2)} \exp \Big[ -\e\qea(t)(y_1+y_2) -\e(y_1+y_2) \p_z \Va(t, \bz)\\  - \frac{\e^2}{2}\Big( y_1^2 \p_z^2 \Va(t, \bz +\e \xi_1y_1)+ y_2^2 \p_z^2 \Va(t,\bz +\e \xi_2 y_2) \Big) \Big] dy_1 dy_2 
	\end{multline*}
	Moreover we can write for some $\theta = \theta(y_1,y_2)\in (0,1)$,
	\begin{align*}
	& \exp(-\e P) = 1 -\e P + \frac{\e^2 P^2}{2} \exp(-\theta\e  P)\, , \\ 
	& \label{eq weight2} P = (y_1+y_2) \Big(\qa(t) + \p_z \Va(t,\bz) \Big) + \dfrac{\e}{2}\Big( y_1^2 \p_z^2 \Va(t,\bz + \e \xi_1 y_1) + y_2^2 \p_z^2 \Va(t,\bz+\e \xi_2 y_2) \Big).
	\end{align*}
	such that 
	\begin{equation}\label{eq boundP}
	\quad |P|\leq K^\ast \left( \abs{y_1}+\abs{y_2} + \frac{\e (y_1^2 + y_2^2)}{2} \right),
	\end{equation}
	Combining the expansions, we find:
	\begin{align}
	\nonumber N(t,z) & =  \frac{1}{\sqrt{2 }\pi }\iint_{\R^2} e^{-Q(y_1,y_2)} \Big(1 -\e P + \frac{\e^2 P^2}{2} \exp(-\theta\e  P)\Big) dy_1dy_2, \\
	&\label{eq P2} = 1 -\e \frac{1}{\sqrt{2} \pi }\iint_{\R^2} e^{-Q(y_1,y_2)}  Pdy_1dy_2 + \frac{\e^2}{2 \sqrt{2 }\pi }\iint_{\R^2} e^{-Q(y_1,y_2)} P^2 \exp(-\theta\e  P) dy_1dy_2
	\end{align}
	The key part is the cancellation of the terms $O(\e)$ due to the symmetry of $Q$ :
	\begin{equation*}
	\frac{1}{\sqrt{2 }\pi }\iint_{\R^2} e^{-Q(y_1,y_2)}  (y_1+y_2) dy_1dy_2 =0.
	\end{equation*}
	Therefore :
	\begin{multline*}
	\frac{\e}{\sqrt{2 } \pi }\iint_{\R^2} e^{-Q(y_1,y_2)}  Pdy_1dy_2 =   \frac{\e^2 }{2 \sqrt 2 \pi }\iint_{\R^2} e^{-Q(y_1,y_2)}  \Big( y_1^2 \p_z^2 \Va(t,\bz + \e  \xi_1 y_1) + y_2^2 \p_z^2 \Va(t,\bz)+\e \xi_2 y_2)\Big) dy_1dy_2.
	\end{multline*}
	And we get the estimate
	\begin{align*}
	\abs{\frac{\e}{\sqrt{\pi}}\iint_{\R^2} e^{-Q(y_1,y_2)}  Pdy_1dy_2 } & \leq  \frac{\e^2 }{2 \sqrt 2 \pi }\iint_{\R^2} e^{-Q(y_1,y_2)}  \Big( y_1^2 + y_2^2\Big) K^\ast dy_1dy_2 \\
	& \leq O^\ast(\e^2).
	\end{align*}
	Thanks to \cref{eq boundP} it is easy to verify that the last term of \cref{eq P2} behaves similarly:
	\begin{align*}
	\frac{\e^2}{2 \sqrt{2 } \pi}\iint_{\R^2} e^{-Q(y_1,y_2)} P^2 \exp(-\theta\e  P) dy_1dy_2 = O^\ast(\e^2).
	\end{align*}
	Indeed, it  states that the term $P$ is at  most quadratic with respect to $y_i$ so $Q + \theta\e P$ is uniformly bounded below by a positive quadratic form for $\e$ small enough.   
	This shows that 
	\begin{align*}
	N(t,z)= 1 + O^\ast(\e^2).
	\end{align*}
	The denominator is easier, with the same arguments, using the Gaussian density :
	\begin{align*}
	D(t) = 1 + O^\ast(\e^2).
	\end{align*}
	Combining the estimates of $N$ and $D$, we get the desired result.
\end{proof}
There exists a link between $\qa$ and $\p_z \Iea(t,\zs)$, which is in fact the motivation behind the choice of $\qa$. 
\begin{prop}[Link between $\qa$ and $\p_z \Iea(t,\zs)$]\label{link qa pzIea}$ $\\
	\begin{align*}
	\p_z \Iea(t,\zs(t)) = \e^2 \left( m''(\zs(t)) \qa(t) - \frac{m^{(3)}(\zs(t))}{2}\right) + O^\ast(\e^4),
	\end{align*}
	where the constants of $O^\ast(\e^4)$ only depend on $K^\ast$.
\end{prop}
The proof of this result was the content of \cite[Lemma 3.1]{spectralsex} and only requires that $(\qa,\Va)$ is uniformly bounded, as stated in \cref{prop Vea bound}. Its proof follows the same procedure of exact Taylor expansions as in the one of \cref{estim Iea}.  

It will be useful to dispose of estimates of $\p_z \Iea$ not only at the point $\zs$. They are less precise, as stated in the following proposition:
\begin{prop}[Estimates of the decay of the derivatives of $\Iea$]\label{decay pzIea}$ $\\
	There exists a constant $\e_\ast$ that depends only on $K^\ast$ such that for all $\e \leq \e\ast$, for $j=1,2,3$:
	\begin{align*}
	\sup_{(t,z) \in\R_+ \times \R} \vphia(t,z) \abs{\p_z^{(j)} \Iea(t,z)}   = O^\ast (\e^2).
	\end{align*}	
\end{prop}
To shortcut notations, we introduce the following difference operator that appears  in the integral $\I_\e$ see \cref{def_Ieps} :
\begin{align}\label{not D}
\D(V)(Y,t,z) & := V(t, \bz) - \frac12 V(t,\bz + \e y_1)- \frac12 V(t,\bz + \e y_2), \quad Y =(y_1,y_2),\\
\D^\ast (V) (y,t) & := V(t, \zs) -  V(t,\zs +\epsilon y).
\end{align}
We will use the following technical lemma giving an estimate of the weight function against the derivatives of a given function.
\begin{lem}[Influence of the weight function.]\label{decay norm}$ $\\	There exists a constant $C$ such that for each ball $B$ of $\East$ or $\F$,  there exists $\e_B$ such that for every $W \in B$, for every $y \in \R $ and $\e \leq \e_B$, for $j=2,3,4$ or $5$:
	
	\begin{align*}
\vphia(t,z) \abs{\p_z^{(j)} W(t, \bz(t) + \e y)} 	& \leq C   \norm{W} & \text{ if} \abs{y}\leq \abs{z-\zs(t)}, \\
	& \leq(1+\abs{y}^\alpha ) \norm{W}. & \text{ otherwise,}
	\end{align*}
	with $\norm{W}= \nast{W}$ or $\normf{W}$ depending on the case. 
\end{lem}
The \cref{decay pzIea} is a prototypical result. It will be followed by a series of similar statements. Therefore, we propose two different proofs. In the first one, we write exact Taylor expansions. However the formalism is heavy, which is why we propose next a formal argument, where the Taylor expansions are written without exact rests for the sake of clarity. 

In the rest of this paper  more complicated estimates will be proved, in the spirit of \cref{decay pzIea},  see \cref{control pgVIe,control pzgV2Ieps} for instance. The notations and formulas will be very long, so we shall only write the "formal" parts of the argument. However it can all be made rigorous, as  below.
\begin{proof}[Proof of \cref{decay pzIea}]$ $\\
	First, write the expression for the derivative, using our notation $\D$ introduced in \cref{not D} : 	
	\begin{align}\label{explicit pzIea}
	\p_z \Iea(t,z) & = \frac{\ds \iint_{\R^2} e^{-Q(y_1,y_2)} \exp \Big( -\e\qea(y_1+y_2) +2 \D(\Va)(Y,t,z) \Big) \D( \p_z \Vea)(Y,t,z)   dy_1 dy_2 }{\ds \sqrt{\pi}\int_{\R} e^{-\frac12 \abs{y}^2} \exp \Big(-\e\qea y  +  \D^\ast(\Va)(y,t) \Big) dy} ,\\
	\nonumber & := \frac{N(t,z)}{D(t)}.
	\end{align}
	We only focus on the numerator. The denominator $D$ can be handled similarly as in the proof of \cref{estim Iea}, where we show that it is essentially $1 + O^\ast(\e^2)$. We perform two Taylor expansions in the numerator $N$, namely: 
	\begin{equation}\label{Taylor1}
	\begin{cases}
	2 \D(\Va)(Y,t,\bz)   = -\e(y_1+y_2)\p_z \Va(t,\bz)- \dfrac{\e^2}{2}\Big( y_1^2 \p_z^2 \Va(t,\bz + \e \xi_1 y_1) + y_2^2 \p_z^2 \Va(t,\bz + \e \xi_2 y_2) \Big),\medskip\\
	\ds \D(\p_z \Va)(Y,t,\bz)  = - \dfrac{\e (y_1+y_2)}{2} \p_z^2 \Va(t,\bz)- \dfrac{\e^2}{4} \Big(y_1^2 \p_z^3 \Va(t,\bz + \e \xi_1 y_1) +   y_2^2 \p_z^3 \Va(t,\bz + \e \xi_2 y_2)\Big),
	\end{cases}
	\end{equation}
	where $\xi_i$ denote some generic number such that $0<\xi<1$ for $i=1,2$. 
	Moreover, we can write
	\begin{multline}\label{P def}
	\exp(-\e P) = 1 -\e P \exp(-\theta\e P)\, \quad \text{ with } \\
	P := (y_1+y_2)\Big( \p_z \Va(t,\bz) + \qa \Big) +  \dfrac{1}{2}\Big( y_1^2 \p_z^2 \Va(t,\bz + \e \xi_1 y_1) + y_2^2 \p_z^2 \Va(t,\bz + \e \xi_2 y_2) \Big)\,
	\end{multline}
	for some $\theta = \theta(y_1,y_2)\in (0,1)$. Combining the expansions, we find:
	\begin{multline*}
\vphia(t,z) \p_z \Iea(t,z) = \dfrac{\vphia(t,z)}{ \sqrt{2} \pi } \ds  \iint_{\R^2}e^{-Q(y_1,y_2)} \left ( 1  -\e P \exp(-\theta\e P)  \right ) \\  \times\left ( - \dfrac{\e (y_1+y_2)}{2} \p_z^2 \Va(t,\bz)- \dfrac{\e^2}{4} \Big(y_1^2 \p_z^3 \Va(t,\bz + \e \xi_1 y_1) +   y_2^2 \p_z^3 \Va(t,\bz + \e \xi_2 y_2)\Big)  \right )  d y_1 d y_2 \, .
	\end{multline*}
	The crucial point is the cancellation of the $O(\e)$ contribution due to the symmetry of $Q$, as already observed above:
	\begin{equation*}
	\iint_{\R^2}e^{-Q(y_1,y_2)} (y_1 + y_2) d y_1 d y_2 = 0\, .
	\end{equation*}
	So, it remains 
	\begin{multline*}
\vphia(t,z)N(t,z)  = \\
	- \e^2 \dfrac{\vphia(t,z)}{4 \sqrt{2} \pi } \iint _{\R^2}e^{-Q(y_1,y_2)} \left[ y_1^2 \p_z^3 \Va(t,\bz + \e \xi_1 y_1) +   y_2^2 \p_z^3 \Va(t,\bz + \e \xi_2 y_2))  \right] dy_1dy_2 \\
	+ \e^2  \frac{\vphia(t,z)}{2 \sqrt{2 }\pi } \iint _{\R^2} e^{-Q(y_1,y_2)} P \exp(-\theta \e P)  (y_1+y_2)  \p_z^2 \Va(t,\bz)    dy_1dy_2  \\
+ 	\e^3\frac{\vphia(t,z)}{4\sqrt{2} \pi} \iint _{\R^2}    e^{-Q(y_1,y_2)} P \exp(-\theta \e P) \left ( y_1^2 \p_z^3 \Va\!(t,\bz + \e \xi_1 y_1)   +  y_2^2 \p_z^3 \Va\!(t,\bz  + \e \xi_2 y_2) \right )         dy_1dy_2.
	\end{multline*}
	If we forget about the weight in front of each term, clearly the last two contributions are uniform $O^\ast(\e)$ since $\e \leq \e_\ast$ small enough and $\Va$ and $\qa$ are uniformly bounded by $K^\ast$, see \cref{prop Vea bound}. The term $P$ is at most quadratic with respect to $y_i$, see \cref{P def}, so $Q + \theta\e P$ is uniformly bounded below by a positive quadratic form for $\e$ small enough.    
	
	A difficulty is to add the weight to those estimates. To do so, we use \cref{decay norm}, for each integral term appearing in the previous formula, because each time appears a term of the form :
	\begin{align*}
	\vphia(t,z) \p_z^{(j)} \Va(t, \bz + \e \xi_i y_i).
	\end{align*}
	Since  every $\xi_i$ verifies $0< \xi_i<1$, the bounds given by \cref{decay norm} ensure that each integral remains bounded by moments of the bivariate Gaussian defined by $Q$, as if there were no weight function. This concludes the proof of the first estimate \cref{decay pzIea}. 
	
	Bounding the quantity $ \vphia(t,z)   \abs{\p_z^{(j)} \Iea (t,z)}$ for $j=2,3$  follows the same steps, as it can be seen on the explicit formulas :
	\begin{multline}\label{pzzIea}
	\p_z^2 \Iea(t,z) = \\ \frac{\ds \iint_{\R^2} \!\!\!  \exp \Big(\! -Q(y_1,y_2) -\e g (y_1+y_2) + 2 \D (\Va)(Y,t,z)\Big) \!\! \left[  \D(\p_z \Va)^2  + \frac12  \D(\p_z^2 \Va) \right](Y,t,z)  dy_1 y_2 } {\ds \sqrt{\pi} \int_{\R}e^{-\frac12\abs{y}^2} \exp \Big(-\e \qa  y +  \D^\ast(\Va)(y,t) \Big) dy},
	\end{multline}
	\begin{multline}\label{pzzzIea}
	\p_z^3 \Iea(t,z) = \ds \iint_{\R^2} \exp \Big(-Q(y_1,y_2) -\e g (y_1+y_2) + 2 \D (\Va)(Y,t,z)\Big) \times  \\  \frac{ \left[  \D(\p_z \Va)^3  + \dfrac{3}{2}  \D(\p_z \Va) \D(\p_z^2 \Va) + \dfrac{1}{4} \D(\p_z^3 \Va) \right](Y,t,z)   } {\ds \sqrt{\pi} \int_{\R}e^{-\frac12\abs{y}^2} \exp \Big(-\e \qa y  + \D^\ast(\Va)(y,t) \Big) dy}dy_1 y_2.
	\end{multline}
The motivation behind going up to the order $5$ of differentiation for $\Va$ in the \cref{def East} lies in the terms 
\begin{align*}
\frac12  \D(\p_z^2 \Va)  \text{ and } \dfrac{1}{4} \D(\p_z^3 \Va) .
\end{align*}
To gain an order $\e^2$ as needed in \cref{decay pzIea} for the estimates, one needs to go up by two orders in the Taylor expansions, which involve fourth and fifth order derivatives. The importance of the order $\e^2$ will later appear in \cref{prop lin Ieps} and the \cref{sec lin stab}.
\end{proof}
We now  propose a formal argument, much simpler to read. 
\begin{proof}[Formal proof of \cref{decay pzIea}.]$ $\\
We tackle the first derivative.	We use the same notations as previously, see \cref{explicit pzIea}, and again focus on the numerator $N$. Formally,
	\begin{multline*}
	N(t,z)=  \iint_{\R^2} e^{-Q(y_1,y_2)} \exp \Big[ -\e  (y_1+y_2) \Big( \qa + \p_z \Va(t,\bz) \Big)+  (y_1^2+y_2^2)O^\ast(\e^2) \Big] \\
	\times \Big[ -\e(y_1+y_2)\p_z^2 \Va(t,\bz) + (y_1^2+y_2^2)O^\ast(\e^2) \Big]   dy_1 dy_2.
	\end{multline*}
Thanks to the linear approximation of  the exponential, we find:
	\begin{multline}
	\label{eq decayIea}N(t,z) =  \iint_{\R^2} e^{-Q(y_1,y_2)} \Big[1-\e  (y_1+y_2) \Big( \qa + \p_z \Va(t,\bz) \Big)  + (y_1^2+y_2^2)O^\ast(\e^2) \Big] \times \\   \Big[ -\e(y_1+y_2)\p_z^2 \Va(t,\bz) + (y_1^2+y_2^2)O^\ast(\e^2) \Big]   dy_1 dy_2.
	\end{multline}
By sorting out the orders in $\e$, this can be rewritten:
	\begin{align*}
	N(t,z) = \e N_1+ O^\ast(\e^2).
	\end{align*}
	By symmetry :
	\begin{align*}
	N_1 : = -\iint_{\R^2} e^{-Q(y_1,y_2)} \Big[ \e(y_1+y_2)\p_z^2 \Va(t,\bz)  \Big] dy_1 dy_2 = 0.
	\end{align*}
	To conclude, we notice that we can add the weight function to those estimates and make the same arguments as in the previous proof.  
\end{proof}
\begin{proof}[Proof of \cref{decay norm}]$ $\\
	If $\abs{z-\zs} \leq 1$, then $ 1 + \abs{z-\zs}\leq 2$ and the result is immediate by the \cref{def East,def F} of the adequate functional spaces. Therefore, one can suppose that $\abs{z-\zs}>1$. We first look at the regime  $\abs{y} \leq \abs{z-\zs}$. 
	Then, by definition of the norms,  
	\begin{align}
	\nonumber \vphia(t,z) \abs{\p_z^{(j)} W(t, \bz + \e y)}   \leq &  	 \, 2 \dfrac{\abs{z-\zs}^\alpha}{ \abs{\bz + \e y - \zs}^\alpha} \Big(  \abs{\bz + \e y - \zs}  \Big)^\alpha \abs{\p_z^{(j)} W(t, \bz + \e y)} \\
	& \label{eq intdecay} \leq 2 \dfrac{\abs{z-\zs}^\alpha}{\abs{\bz + \e y - \zs}^\alpha} \norm{W}. 
	\end{align}
	To bound the last quotient, we use the following inequality, that holds true because we are in the regime $\abs{y} \leq \abs{z-\zs}$ :
	\begin{align*}
	\abs{\bz+ \e y-\zs} \geq -\abs{\e y} + \abs{\bz-\zs} \geq \frac12\abs{z-\zs}- \e \abs{z-\zs}.
	\end{align*}
	This yields 
	\begin{align}\label{eq intdecay2}
	2 \dfrac{\abs{z-\zs}}{\abs{\bz + \e y - \zs}} \leq \frac{2}{\frac12 - \e}.
	\end{align}
	Bridging together \cref{eq intdecay,eq intdecay2}, one gets the \cref{decay norm} in the regime $\abs{y}\leq \abs{z-\zs}$; on the condition that $\e < \frac12$. 
	
	On the contrary, when $\abs{z-\zs} \leq \abs{y}$, we have immediately:
	\begin{align*}
	\Big(1+\abs{z-\zs}^\alpha \Big) \abs{\p_z^{(j)} W(t, \bz + \e y)}   \leq  (1+\abs{y}^\alpha ) \norm{W}.
	\end{align*}
\end{proof}

\section{Linearization of $\I_\e$ and its derivatives}\label{sec estim Ie}
The first step to obtain a linearized equation on $W_\e$ is to study the nonlinear terms of \cref{perturb Weps}. A key point is the study of the functional $\I_\e$ defined in \cref{def_Ieps}, which plays a major role in our study. We will show that it converges uniformly to $1$, as we claimed in the \cref{sec intro} and that its derivatives are uniformly small, with some decay for large $z$, similarly to what we proved for the function $\Iea$ in the previous section. This will enable us to linearize $\I_\e$ and its derivatives in \cref{prop lin Ieps,prop lindecay}.
\subsection{Linearization of $\I_\e$}$ $\\
We first bound uniformly all the terms that appear during the linearization of $\I_\e$ by Taylor expansions. One starts by measuring the first order directional derivatives. 
\begin{prop}[Bounds on the directional derivatives of $\I_\e$]\label{control pgVIe}$ $\\
	For any ball $B$ of $\R \times \E$, there exists a constant $\e_B$ that depends only on $B$ such that for all $\e\leq \e_B$ we have for all $(g,W) \in B$, and $H \in \E$ :
	\begin{align}
	\label{bound pgIe} \sup_{(t,z)\in\R_+\times \R} \abs{ \p_g  \I_\e(g,W)(t,z)} & \leq  \norm{(g,W)}O(\e^2)  , \\
	\label{bound pVIe} \sup_{(t,z)\in\R_+\times \R} \abs{ \p_V  \I_\e(g,W) \cdot H(t,z)} & \leq \norm{(g,W)} \norm{H}_\E O(\e^2)  .
	\end{align}
\end{prop}

\begin{proof}[Proof of \cref{control pgVIe}]$ $\\
	As in the estimates of $\Iea$ and its derivatives in the previous section, the argument to obtain the result will be to perform exact Taylor expansions. As explained before we will not pay attention to  the exact rests that can be handled exactly as before, and we refer to the proof of \cref{estim Iea,decay pzIea} to see the details. However   our computations  will make clear  the order $\e^2$ of \cref{bound pgIe,bound pVIe}. First, thanks to the derivation with respect to $g$ an order of $\e$ is gained straightforwardly:
	\begin{multline}\label{eq pgIe1}
	\p_g \I_\e(g,W)(t,z)  = \\
	- \e  \left( \frac{\ds \iint_{\R^2} \exp\Big[ -Q(y_1,y_2)  -\e g(y_1+y_2) +2\D(W)(Y,t,z) \Big] (y_1 +y_2) dy_1 dy_2 }{\ds \sqrt{\pi} \int_{\R} e^{-\frac12\abs{y}^2} \exp \Big(-\e g  y  + \D^\ast(W)(y,t) \Big) dy} \right.\\
	\left. - \I_\e(g,W)(t,z)  \frac{\ds \int_{\R} e^{-\frac12\abs{y}^2} \exp \Big(-\e g  y  + \D^\ast(W)(y,t) \Big) y \, dy}{ \sqrt{\pi} \ds \int_{\R} e^{-\frac12\abs{y}^2} \exp \Big(-\e g  y  +  \D^\ast(W)(y,t) \Big) dy} \right) .
	\end{multline}
	The common denominator is bounded : 
	\begin{align*}
	\ds \int_{\R} e^{-\frac12\abs{y}^2} \exp \Big(-\e g  y  +  \D^\ast(W)(y,t) \Big) dy \geq \ds\int_{\R} \exp \left[ - \frac12 |y|^2  -2 \e \abs{y} \norm{(g,W)}\right] d y.
	\end{align*}
	For the numerators, a supplementary order in $\e$ is gained by symmetry of $Q$, as in other estimates, see \cref{decay pzIea} for instance. For the single integral we write:
	\begin{multline*}
	\ds \int_{\R}  e^{-\frac12\abs{y}^2} y \exp \Big(-\e g y  + \D^\ast(W)(y,t) \Big) dy \leq  \ds \int_{\R}  e^{-\frac12\abs{y}^2} y \exp \Big(-\e g y  + 2 \e \abs{y} \norm{(g,W)} \Big) dy \\ 
	\leq  \ds \int_{\R}  e^{-\frac12\abs{y}^2} y \Big[1- \e g y  + O(\e) \abs{y} \norm{(g,W)} \Big] dy.  
	\end{multline*}
	Finally
	\begin{align}
	\ds \int_{\R}  e^{-\frac12\abs{y}^2} y \exp \Big(-\e g y  + \D^\ast(W)(y,t) \Big) dy \leq \norm{(g,W)}  O (\e)  .\label{eq pgIe2}
	\end{align}
	For the first numerator of \cref{eq pgIe1}, computations work in the same way :
	\begin{align}
\nonumber	\ds \iint_{\R^2} \exp\Big( -Q(y_1,y_2)  -\e &g(y_1+y_2)  +2 \D(W)(Y,t,z)  \Big) (y_1 +y_2) dy_1 dy_2  \\
\nonumber& 	\leq \ds \iint_{\R^2} \exp\Big( -Q(y_1,y_2) + O(\e)  (y_1+y_2)\norm{(g,W)} \Big) (y_1+y_2) dy_1dy_2, 
 \\ 
 \nonumber&  \leq  \ds \iint_{\R^2}\exp \Big( -Q(y_1,y_2)\Big) \Big[1 + O(\e)  (y_1+y_2)\norm{(g,W)} \Big] (y_1+y_2) dy_1dy_2 \\
 &	\leq  \norm{(g,W)} O(\eps). \label{eq pgIe3}
		\end{align}
	Therefore, combining \cref{eq pgIe1,eq pgIe2,eq pgIe3} we have proven the bound upon the first derivative of $\I_\e$ in \cref{bound pgIe}. 
	
	Concerning \cref{bound pVIe}, one starts by writing the following formula for the Fréchet derivative:
	\begin{multline}\label{eq pVIe1}
	\p_V \I_\e(g,W)\cdot H(t,z)  = \\  \frac{\ds \iint_{\R^2} \exp\Big[ -Q(y_1,y_2)  -\e g(y_1+y_2) +2 \D(W)(Y,t,z) \Big] 2 \D(H)(Y,t,z) dy_1 dy_2 }{\ds \sqrt{\pi} \int_{\R}  e^{-\frac12\abs{y}^2} \exp \Big(-\e g  y  +  \D^\ast (W)(y,t)   \Big) dy} \\
	- \I_\e(g,W)(t,z)  \frac{\ds \int_{\R}  e^{-\frac12\abs{y}^2} \exp \Big(-\e g  y  +  \D^\ast (W)(y,t) \Big) \D^\ast (H)(y,t) dy}{ \sqrt{\pi} \ds \int_{\R}  e^{-\frac12\abs{y}^2}  \exp \Big(-\e g  y  + \D^\ast (W)(y,t) \Big) dy}  .
	\end{multline}
	The claimed order $\e^2$ holds true, by similar symmetry arguments. For instance, when we do the Taylor expansions on the numerator of the first term of \cref{eq pVIe1}: 
	\begin{multline*}
	\ds \iint_{\R^2}\exp\Big[ -Q(y_1,y_2)  -\e g(y_1+y_2) +2 \D(W)(Y,t,z) \Big] 2 \D(H)(Y,t,z) dy_1 dy_2  \\
	= 2 \ds \iint_{\R^2} \exp\Big( -Q(y_1,y_2) \Big) \Big[1 -\e  (y_1+y_2) \Big( g + \p_z W(t,\bz) \Big)   + O(\e^2) \norm{W}_\E \Big] \\ \times  \Big[-\e(y_1+y_2) \p_z H (t, \bz ) + O(\e^2)(y_1^2+y_2^2) \norm{H}_\E  \Big] dy_1dy_2,
	\end{multline*}
	\begin{multline} \label{eq pVIe2}
	\quad  \quad = -2 \e \p_z H (t, \bz ) \ds \iint_{\R^2} \exp\Big( -Q(y_1,y_2)\Big)(y_1+y_2) dy_1dy_2    +   \e^2 \p_z H (t, \bz)  \Big( g + \p_z W(t,\bz) \Big) \times \\ \left( \iint_{\R^2} \exp\Big( -Q(y_1,y_2)\Big)  (y_1+y_2)^2  dy_1dy_2 \right) + O(\e^2) \norm{H}_\E \norm{(g,W)} \\ \leq  \norm{(g,W)} \norm{H}_\E  O(\e^2). 
	\end{multline}
	For the second term of \cref{eq pVIe1}, we also gain an order $\e^2$ when making Taylor expansions of $\D^\ast(W)$,  since $\p_z H(t,\zs)=0$:
	\begin{align}
	\nonumber \ds \int_{\R}  e^{-\frac12\abs{y}^2} & \exp \Big(-\e g y  + \D^\ast(W)(y,t) \Big)\D^\ast(H)(y,t) dy \\
	\nonumber & = - \ds \int_{\R}  e^{-\frac12\abs{y}^2} \exp \Big(-\e g y   + 2 \e \abs{y} \norm{(g,W)} \Big) y^2  O(\e^2) \norm{H}_\E dy,  \\
	\nonumber & =  -\ds \int_{\R}   e^{-\frac12\abs{y}^2}  \Big[1-\e g y  +  2 \e \abs{y} \norm{(g,W)} \Big] y^2  O(\e^2) \norm{H}_\E dy  \\ & \leq  \norm{(g,W)}  \norm{H}_\E  O(\e^2).  \label{eq pVIe3}
	\end{align}
	As before the denominator of \cref{eq pVIe1} has a universal lower bound, therefore combining  \cref{eq pVIe1,eq pVIe2,eq pVIe3}  concludes the proof. 
\end{proof}
We have proven all the tools to linearize $\I_\e$ as follows, thanks to the previous estimates on the directional derivatives of $\I_\e$.
\begin{prop}[Linearization of $\I_\e$]\label{prop lin Ieps}$ $\\
	For any ball $B$ of $\R \times \E$, there exists a constant $\e_B$ that depends only on $B$ such that for all $\e\leq \e_B$ we have for all $(g,W) \in B$ : 
	\begin{align}
	\I_\e(\qea + \e^2 g, \Vea + \e^2 W)(t,z) & = \Iea(t,z)+ O(\e^3)\norm{(g,W)} ,\label{Iea O3} \\ 
	&  =1+ O^\ast(\e^2) + O(\e^3)\norm{(g,W)},\label{Iea 1Oast}
	\end{align}
	where $O(\e^3)$ only depends on the ball $B$.
\end{prop}

\begin{proof}[Proof of proposition \cref{prop lin Ieps}]$ $\\
	We  write an exact Taylor expansion:
	\begin{multline*}
	\I_\e(\qea + \e^2 g, \Vea + \e^2 W) =  \I_\e^\ast + \\\e^2 \Big[ \p_g   \I_\e(\qea + \e^2 \xi g, \Vea + \e^2 \xi W) + \p_V  \I_\e(\qea + \e^2 \xi g, \Vea + \e^2 \xi W) \cdot W \Big] .
	\end{multline*}
	for some $0 < \xi <1$. Therefore \cref{Iea O3} is a direct application  of \cref{control pgVIe} to $g'=\qa+ \e^2 \xi g$, $W'=\Va+  \e^2 \xi W$ and $H=W$. One deduces the estimation of \cref{Iea 1Oast} from \cref{estim Iea}.
\end{proof}
As a matter of fact, in \cref{Iea 1Oast}, we have even shown an estimate $1+ O^\ast(\e^2) + O(\e^4)\norm{(g,W)}$. However we choose to reduce arbitrarily the order in $\e$ for consistency reasons with further  estimates of this article. It suffices to our purposes.
\subsection{Linearization of $\p_z \I_\e$ and decay estimates}$ $\\
To prove \cref{main theo},  we need to bound uniformly $\norm{W_\e}_\F$, and this implies $L^\infty$ bounds of the derivatives of $W_\e$. To obtain those, our method is to work on the linearized equations they verify.  Therefore, linearizing  $\I_\e$ is not enough, we need to linearize $\p_z^{(j)} \I_\e$ as well, for $j=1,2$ and $3$.  For that purpose we need more details than previously upon the nature of the negligible terms. More precisely we need to know how it behaves relatively to the weight function of the space $\E$ and $\F$, that acts by definition upon the second and third derivatives. The objective of this section is to linearize $\p_z^{(j)} \I_\e$ to obtain similar results to \cref{prop lin Ieps}.  We first prove the following estimates on the derivatives of $\I_\e$:
\begin{prop}[Decay estimate of $\p_z \I_\e$]\label{decay pzIe}$ $\\
	For any  ball $B$ of $\R \times \E$,  there exists a constant $\e_B$ that depends only on $B$ such that for any  pair $(g,  W)$ in $B$ , for all $\e \leq \e_B$ :
	\begin{align*}
	\sup_{(t,z) \in \R_+ \times \R} \vphia(t,z)\abs{\p_z \I_\e(g,W)(t,z)}  & \leq  \norm{(g,W)}O(\e), \\
	\sup_{(t,z) \in \R_+ \times \R} \vphia(t,z)  \abs{\p_z^2 \I_\e(g,W)(t,z)}   & \leq \norm{(g,W)}O(\e), \\
	\sup_{(t,z) \in \R_+ \times \R} \vphia(t,z)  \abs{\p_z^3 \I_\e(g,W)(t,z)}  & \leq  \norm{(g,W)}O(\e^\alpha) + \frac1{2^{1-\alpha}} \norm{\vphia \p_z^3 W}_\infty.
	\end{align*}	
	where all $O(\e)$ depend only on the ball $B$, and $\ds \norm{\vphia \p_z^3 W}_\infty = \sup_{ (t,z) \in\R_+ \times \R} \vphia(t,z) \abs{\p_z^3 W_\e(t,z)}$.
\end{prop}
This proposition has to be put in parallel with \cite[Proposition 4.6]{spectralsex}. We are not able to propagate an order $\e$ for all derivatives, but the factor $\frac1{2^{1-\alpha}} \norm{\vphia \p_z^3 W}_\infty$ that we pay can, and will, be involved in a contraction argument, just as in \cite{spectralsex}, mostly since $2^{\alpha-1}<k(\alpha)<1$, where $k(\alpha)$ plays the same role in \cref{stab We}. This is the core of the perturbative analysis strategy we use.  
\begin{proof}[Proof of \cref{decay pzIe}]$ $\\
	We focus on the first derivative, the proof for the second one can be straightforwardly adapted. 
	\begin{multline}\label{eq pzIe1}
	\vphia(t,z) \p_z \I_\e(g,W)(t,z)  = 
\vphia(t,z) \times \\  \frac{\ds \iint_{\R^2} \exp\Big[ -Q(y_1,y_2)  -\e g(y_1+y_2) + 2 \D(W)(Y,t,z) \Big]  \D\left(\p_z W \right)(Y,t,z)  dy_1 dy_2 }{\ds \sqrt{\pi} \int_{\R}  e^{-\frac12\abs{y}^2} \exp \Big(-\e g  y  +  \D^\ast(W)(y,t) \Big) dy}.
	\end{multline}
As before the following formal Taylor expansions hold true for the numerator, ignoring the weight in a first step:
	\begin{multline} \label{eq pzIe2}
	\ds \iint_{\R^2}\exp\Big[ -Q(y_1,y_2)  -\e g(y_1+y_2) +2 \D(W)(Y,t,z) \Big] \D(\p_z W)(Y,t,z) dy_1 dy_2  \\
	= \ds \iint_{\R^2} \exp\Big[ -Q(y_1,y_2) \Big]\Big[1 - O(\e)  (y_1+y_2) \norm{(g,W)} \Big]  \Big[-O(\e)(y_1+y_2) \norm{(g,W)} \Big] dy_1dy_2,\\
	\leq   O(\e)\norm{(g,W)}. 
	\end{multline}
	Meanwhile the denominator has a uniform lower bound :
	\begin{align*}
	\ds \int_{\R} e^{-\frac12\abs{y}^2} \exp \Big(-\e g  y  +   \D^\ast(W)(y,t) \Big) dy \geq \ds\int_{\R} \exp \left[ - \frac12 |y|^2  -2 \e \abs{y} \norm{(g,W)} \right] d y.
	\end{align*}
	The estimate of \cref{eq pzIe2} can be made rigorous as in the proof of \cref{decay pzIea} for instance. Moreover, one can add the weight to bound \cref{eq pzIe1} thanks to \cref{decay norm}, as explained in the proof of \cref{decay pzIea}.  Therefore, the proof of the first estimate of \cref{decay pzIe} is achieved.
	
	For the second term of \cref{decay pzIe}, involving the second order derivative, 
	the arguments and decomposition of the space are the same, we follow the same steps, with the formula
	\begin{multline*}
	\p_z^2 \I_\e(g,W)(t,z) = \\ \frac{\ds \iint_{\R^2} \exp \Big(-Q(y_1,y_2) -\e g (y_1+y_2) + 2 \D (W)(Y,t,z)\Big) \Big[  \D(\p_z W)^2  + \frac12  \D(\p_z^2 W) \Big](Y,t,z)  dy_1 y_2 } {\ds \sqrt{\pi} \int_{\R}e^{-\frac12\abs{y}^2} \exp \Big(-\e g  y  +  \D^ \ast(W)(y,t) \Big) dy}.
	\end{multline*}
	Things are a little bit different for the third derivative, as can be seen on the following explicit formula:
	\begin{multline}\label{pzzzIe}
	\p_z^3 \I_\e(t,z) = \frac{\ds \iint_{\R^2} \exp \Big(-Q(y_1,y_2) -\e g (y_1+y_2) + 2 \D (W)(Y,t,z)\Big)} {\ds \sqrt{\pi} \int_{\R}e^{-\frac12\abs{y}^2} \exp \Big(-\e g  y  + \D^\ast(W)(y,t) \Big) dy} \times  \\   \left[  \D(\p_z W)^3  + \dfrac{3}{2}  \D(\p_z W) \D(\p_z^2 W) + \dfrac{1}{4} \D(\p_z^3 W) \right](Y,t,z) dy_1 y_2.  
	\end{multline}
	All terms in this formula will provide an order $\e$ exactly as before, except for the linear contribution $\D(\p_z^3 W)$ since we lack a priori controls of the fourth derivative of $W$ in $\F$. Therefore, for this term we proceed as follows : 
	\begin{align}
\nonumber	\vphia(t,z) \abs{ \D(\p_z^3 W)(Y,t,z) } &  =   \Big( 1 + \abs{z-\zs}  \Big)^\alpha   \abs{ \p_z^3 W(t,\bz) - \frac12 \p_z^3 W(t,\bz + \e y_1 )-\frac12 \p_z^3 W(t,\bz + \e y_2 )} 
	 \\ & \nonumber \leq  \Big( 1 + \abs{z-\zs}  \Big)^\alpha \left( \abs{ \p_z^3 W(t,\bz)} +  \frac12\abs{ \p_z^3 W(t,\bz + \e y_1 )} + \frac12 \abs{\p_z^3 W(t,\bz + \e y_2 )} \right), \\
	&  \label{eq D}\leq 2^{\alpha+1} \norm{\vphia \p_z^3 W}_\infty  \Big( 1 +  {\e^\alpha}  (\abs{y_1}^\alpha + \abs{y_2}^\alpha ) \Big).
	\end{align}
	For this computation, we used the following property of the weight function, which was also of crucial importance in \cite[Lemma 4.5]{spectralsex}:
	\begin{align*}
\sup_{ (t,z) \in\R_+ \times \R}	\frac{\vphia(t,z)}{\vphia(t,\bz)} \leq 2^\alpha.
	\end{align*}
	As a consequence, take $i=1$ or $2$, then:
	\begin{align*}
\vphia(t,z) \abs{  \p^3_z W(\bz+\e y_i)} &  \leq  \dfrac{ 2^\alpha \vphia(t,\bz)  } { \Big(1 + \abs{\bz+\e y_i-\zs} \Big)^\alpha  }   \norm{\vphia \p_z^3 W}_\infty, \\
	& \leq 2^\alpha \left( 1 + \dfrac{|\e y_i |}{1+|\bz+\e y_i-\zs|} \right)^\alpha  \norm{\vphia \p_z^3 W}_\infty \leq 2^\alpha( 1 + \e^\alpha |y_i|^\alpha  ) \norm{\vphia \p_z^3 W}_\infty.
	\end{align*}
	As a consequence, we deduce that 
	\begin{multline*}
	\vphia(t,z) \frac{\ds \iint_{\R^2} \exp \Big(-Q(y_1,y_2) -\e g (y_1+y_2) + 2 \D (W)(Y,t,z)\Big)   \left[ \frac14 \D(\p_z^3 W)(Y,t,z) \right] dy_1 y_2}{\ds \sqrt{\pi} \int_{\R}e^{-\frac12\abs{y}^2} \exp \Big(-\e g  y  + \D^\ast(W)(y,t) \Big) dy} \\
	\leq   \frac{1} {2^{1-\alpha}} \norm{ \vphia \p_z^3 W}_\infty + O(\e^\alpha ) \norm{(g,W)},
	\end{multline*}
by sub-additivity of $\abs{\cdot}^\alpha$. This justifies \cref{eq D}. Once added to other estimates of the terms of \cref{pzzzIe}, obtained by Taylor expansions of $\D$ as before we get the desired estimate.
\end{proof}
One can notice in the proof that the order $O(\e)$ is not the sharpest one can possibly get for the first derivative, see \cref{eq pzIe2}. However it is sufficient for our purposes. 
We now detail the control upon the directional derivatives of $\I_\e$.
\begin{prop}[Bound of the directional derivatives of $\I_\e$]\label{decay pzgVIe}$ $\\
	For any  ball $B$ of $\R \times \E$,  there exists a constant $\e_B$ that depends only on $B$ such that for any  pair $(g,  W)$ in $B$ and any function $H \in \E$, for every $\e \leq \e_B$ :
	\begin{align}
		\label{bound pgzIe} \sup_{(t,z) \in\R_+ \times \R} \left( \vphia(t,z)  \abs{\p_g \p_z^{(j)} \I_\e(g,W)(t,z)}  \right) & \leq  O(\e)\norm{(g,W)}_\E , & (j=1,2,3),\\
		\sup_{(t,z) \in\R_+ \times \R} \left(\vphia(t,z)   \abs{\p_V \p_z^{(j)} \I_\e(g,W) \cdot H (t,z)}   \right) & \leq  O(\e)\norm{H}_\E, & (j=1,2),\\
	 \label{bound pgzzzIe} \sup_{(t,z) \in\R_+ \times \R}  \Big( \vphia(t,z) \abs{\p_V \p_z^3 \I_\e(g,W) \cdot H (t,z)} \Big)  & \leq   O(\e^\alpha) \norm{H}_\E+ \frac{1}{2^{1-\alpha}} \norm{\vphia \p_z^3 H}_{\infty}.
	\end{align}
where the $O(\e)$ depends only on the ball $B$.
\end{prop}
As for \cref{decay pzIe}, in those estimates, the order of precision $O(\e)$ is not optimal and we could improve it without it being useful. We will not give the full proof for each estimate of this Proposition. However, we see that it follows the same pattern than in \cref{decay pzIe}, and we will even use those results for the proof. In particular for the third derivative, it is not possible to completely recover an order $\e$, hence the term $\frac{1}{2^{1-\alpha}} \norm{\vphia \p_z^3 H}_{\infty}$. It comes from the \emph{linear} part $\D(\p_z^{3}W)$ that appears in $\p_z^3 \I_\e$, see \cref{pzzzIe}. However, it does not prevent us from carrying our analysis since the factor $\frac{1}{2^{1-\alpha}}$ will be absorbed by a contraction argument, see \cref{sec proof}.
\begin{proof}[Proof of \cref{decay pzgVIe}]$ $\\
	We first detail the proof of \cref{bound pgzIe}, because derivatives in $g$ are somehow easier to estimate.
	The formula for the first derivative is:
	\begin{multline}\label{eq pgzIe1}
	\p_g \p_z \I_\e(g,W)(t,z)  = \\
	- \e \left( \frac{\ds \iint_{\R^2} \exp\Big[ -Q(y_1,y_2)  -\e g(y_1+y_2) +2 \D(W)(Y,t,z)\Big] (y_1 +y_2) \D(\p_z W)(Y,t,z) dy_1 dy_2 }{\ds \sqrt{\pi} \int_{\R} e^{-\frac12 \abs{y}^2} \exp \Big(-\e g  y  + \D^\ast(W)(y,t) \Big) dy} \right.\\
	\left. - \p_z \I_\e(g,W)(t,z)  \frac{\ds \int_{\R} e^{-\frac12 \abs{y}^2} \exp \Big(-\e g  y  +  \D^ \ast(W)(y,t) \Big) y dy}{ \sqrt{\pi} \ds \int_{\R} e^{-\frac12 \abs{y}^2} \exp \Big(-\e g  y  +  \D^\ast(W)(y,t) \Big) dy} \right) .
	\end{multline}
	The first term of this formula closely resembles the one for $\p_z I_\e(g,W)$, with an additional factor $\e(y_1+y_2)$. We do not detail how to bound it, as it follows the same steps, see the work done following \cref{eq pzIe1}.  For the second term we first use the following  bound : 
	\begin{align}\label{eq fracpgzIe}
	\frac{\ds \int_{\R} e^{-\frac12 \abs{y}^2} \exp \Big(-\e g  y  +  \D^ \ast(W)(y,t) \Big) y dy}{ \sqrt{\pi} \ds \int_{\R} e^{-\frac12 \abs{y}^2} \exp \Big(-\e g  y  +  \D^\ast(W)(y,t) \Big) dy}  \leq \frac{\ds\int_{\R} \exp \left[ - \frac12 |y|^2  + 2 \e \abs{y}  \norm{(g,W)}\right] y  d y}{\ds\int_{\R} \exp \left[ - \frac12 |y|^2  -2 \e \abs{y}  \norm{(g,W)}\right] d y}.
	\end{align}
	For $\e$ sufficiently small that depends only on $\norm{(g,W)}$ we deduce an uniform bound with moments of the Gaussian distribution. We then use the estimate from \cref{decay pzIe} on $\p_z I_\e(g,W)$, which takes the weight into account, to conclude. 
	
	Every other estimate of \cref{decay pzgVIe} works along the same lines. We illustrate this with the second derivative in $g$ and $z$:
	\begin{multline}\label{eq pgzzIe}
	\p_g \p_z^2 \I_\e(g,W)(t,z) =  - \e \left( \frac{\ds \iint_{\R^2} \exp \Big(-Q(y_1,y_2) -\e g (y_1+y_2) + 2 \D (W)(Y,t,z)\Big)}{\ds \sqrt{\pi} \int_{\R}e^{-\frac12\abs{y}^2} \exp \Big(-\e g  y  + \D^\ast(W)(y,t) \Big) dy}  \right.  \\ \times\left.(y_1 + y_2 ) \left[  \D(\p_z W)^2  + \frac12 \D(\p_z^2 W) \right](Y,t,z)  dy_1 y_2   \right.  \\ \left. 
	-  \p_z^2 \I_\e(g,W) \dfrac{\ds \int_{\R}e^{-\frac12\abs{y}^2} y  \exp \Big(-\e g  y  +  \D^\ast(W)(y,t) \Big)  \,  dy}{\ds \int_{\R}e^{-\frac12\abs{y}^2} \exp \Big(-\e g  y  +  \D^\ast(W)(y,t) \Big) dy} \right).
	\end{multline}
	This is very close to $\p_z^2 I_\e$ that has already been estimated in \cref{decay pzIe}, and therefore the same arguments as before hold.
	
	The structure is different for the derivatives in $V$, as can be seen  for $\p_V \p_z \I_\e(g,W)\cdot H  $  :
	\begin{multline}
	\p_V \p_z \I_\e(g,W)\cdot H  (t,z) = \\
	\frac{\ds \iint_{\R^2} \! \!\! \exp\Big(\! \! - \! Q(y_1,y_2) \!   - \! \e g(y_1+y_2) \!+\!2 \D(W)(Y,t,z)\Big) \! \Big[ \! 2 \D(\p_z W)\D(H) \! + \! \D(\p_z H) \Big]\!(Y,t,z) dy_1 dy_2 }{\ds \sqrt{\pi} \int_{\R} e^{-\frac12 \abs{y}^2} \exp \Big(-\e g  y  +  \D^\ast(W)(y,t) \Big) dy} \\
	- \p_z \I_\e(g,W)(t,z)  \frac{\ds \int_{\R} e^{-\frac12 \abs{y}^2} \exp \Big(-\e g  y  +  \D^\ast(W)(y,t) \Big) \D^\ast(H)(y,t)  dy}{ \sqrt{\pi} \ds \int_{\R} e^{-\frac12 \abs{y}^2} \exp \Big(-\e g  y  +  \D^\ast(W)(y,t) \Big) dy} .
	\end{multline}
	The second term can still be bounded using \cref{decay pzIe} and the estimate \cref{eq fracpgzIe}, and the following immediate result :
	\begin{equation*}
	\ds \abs{\D^\ast(W)(y,t)}  \leq \e \abs{y} \norm{W}_\E.
	\end{equation*}
	For the first term, we must do Taylor expansions of $2 \D(\p_z W)\D(H) + \D(\p_z H)$ to control them with the weight. One ends up with moments of the multidimensional Gaussian distribution just as in all the previous proofs. For instance, 
	\begin{align*}
2	\vphia(t,z ) \abs{\D(\p_z W)\D(H)}(t,z)&  \leq  \vphia(t,z) \abs{\D(\p_z W)(t,z)} O(\e) (\abs{y_1}+\abs{y_2}) \norm{H}_\E, \\
& \leq  O(\e) (\abs{y_1}+\abs{y_2}+\abs{y_1}^{\alpha+1}+\abs{y_2}^{1+\alpha})(\abs{y_1}+\abs{y_2}) \norm{H}_\E \norm{W}_\E.
	\end{align*}
The same method holds for the second derivative in $V$ and $z$.
	
The estimate of the third derivative in $g$ and $z$ is similar to the previous computations with the following formula: 
		\begin{multline}\label{pgzzzIe}
\p_g 	\p_z^3 \I_\e (t,z) = - \e \left( \frac{\ds \iint_{\R^2} \exp \Big(-Q(y_1,y_2) -\e g (y_1+y_2) + 2 \D (W)(Y,t,z)\Big)} {\ds \sqrt{\pi} \int_{\R}e^{-\frac12\abs{y}^2} \exp \Big(-\e g  y  + \D^\ast(W)(y,t) \Big) dy} \times \right.  \\   (y_1+y_2) \left[  \D(\p_z W)^3  + \dfrac{3}{2}  \D(\p_z W) \D(\p_z^2 W) + \dfrac{1}{4} \D(\p_z^3 W) \right](Y,t,z) dy_1 y_2     \\\left. + \p_z^3 \I_\e(t,z) \dfrac{\ds \int_{\R}e^{-\frac12\abs{y}^2} y  \exp \Big(-\e g  y  +  \D^\ast(W)(y,t) \Big)  \,  dy}{\ds \int_{\R}e^{-\frac12\abs{y}^2} \exp \Big(-\e g  y  +  \D^\ast(W)(y,t) \Big) dy}\right) .
	\end{multline}
However, to get the bound \cref{bound pgzzzIe}, things are a little bit different, because of the linear term of higher order  $\D(\p_z^3 H)$: 
	\begin{multline*}
	\p_V \p_z^3 \I_\e(g,W)\cdot H  (t,z) = \frac{ \ds \iint_{\R^2} \exp \Big(-Q(y_1,y_2) -\e g (y_1+y_2) + 2 \D (W)(Y,t,z)\Big)} {\ds \sqrt{\pi} \int_{\R}e^{-\frac12\abs{y}^2} \exp \Big(-\e g  y  + \D^\ast(W)(y,t) \Big) dy} \times  \\ \left[ \D(H) \left( 2 \D(\p_z W)^3  + 3  \D(\p_z W) \D(\p_z^2 W) + \dfrac{1}{2} \D(\p_z^3 W) \right)   + \right. \\
	\left.\left(  3 \D(\p_z H) \D(\p_z W)^2 \! + \! \frac32  \Big( \! \D(\p_z W) \D(\p_z^2 H) \!+ \! \D(\p_z H) \D(\p_z^2 W) \Big) \! + \! \dfrac{1}{4} \D(\p_z^3 H) \right) \right]   (Y,t,z) dy_1 y_2 \\
	+ \p_z^3 \I_\e(t,z) \frac{\ds \int_{\R} e^{-\frac12 \abs{y}^2} \exp \Big(-\e g  y  +  \D^\ast(W)(y,t) \Big) \D^\ast(H)(y,t)  dy}{ \sqrt{\pi} \ds \int_{\R} e^{-\frac12 \abs{y}^2} \exp \Big(-\e g  y  +  \D^\ast(W)(y,t) \Big) dy}.
	\end{multline*}
	We do not get an order $\e$ from the linear part $\D(\p_z^3 H)$, since we do not control the fourth derivative in $\E$. We then proceed with arguments following \cref{pzzzIe} in the proof of \cref{decay pzIe}. 
\end{proof}
Thanks to those estimates we are able to write our main result for this part, which is a precise control of the linearization of the derivatives of $\I_\e$ :
\begin{prop}[Linearization with weight]\label{prop lindecay}$ $\\
	For any ball $B$ of $\R \times \E$, there exists a constant $\e_B$ that depends only on $B$ such that for all $\e\leq \e_B$ we have for all $(g,W) \in B$ : 
	\begin{align}
	\label{lindecay pzIe}	 & \p_z \I_\e(\qea + \e^2 g, \Vea + \e^2 W)(t,z)   =    \p_z \Iea(t,z) + \frac{\norm{(g,W)}}{\vphia(t,z)}O(\e^3), \\
	\label{lindecay pzzIe} & \p_z^2 \I_\e(\qea + \e^2 g, \Vea + \e^2 W)(t,z)  =   \p_z^2 \Iea(t,z) + \frac{ \norm{(g,W)}}{\vphia(t,z)}O(\e^3), \\
	\label{lindecay pzzzIe} & \p_z^3 \I_\e(\qea + \e^2 g, \Vea + \e^2 W)(t,z)  =   \p_z^3 \Iea(t,z) + \frac{\e^2\norm{\vphia \p_z^3 W}_\infty}{2^{1-\alpha}\vphia(t,z)} +\frac{ \norm{(g,W)}}{\vphia(t,z)}O(\e^{2+\alpha})  .
	\end{align}
	where $O(\e^3)$ only depends on the ball $B$.
\end{prop}
\begin{proof}[Proof of \cref{prop lindecay}]$ $\\
	The methodology for \cref{lindecay pzIe,lindecay pzzIe,lindecay pzzzIe} is the same. We detail for instance how to prove \cref{lindecay pzIe}.  One begins by writing the following exact Taylor expansion up to the second order:
	\begin{multline*}
	\p_z \I_\e(\qea + \e^2 g, \Vea + \e^2 W)(t,z) = \\
	\p_z \I_\e^\ast(t,z) + \e^2 \Big[ \p_g  \p_z \I_\e(\qea + \e^2 \xi g, \Vea + \e^2 \xi W)(t,z) + \p_V \p_z \I_\e(\qea + \e^2 \xi g, \Vea + \e^2 \xi W) \cdot W(t,z) \Big] .
	\end{multline*}
	with $0 < \xi <1$. The result for \cref{lindecay pzIe} is then given by the directional decay estimates  of \cref{decay pzgVIe} applied to $g'=\qa+ \e^2 \xi g$, $W'=\Va+ \e^2 \xi W $, $H = W$.
\end{proof}
Together with \cref{decay pzIea}, we know exactly how  $\p_z^{j} I_\e$ behaves when $\e$ is small:
\begin{align*}
\p_z^{(j)} \I_\e(\qea + \e^2 g, \Vea + \e^2 W)(t,z) = O^\ast(\e^2) + \frac{\norm{(g,W)}}{\vphia(t,z)}O(\e^3),  
\end{align*}
where  $j=1,2$, and only slightly different for $j=3$.
\subsection{Refined estimates of $\Iea$ at $z=\zs$}\label{sec estim kappa}$ $\\
To conclude this section dedicated to estimates of $\I_\e$, we now show that  our estimates above can be made much more precise when looking at the particular case of the function  $\Iea$ evaluated at the point $\zs$. In particular we will gain information upon the sign of the derivatives, that will prove crucial regarding the stability of $\kappa_\e$. This additional precision is similar to what was  needed in the stationary case, \cite[Lemma 3.1]{spectralsex}, where detailed expansions of $\I_\e$ were needed for the study of the affine part, thereby named $\gamma_\e$. We will find convenient to use the following notations, as in \cite{spectralsex}:
\begin{definition}[Measures notation]$ $\\
	We introduce  the following measures  : \begin{multline}\label{def_dG}
	d G^\ast_\e (Y,z,t) := \dfrac{G_\e^\ast (Y,t,z)}{\ds \iint_{\R^2} G_\e^\ast(Y, t,z)dy_1dy_2 } , \quad \text{ with } Y= (y_1,y_2), \\
	=  \dfrac{\ds \exp \Big[-Q(y_1,y_2)- \e \qea (y_1+y_2) +2 \D(\Va)(Y,t,z) \Big] }                               {\ds\iint_{\R^2} \!\! \! \exp \! \Big[\!- \! Q(y_1,y_2)- \e \qea (y_1+y_2) +2  \D(\Va)(Y,t,z)  \Big]  d y_1 d y_2 }. 
	\end{multline}
	And :
	\begin{align}\label{def_dN}
	d N^\ast_\e (y,t) := \dfrac{N_\e^\ast (y,t)}{\ds \int_{\R^2} N_\e^\ast(\cdot,t) } : =  \dfrac{\ds \exp \left[ - \frac12 |y|^2 -\e \qea y + \D^\ast(\Va)(y,t) \right] }                               {\ds\int_{\R} \exp \left[ - \frac12 |y|^2  -\e \qea y +  \D^\ast(\Va)(y,t)\right] d y }. 
	\end{align}
\end{definition}

\begin{prop}[Uniform control of the directional derivatives of $\p_z \I_\e^\ast$ ]\label{control pzgVIea}$ $\\
	There exist a function of time $R_\e^\ast$, such that for any ball $B$ of $\E$,	there exists a constant $\e^\ast$ that depends only on $K^\ast$, that verifies for all $\e\leq \e^\ast$, for all $H \in \E$:
	\begin{align}
	\label{bound pgzIea} &\p_g \p_z \I_\e^\ast (t,\zs ) =  \e^2 R_\e^\ast(t) +O^\ast(\e^3) ,\\
	& \label{bound pzVIea}\p_V \p_z \I_\e^\ast \cdot H(t,\zs) = 
	 O^\ast(\e^2)\norm{H}_\E.
	\end{align}
	where all $O^\ast(\e^j)$ depends only on $K^\ast$ defined in \cref{prop Vea bound} and $R_\e^\ast$ is given by the following formula: 
\begin{align*}
R_\e^\ast(t) :=	m''(t,\zs) \iint_{\R^2}   dG_\e^\ast(Y,t,\zs)  (y_1+y_2)^2 dy_1dy_2,
\end{align*}	
Finally, $R_\e^\ast$ is uniformly bounded and there exists a constant $R_0$ and a time $t_0$ such that $R_\e^\ast\geq  R_0>0$ for all $t\geq t_0$.
\end{prop}
The sign of $R_\e^\ast$ is directly connected to the behavior of $\zs$ we assumed in the introduction, see \cref{sign pzzm}. The derivative in $V$ admits a lower order in $\e$ as in previous estimates, see \cref{lindecay pzzzIe} and \cref{bound pgzzzIe} for instance. This lower order term will be absorbed by a contraction argument, see \cref{sec proof}, once we have a definitive estimate of $\norm{W_\e}_\F$, see the estimate \cref{eq cont}.  
\begin{proof}[Proof of proposition \cref{control pzgVIea}]$ $\\
	First we focus on the bound of \cref{bound pgzIea}. Similarly to \cref{eq pgzIe1}, the explicit formula for the derivative is: 
	\begin{multline}\label{eq pgzIea1}
	\p_g \p_z \Iea(t,\zs) : = -\e ( I_1 + I_2) =\\
	- \e \! \left( \! \frac{\ds \iint_{\R^2} \!\! \exp\Big[\! \! -\!Q(y_1,y_2)  -\e \qa(y_1+y_2) +2 \D(\Va)(Y,t,\zs)\Big]\! (y_1 +y_2) \D(\p_z \Va)(Y,t,\zs) dy_1 dy_2 }{\ds \sqrt{\pi} \int_{\R} e^{-\frac12\abs{y}^2} \exp \Big(-\e \qa  y  +   \D^\ast(\Va)(y,t) \Big) dy} \right.\\
	\left. - \p_z \Iea(t,\zs)  \frac{\ds \int_{\R} e^{-\frac12\abs{y}^2} y \exp \Big(-\e \qa  y  + \D^\ast(\Va)(y,t)  \Big) dy}{ \sqrt{\pi} \ds \int_{\R} e^{-\frac12\abs{y}^2} \exp \Big(-\e \qa  y  +  \D^\ast(\Va)(y,t) \Big) dy} \right) .
		\end{multline}
	Thanks to the \cref{decay pzgVIe}, we already know that $\abs{\p_z \Iea(t,\zs)} = O^\ast(\e^2)$. Moreover we bound uniformly the second term:
	\begin{align*}
	\abs{\frac{\ds \int_{\R} e^{-\frac12\abs{y}^2} y  \exp \Big(-\e \qa  y  + \D^\ast(\Va)(y,t) \Big)  dy}{ \sqrt{\pi} \ds \int_{\R} e^{-\frac12\abs{y}^2} \exp \Big(-\e \qa  y  +  \D^\ast(\Va)(y,t) \Big) dy}} &  \leq   \frac{\ds \int_{\R}\! \exp \Big(\!-\!\frac12\abs{y}^2 +   2 \e K^\ast \abs{y} \! \Big) \abs{y} dy}{ \sqrt{\pi} \ds \int_{\R}\!\! \exp \Big(-\frac12\abs{y}^2 - 2 \e K^\ast \abs{y}  \!\Big) dy} \\ & \leq O^\ast(1), \end{align*}
	where $K^\ast$ was defined in \cref{prop Vea bound}. This shows that $I_2=O^\ast(\e^2)$. Therefore one can focus on $I_1$. In order to gather information upon the sign of this quantity and not only get a bound in absolute value, we perform exact Taylor expansions of $\D (\p_z \Vea)$. We divide $I_1$ by $\Iea(t,\zs)$, and thanks to the definitions of \cref{def_dG,def_dN} we get :
	\begin{align*}
	\frac{ I_1 }{\I_\e^\ast(t,\zs)} = \ds \iint_{ \R^2 } dG_\e^\ast(Y,t,\zs) (y_1+y_2) \D(\p_z \Va)(Y,t,\zs) dy_1dy_2.	\end{align*}
	As usual, we make Taylor expansions : there exists $0 < \xi_1,\xi_2<1$ such that 
	\begin{multline}\label{eq pgzIeps int}
	\frac{I_1}{\I_\e^\ast(t,\zs)}   = \ds \iint_{ \R^2 } dG_\e^\ast(Y,t,\zs)  \left[  -\e \dfrac{(y_1+y_2)^2}{2}\p_z^2 \Vea(t,\zs)      -\frac{\e^2 y_1^2 (y_1+y_2)}{4} \p_z^3 \Vea(t,\zs + \e\xi_1 y_1) \right.   \\\left. -\frac{\e^2 y_2^2 (y_1+y_2)}{4} \p_z^3 \Vea(t,\zs + \e\xi_2 y_2) \right] dy_1dy_2.
	\end{multline}
We next define  $R_\e^\ast$ as :
	\begin{align*}
	\e \p_z^2 \Vea(t,\zs)  \iint_{\R^2}   dG_\e^\ast(Y,t,\zs)\dfrac{(y_1+y_2)^2}{2}  dy_1dy_2 =: \e R_\e^\ast(t),
	\end{align*}
	with the following uniform bounds, that come from bounding by  moments of a Gaussian distribution:  
	\begin{align*}
0 <R_0 \leq R_\e^\ast (t) \quad 	\forall t \geq t_0.
	\end{align*}
	Moreover, it is easy to see that $R_\e^\ast$ is uniformly bounded. 
	The next terms of \cref{eq pgzIeps int} are of order superior to $\e^2$ , and can be bounded uniformly by : 
	\begin{align*}
	\frac{\e^2}{4}  \ds \abs{ \iint_{\R^2} dG_\e^\ast(Y,t,\zs) \Big[   y_1^2 (y_1+y_2) +  y_2^2 (y_1+y_2) \Big]K^\ast dy_1dy_2} \leq  O^\ast(\e^2).
	\end{align*}
	Therefore one can rewrite \cref{eq pgzIeps int} as
	\begin{align*}
	\frac{I_1}{\I_\e^\ast(t,\zs)} = -\e  R_\e^\ast(t) + O^\ast(\e^2).
	\end{align*}
	Thanks to \cref{estim Iea} we recover a similar estimate for $I_1$ :
	\begin{align*}
	I_1 =-\e R_\e^\ast(t) + O^\ast(\e^2).
	\end{align*}
	Finally coming back to \cref{eq pgzIea1}, we have shown that
	\begin{align*}
	\p_g \p_z \Iea(t,\zs) = \e^2  R_\e^\ast(t) + O^\ast(\e^3).
	\end{align*}
	This concludes the proof of the estimate \cref{bound pgzIea}. Next, we tackle the proof of the estimate upon the Fréchet derivative \cref{bound pzVIea}, where, again, we  first divide by $\Iea(t,\zs)$ : 
	\begin{multline}\label{eq pzIea}
	\frac{ \p_V \p_z \I_\e^\ast \cdot H(t,\zs)}{\I_\e^\ast(t,\zs)} = \ds \iint_{\R^2} dG_\e^\ast(Y,t,\zs) \Big[ \D(\p_z \Vea)  2 \D(H) + \D(\p_z H)  \Big] (Y,t,\zs)dy_1dy_2   \\ 
	-\frac{ \p_z \Iea(t,\zs)}{\Iea(t,\zs)} \ds \int_\R  d N_\e^\ast(y,t)   \D^\ast(H)(y,t)  dy.
	\end{multline}
	Thanks to \cref{decay pzIea,estim Iea}, and a uniform bound on $\D^\ast(W) $ :
	\begin{align}\label{eq pzIea2}
	\abs{ \frac{ \p_z \Iea(t,\zs)}{\Iea(t,\zs)} \ds \int_\R  d N_\e^\ast(y,t)   \D^\ast(H)(y,t)  dy.} \leq O^\ast (\e^3)\norm{H}_\E.
	\end{align}
	For the first term of \cref{eq pzIea}, we first make a bound based on Taylor expansions of $\D(H) $ :
	\begin{align*}
	\abs{ \D( H)(Y,t,\zs) } \leq \frac{\e^2}{2} (\abs{y_1}^2 + \abs{y_2}^2)\norm{H}_\E.
	\end{align*}
	The key element here is that since $\D$ is evaluated at $\zs$ one gains an order in $\e$ because $\p_z H(t,\zs)=0$, by definition of $\E$. 
	Therefore, one gets 
	\begin{align}\label{eq pzIea3}
	\abs{ \ds \iint_{\R^2} dG_\e^\ast(Y,t,\zs) \Big[ \D(\p_z \Vea)  2 \D(H) \Big] (Y,t,\zs)dy_1dy_2 }\leq  O^\ast(\e^3) \norm{H}_\E,
	\end{align}
where the additional order in $\e$ is gained through a Taylor expansion of $\D(\p_z \Va)$.
We finally tackle the last term of \cref{eq pzIea} we did not yet estimate, involving $\D(\p_z H)$. Based only on Taylor expansions in $\E$, we do not gain an order $\e^3$ as in the previous terms, which explains our estimate of order $\e^2$ in \cref{eq pzIea}. Rather, we obtain, for some $0<\xi<1$ :
\begin{multline} \label{eq pzIea4}
 \ds \iint_{\R^2} dG_\e^\ast(Y,t,\zs)  \D(\p_z H)  (Y,t,\zs)dy_1dy_2 = \e \frac{\p_z^2 H(t,\zs)}{2} \ds \iint_{\R^2} dG_\e^\ast(Y,t,\zs)  (y_1+y_2)  dy_1dy_2 \\
 + \frac{\e^2 }{4} \ds \iint_{\R^2} dG_\e^\ast(Y,t,\zs) \Big[ y_1^2\p_z^3 H(t,\zs+ \e\xi y_1) +y_2^2\p_z^3 H(t,\zs+ \e\xi y_2) \Big] dy_1dy_2
\end{multline}
It is straightforward, based on multiple similar computations,  to deduce that the first moment of $dG_\e^\ast$ is zero at the leading order. Therefore,
\begin{align}\label{eq pzIea5}
\e \frac{\p_z^2 H(t,\zs)}{2} \ds \iint_{\R^2} dG_\e^\ast(Y,t,\zs)  (y_1+y_2)  dy_1dy_2 =  \e \frac{\p_z^2 H(t,\zs)}{2}O^\ast(\e) = O^\ast(\e^2) \norm{H}_\E.
\end{align}
See for instance the proof of \cref{estim Iea} for similar computations. In the second term of \cref{eq pzIea4},  we also cannot do better than an order $\e^2$.
\begin{multline*}
\frac{\e^2 }{4} \ds \iint_{\R^2} dG_\e^\ast(Y,t,\zs) \Big[ y_1^2\p_z^3 H(t,\zs+ \e\xi y_1) +y_2^2\p_z^3 H(t,\zs+ \e\xi y_2) \Big] dy_1dy_2 \\ \leq  \frac{\e^2 \norm{H}_\E}4 \ds \iint_{\R^2} dG_\e^\ast(Y,t,\zs) \Big[ y_1^2 +y_2^2  \Big] dy_1dy_2 =  O^\ast(\e^2) \norm{H}_\E.
\end{multline*}
Finally, by putting together \cref{eq pzIea2}, \cref{eq pzIea3}, \cref{eq pzIea4} and finally \cref{eq pzIea5}, the estimate \cref{bound pzVIea} is proven.
\end{proof}
The order $\e^3$ of \cref{bound pzVIea} will be crucial in our analysis  around $\kappa_\e$ the perturbation of the linear part $q_\e$ defined in \cref{perturb kappa}. Next, we provide an accurate linearization of $\p_z \I_\e$ compared to the one provided before in \cref{prop lindecay,lindecay pzIe}. This is possible thanks to an evaluation at $z=\zs$, and it will prove useful when tackling the perturbation of the linear part $\kappa_\e$. 
This is the content of the following lemma.
\begin{lem}[Uniform control of the second Fréchet derivative of $\p_z \I_\e $]\label{control pzgV2Ieps}$ $\\
	For any ball $B$ of $\R \times \E$, there exists a constant $\e_B$ that depends only on $B$ such that for all $\e\leq \e_B$ we have for all $(g,W) \in B$ : 
	\begin{multline}\label{lin pzIe zs}
	\p_z \I_\e(\qa +  \e^2  g, \Va +\e^2  W )(t,\zs)  = \p_z \Iea(t,\zs)  + \\  \e^2 \Big[  \p_g \p_z \I_\e^\ast (t,\zs) g +  (\p_V \p_z \I_\e^\ast \cdot W)(t,\zs) \Big] +  O(\e^5) \norm{(g,W)}.
	\end{multline}
\end{lem}
\begin{proof}[Proof of \cref{control pzgV2Ieps}]$ $\\
	We will denote  $f(p):=\p_z \I_\e(\qea + p g, \Vea + p W)(t,z)$. We recognize in the formula \cref{lin pzIe zs} a Taylor expansion of $f$. Then, to prove the estimate of \cref{lin pzIe zs} it is sufficient to  bound $ f''(\e^2)$ uniformly: 
	\begin{align*}
	f''(\e^2) \leq O(\e) \norm{(g,W)}.
	\end{align*}
	The formula for $f''$ is very long, so for clarity we will denote respectively $A_\e(p)$ the numerator  and $B_\e(p)$ the denominator of $f(p)$, so that when we differentiate we  have the structure :
	\begin{align}\label{eq g''}
	f''(p) = \frac{A_\e''(p)}{B_\e(p)} - 2 \frac{A_\e'(p)B_\e'(p)}{B_\e(p)^2}- \frac{A_\e(p) B_\e''(p)}{B_\e(p)^2} + 2 \frac{A_\e(p) B_\e'(p)^2 }{B_\e(p)^3}.
	\end{align}
	%
	The numerator is defined as : 
	\begin{multline*}
	A_\e(p) := \ds\iint_{\R^2} \exp \Big[-Q(y_1,y_2) +2 \D(\Vea+  p W)(Y,t,\zs)- \e (\qea+  p g) (y_1+y_2) \Big] \\  \times \Big( \D(\p_z \Va + p W)  \Big)(Y,t,\zs)  d y_1 d y_2,
	\end{multline*}
	while the denominator reads :
	\begin{align*}
	B_\e(p):= \int_{\R} e^{-\frac12\abs{y}^2}  \exp \Big(-\e (\qa+pg )  y  +  \D^\ast(\Va+p W)(y,t) \Big)  \,  dy.
 	\end{align*}
	Therefore we will divide each term by $\Iea$ to simplify the notations, this will make appear  the measures $dG_\e^\ast, d N_\e^\ast$ introduced in \cref{def_dG,def_dN}. For instance  :
	\begin{multline*}
	\frac{A_\e(p)}{\I_\e ^ \ast(t,\zs) B_\e(p) } := \frac{\ds \iint_{\R^2 } d G_\e^ \ast(Y,t,\zs)\exp \Big(-\e pg (y_1+y_2)  +2  p \D (W)(Y,t,\zs) \Big)} {\ds \int_{\R} d N_\e^\ast(y,t)  \exp \Big(  p \D^\ast(W)(y,t)-\e  p g y   \Big)dy} \times \\  \Big[  \D (\p_z \Vea + p \p_z \Va)(Y,t,\zs) \Big]dy_1dy_2 .
	\end{multline*}
	We notice that any factor of the sum in \cref{eq g''} (divided by  $\I_\e^ \ast$) is a sum (and  a product) of terms of the form 
	\begin{align*}
	\frac{A_\e^{(j)}(p)B_\e^{(k)}(p)}{B_\e(p)  \I_\e^\ast(t,\zs)} = \frac{A_\e^{(j)}(p)}{\I_\e^\ast(t,\zs) B_\e(p)} \frac{B_\e^{(k)}(p)}{B_\e(p)},
	\end{align*}
	with $j=0,1,2$, $k=1,2$ and the constraint $j+k=2$.  It is rather convenient to bound separately each of those terms. For instance we deal with the second  one:
	\begin{align}
	\frac{A_\e ' (p) B_\e'(p)}{B_\e(p)^2 \I_\e^\ast(t,\zs)} = \frac{A_\e'(p)} {\I_\e^\ast(t,\zs) B_\e(p)} \frac{B_\e'(p)}{B_\e(p)},
	\end{align}
	The first term of this product is
	\begin{multline*}
	\frac{A_\e'(p)}{\I_\e ^ \ast(t,\zs) B_\e(p) } := \\
	\frac{\ds \iint_{\R^2} d G_\e^ \ast(Y,t,\zs) \Big[  \exp \Big(2 p \D (W)-\e g p (y_1+y_2)   \Big)  \D (  \p_z W)\Big] dy_1dy_2} { \ds \int_{\R} d N_\e^\ast(y,t)  \exp \Big( 2  \D^\ast(W)(y,t)-\e  g y   \Big) dy} + \\ \frac{\ds \iint_{\R^2 } \!\!\! d G_\e^ \ast(Y,t,\zs )  \exp \Big(\!2 p \D (W) \! -\! \e g p (y_1\! + \! y_2)   \Big) 2 \D (\p_z \Vea\! + \! p \p_z W) \! \Big( \!\D (W) \! - \! \e g  (y_1+y_2)  \Big) dy_1dy_2 } {\ds \int_{\R} d N_\e^\ast(y,t)  \exp \Big( 2  \D^\ast(W)(y,t)-\e  g y   \Big) dy}. 
	\end{multline*}
	The numerator and denominator can be bounded  by estimating naively $\D$ :
	\begin{multline}
	\abs{\frac{A_\e'(p)}{B_\e(p) \I_\e^\ast(t,z) }} \leq  \\ 
	\frac{ \ds \iint_{\R^2} d G_\e^ \ast(Y,t,\zs)  \exp \Big( 3 \e \norm{(g,W)} (\abs{y_1}+\abs{y_2})   \Big)   \e (\abs{y_1}+\abs{y_2}) \norm{(g,W)}  dy_1dy_2} {\ds  \int_{\R} d N_\e^\ast(y,t)  \exp \Big( -3 \e \norm{(g,W)} \abs{y}   \Big) dy}  + \\  \frac{ \ds \iint_{\R^2}  d G_\e^ \ast (Y,t,\zs )  \exp  \Big( 3 \e \norm{(g,W)} \! (\abs{y_1} +  \abs{y_2})    \Big)   \e^2  (\abs{y_1} + \abs{y_2})^2  (3\norm{(g,W)} + \! 2 K^\ast)  3 \! \norm{(g,W)} \! dy_1dy_2} {\ds  \int_{\R} d N_\e^\ast(y,t) \exp \Big(-3 \e \norm{(g,W)} \abs{y}   \Big) dy}.
	\end{multline}
	Therefore, we only get moments of a Gaussian distribution, so the previous bound is  in fact
	\begin{align*}
	\abs{\frac{A_\e'(p)}{B_\e(p) \I_\e^\ast(t,z) }} \leq O( \e ) \norm{(g,W)}.
	\end{align*}
	With the exact same arguments but more convoluted formulas, one shows that 
	\begin{align}\label{eq ej}
	\abs{\frac{A_\e''(p)}{B_\e(p) \I_\e^\ast(t,z) }} \leq O( \e )\norm{(g,W)}.
	\end{align}
	For the quotients of $B$ in \cref{eq g''}, we loose the structure of the measures $d G_\e^ \ast$ and $d N_\e^ \ast$, but they are replaced by an actual Gaussian measure $ \exp(-y^2/2)$. Therefore, with the same arguments as before, we bound the quotient by the  moments of a Gaussian distribution. For instance,
	\begin{align}
	\nonumber \abs{ \frac{B_\e'(p)}{B_\e(p)} } & = \abs{ \frac{\ds \int_{\R} e^{-\frac12\abs{y}^2}  \exp \Big(2  \D^ \ast (\Vea +  p W)-\e (\qea + g p) y   \Big)  \Big( 2\D^ \ast (W)-\e g y \Big) dy}{\ds \int_{\R} e^{-\frac12\abs{y}^2} \exp \Big(2  \D^ \ast (\Vea +  p W)-\e (\qea + g p) y   \Big) dy}} \\
	&  \nonumber \leq \frac{\ds \int_{\R} e^{-\frac12\abs{y}^2} \exp \Big( 3 \e \abs{y} K^\ast + 3 \e \norm{(g,W)} \abs{y}\Big)  \Big( 3 \e \norm{(g,W)} \abs{y}  \Big) dy }{\ds \int_{\R} e^{-\frac12\abs{y}^2}  \exp \Big(-3 \e K^\ast \abs{y}-3 \e \norm{(g,W)}  \abs{y} \Big) dy}, \\
	&\label{eq ek} \leq  O(\e )\norm{(g,W)}.
	\end{align}
When multiplying each term of \cref{eq ej} by  \cref{eq ek} and then combining them yields the desired estimate  result, given the separation of terms made in \cref{eq g''} : 
	\begin{align*}
	\abs{\frac{f''(p)}{\Iea(t,z) }} \leq O(\e)\norm{(g,W)}.
	\end{align*}
	Thanks to \cref{estim Iea} , \cref{control pzgV2Ieps} is proven.
\end{proof}

\section{Linearized equation for $\kappa_\e$, convergence of $p_\e$}
\subsection{Uniform boundedness of  $\kappa_\e$}\label{sec kappa}$ $\\
Thanks to the estimates of the previous sections, every useful tools to look at  the perturbation $\kappa_\e$ are made available. We recall that our final goal is to show that $\kappa_\e$ is bounded as it is the perturbation from $\qa$, see \cref{perturb kappa}. We show in this section that one gets an approximated Ordinary Differential Equation (ODE) on $\kappa_\e$ with good properties when linearizing, see \cref{prop eq kde}. It is obtained by differentiating \cref{perturb Weps} and evaluating at $z=\zs$. This is exactly what suggested the spectral analysis of the formal linearized operator in \Cref{tab}. Now, thanks to our previous set of estimates of \cref{sec estim Ie}, we are able to carefully justify our  linearization. Finally, the limit ODE we introduced for $\qa$ in \cref{def qda}  will appear clearly when we do our analysis to balance contributions of smaller order. 

To shortcut expressions, we introduce the following alternative notations for all $t,z \in\R_+ \times \R$:
\begin{align}\label{def M}
\Xi_\e(t,z) & := W_\e(t,z)-2 W_\e(t, \bz(t)).
\end{align}
Compared to previous sections, and for the rest of this article, we will work in the space $\F$ that is well suited to measure $W_\e$ and build the linearization results, here for $\kappa_\e$. All our previous estimates that were established in $\E$ remain true in $\F$.
\begin{prop}[Equation on $\kappa_\e$]\label{prop eq kde}$ $\\
	For any ball $B$ of $\R \times \F$ there exists a constant $\e_B$ that depends only on $B$ such that  if $(\kappa_\e,W_\e) \in B$ is a solution of \cref{perturb Weps}, then for all $\e\leq \e_B$, $\kappa_\e$ is a solution of the following ODE: 
	\begin{align}\label{EDO kappa}
		-\kde(t)  = R_\e^\ast(t)  \kappa_\e   +  O^\ast(1) \norm{W_\e}_\F  + O^\ast(1)+ O(\e)\norm{(\kappa_\e,W_\e)}.
	\end{align}
	where the $O(\e)$ depends only on $B$, and $R_\e^\ast$ are defined in \cref{control pzgVIea}.
\end{prop}
\begin{proof}[Proof of \cref{prop eq kde}]$ $\\
	As announced, one starts by differentiating \cref{perturb Weps}. This yields, with the notation $\Xi_\e$ introduced in \cref{def M} : 
	\begin{multline*}
	\p_z M(t,z) -\e^2  \qdea(t) -\e^2 \p_z \p_t \Vea(t,z)  -  \e^4  \kde(t)-\e^4 \p_z \p_t W_\e(t,z)  = \\ 
	M(t,z) 	\p_z  \I_\e(\qea + \e^2 \kappa_\e, \Vea + \e^2 W_\e)(t,z)  \exp( \e^2\Xi_\e(t,z))  \\
	+  \p_z M(t,z)  \I_\e(\qea + \e^2 \kappa_\e, \Vea + \e^2 W_\e)(t,z)  \exp( \e^2\Xi_\e(t,z))   \\
	+  \e^2M(t,z) \I_\e(\qea + \e^2 \kappa_\e, \Vea + \e^2 W_\e)(t,z)  \exp( \e^2\Xi_\e(t,z) ) \p_z \Xi_\e(t,z) .
	\end{multline*}
	When we evaluate the expression at $z = \zs$, the last two terms vanish, since $\p_z M(t,\zs)=\p_z \Xi_\e(t,\zs)=0$. Therefore, the equation becomes, since $\Xi_\e(t,\zs)=0$ and $M(t,\zs)=1$,
	\begin{align}\label{eq kde int}
	-\e^2  \qdea(t) -\e^2 \p_z \p_t \Vea(t,\zs)  -\e^4  \kde(t)-\e^4 \p_z \p_t W_\e(t,\zs)  = 
	\p_z  \I_\e(\qea + \e^2 \kappa_\e, \Vea + \e^2 W_\e)(t,\zs).
	\end{align}
	We then use directly the linearization result of \cref{control pzgV2Ieps} that we prepared for that purpose : 
	\begin{multline}\label{eq kappae}
	\p_z \I_\e(\qea + \e^2 \kappa_\e, \Vea + \e^2 W_\e)(t,\zs ) = \\ \p_z \I_\e^\ast(t, \zs ) + \e^2 \Big[  \p_g \p_z \I_\e^\ast (t,\zs) \kappa_\e +  (\p_V \p_z \I_\e^\ast \cdot W_\e)(t,\zs) \Big] 	+ O(\e^5) \norm{(\kappa_\e,W_\e)}.
	\end{multline}
	We see that for most of the terms, we provided a careful estimate in the previous \cref{sec estim Ie}. First, by \cref{link qa pzIea},
	\begin{align*}
	\p_z \Iea(t,\zs) = \e^2 \left( m''(\zs) \qa(t) - \frac{ m^{(3)}(\zs)}{2}\right) + O^\ast(\e^4).
	\end{align*}
	Plugging this in the asymptotic development of \cref{eq kappae}, we get the following :
	\begin{multline*}
	\p_z \I_\e(\qea + \e^2 \kappa_\e, \Vea + \e^2 W_\e)(t,\zs ) = \e^2 \left( m''(\zs(t)) \qa(t) - \frac{m^{(3)}(\zs(t))}{2}\right) + \\ \e^2 \Big[ \p_g \p_z \I_\e^\ast (t,\zs ) \kappa_\e +  \p_V \p_z \I_\e^\ast \cdot W_\e(t,\zs) \Big]  + O^\ast(\e^4) + O(\e^5)\norm{(\kappa_\e,W_\e)}.
	\end{multline*}
	Combining this with the \cref{control pzgVIea} where we got precise estimates at the point $\zs$, we complete the expansion of $\p_z \I_\e$: 
	\begin{multline*}
	\p_z \I_\e(\qea + \e^2 \kappa_\e, \Vea + \e^2 W_\e)(t,\zs ) = \\\e^2 \left( m''(\zs) \qa(t) - \frac{m^{(3)}(\zs)}{2}\right) +  \e^4  R_\e^\ast(t)  \kappa_\e  +  O^\ast(\e^4) \norm{W_\e}_\F +  O^\ast(\e^4) + O(\e^5)\norm{(\kappa_\e,W_\e)}.
	\end{multline*}
	When we turn back to \cref{eq kde int}, we have shown at this point the following relationship:
	\begin{multline}\label{eq kde int2}
	-\e^2  \qdea(t) -\e^2 \p_z \p_t \Vea(t,\zs)  -\e^4  \kde(t)-\e^4 \p_z \p_t W_\e(t,\zs) =  \\
	\e^2 \left( m''(\zs) \qa(t) - \frac{m^{(3)}(\zs)}{2}\right) +  \e^4  R_\e^\ast(t)  \kappa_\e  +  O^\ast(\e^4)\norm{W_\e}_\F + O^\ast(\e^4) + O(\e^5)\norm{(\kappa_\e,W_\e)}.
	\end{multline}
	To get a stable equation on $\kappa_\e$, the terms of order $\e^2$ must cancel out. This is precisely the role played by the dynamics of $\qa$ defined in \cref{def qda}. To see it, we just rewrite a term of \cref{eq kde int2} using that $\p_z \Va(t,\zs)=0$, and  \cref{prop Va}:
	\begin{align*}
	\p_z \p_t \Vea(t,\zs)  = m'(\zs) \p_z^2 \Vea(t,\zs)= 2 m'(\zs) m''(\zs).
	\end{align*}
	Therefore, we recognize that by definition of $\qa$ in \cref{def qda}, the following terms cancel:
	\begin{align*}
	\e^2 \left( \qda(t) + m''(\zs) \qa(t) - \frac{m^{(3)}(\zs)}{2} +2 m''(\zs) m'(\zs) \right)=0.
	\end{align*}
	We then rewrite the second term of \cref{eq kde int2} of order $\e^4$ :
	\begin{align*}
	\p_z \p_t W_\e(t,\zs)  = m'(\zs) \p_z^2 W_\e(t,\zs) = O^\ast(1) \norm{W_\e}_\F.
	\end{align*}
	Finally, we deduce from \cref{eq kde int2} the following relationship:
	\begin{align*}
	-\kde(t)  = R_\e^\ast(t)  \kappa_\e   +  O^\ast(1) \norm{W_\e}_\F  + O^\ast(1)+ O(\e)\norm{(\kappa_\e,W_\e)}.
	\end{align*} We have proven \cref{EDO kappa}.
\end{proof}
In this ODE solved by $\kappa_\e$, each term play a separate part. First the function $R_\e^\ast$ is what guarantees the stability of $\kappa_\e$ because it is negative for large times. The other terms come from our perturbative analysis methodology. The term $O^\ast(1) + O(\e)\norm{(\kappa_\e,W_\e)}$ measures the error made when linearizing to obtain the ODE, and it ensures that it is of superior order in $\e$ except for the part that comes from the reference point of our linearization : $O^\ast(1)$. Interestingly there is also an error term that is not of superior order when linearizing, $O^\ast(1) \norm{W_\e}_\F $, but what saves our contraction argument of \cref{sec proof} is that this term only involves $W_\e$, which we can bound independently, see \cref{sec lin stab}.   

\subsection{Equation on $p_\e$}\label{equation pe}$ $\\
We did not perturb the number $p_\e$ as we did for $(q_\e,V_\e)$ since it can be straightforwardly computed from our reference \cref{perturb Weps}.  
Given the spectral decomposition of heuristics \cref{sec heuristics}, it is consistent to evaluate \cref{perturb Weps} at $z=\zs$ to gain the necessary information upon $p_\e$. This yields :
\begin{align}\label{def pe}
1 -\e^2\Big(\pde(t)  + m'(\zs) \qea(t)  \Big) 
-\e^4  m'(\zs) \kappa_\e(t)   = \I_\e( \qea+\e^2 \kappa_\e, \Vea + \e^2 W_\e)(t,\zs ).
\end{align}
Thanks to  \cref{prop Vea bound,prop lin Ieps}, and as long as $\kappa_\e$ is bounded, which we will show in \cref{sec proof},
\begin{align*}
\e^2\Big(\pde(t)  +m'(\zs(t)) \qea(t)  \Big) =   O(\e^2)
\end{align*}
In this last equation, the order of precision is not enough to recover the equation on $\pa$ when $\e \to 0$. The problem is that the linearization of $ \I_\e$ made in \cref{Iea 1Oast} is a little too rough. Coming back to \cref{estim Iea}, we make the more precise following estimate : 
\begin{align}\label{precise lin Iea}
\Iea(t,\zs) = 1 - \frac{\e^2}{2} \p_z^2 \Va(t,\zs) + O^\ast(\e^4).
\end{align}
The proof of this result is a direct adaptation of the one of \cref{estim Iea}, by making Taylor expansions up to the fourth derivative of $\Va$, as made possible by the introduction of $\East$, see \cref{def East}. This involves computing  the moments of the Gaussian distribution $\exp(-Q) $ : 
\begin{align}\label{mom Q}
\ds \dfrac{1}{\sqrt 2 \pi} \iint_{\R^2} e^{-Q(y_1,y_2)}(y_1^2+y_2^2)dy_1dy_2= \frac12.
\end{align}
By plugging \cref{precise lin Iea} into \cref{def pe}, and using \cref{Iea O3}, we find 
\begin{align}
\nonumber \pde(t) + m'(\zs) \qea(t) & = \frac{\p_z^2 \Va(t,\zs)}{2} +O(\e^2), \\
& = m''(\zs) + O(\e^2).\label{def pde precise}
\end{align}
We used \cref{link pVea pm} for the last equality. From \cref{def pde precise}, the convergence of $p_\e$ towards $\pa$ defined by \cref{def pda}, stated in \cref{main theo} is straightforward. 

\section{Linearization results}\label{sec lin}
We finally tackle the complete linearization of \cref{perturb Weps}. A foretaste was given when we studied the equation on $\kappa_\e$, however it was local since we had beforehand evaluated at $\zs(t)$. Here, we will provide global (in space) results. 
\subsection{Linearization for $W_\e$}\label{sec:linearization-for-we}$ $\\
A first step is to control the function $\Xi_\e$, which we recall, is a byproduct of $W_\e$, introduced in \cref{def M}.
\begin{lem}[Control of $\Xi_\e$]\label{control M}$ $ \\
	For any ball $B$ of $\F$, there exists a constant $\e_B$ that depends only on $B$ such that for all $\e\leq \e_B,$ if  $W_\e \in B$ , $\Xi_\e$ defined in \cref{def M} verifies
	\begin{align*}
	\exp(\e^2 \Xi_\e(t,z)) = 1 + \e^2 \Xi_\e(t,z) + O(\e^4)\norm{W_\e}_\F.
	\end{align*}	
	where $O(\e^4)$ depends only on the ball $B$. 
\end{lem}
\begin{proof}[Proof of \cref{control M}]$ $\\
	By the choice of the norm in $\F$, and in the setting of $W_\e \in B$ we have the uniform control for all $t,z$ :
	\begin{align*}
	\abs{\Xi_\e(t,z)}\leq \norm{W_\e}_\F.
	\end{align*}
	Then, by performing an exact Taylor expansion, there exists $0 <\xi < 1$ such that 
	\begin{align*}
	\exp(\e^2 \Xi_\e(t,z)) = 1 + \e^2 \Xi_\e(t,z) + \frac{\e^4}{2} \Xi_\e(t,z)^2 \exp \Big( \e^2 \xi \Xi_\e(t,z) \Big).
	\end{align*}
	To conclude we uniformly bound the rest for $\e^2 \leq 1/ \norm{W_\e}_\F$ :
	\begin{align*}
	\abs{  \frac{\e^4}{2} \Xi_\e(t,z)^2 \exp \Big( \e^2 \xi \Xi_\e(t,z) \Big) } \leq e\e^4 \frac{\norm{W_\e}_\F ^2}{2} .
	\end{align*}
\end{proof}
This first result is prototypical of the tools we will employ to linearize  the problem \cref{perturb Weps} solved by $(\kappa_\e,W_\e)$. 
We now write the linearized problem verified by $W_\e$.
\begin{prop}[Linearization for $W_\e$]\label{prop lin Weps}$ $\\
	For any ball $B$ of $\R \times \F$, there exists a constant $\e_B$ that depends only on $B$ such that for all $\e\leq \e_B$, any pair $(\kappa_\e,W_\e) \in B$  solution of \cref{perturb Weps} verifies the following estimate :
	\begin{align}\label{PL Weps}
	-\e^2   \p_t W_\e(t,z) = M(t,z)  \Big( \Xi_\e(t,z) + \Oun \norm{(\kappa_\e,W_\e)} \Big) ,
	\end{align}
	where $O(\e)$ depends only on $B$. 
\end{prop}
\begin{proof}[Proof of \cref{prop lin Weps}]$ $\\
	One starts from the equation \cref{perturb Weps}, 
	\begin{multline}\label{eq linWint1}
	M(t,z)-\e^2\Big( \pde(t) + m'(\zs) \qa(t) + \qdea(t)(z-\zs)  + \p_t \Vea(t,z) \Big)  \\ -\e^4 \Big(  \kde(t)(z-\zs ) + m'(\zs) \kappa_\e(t) + \p_t W_\e(t,z) \Big)  \\
	= M(t,z)  \I_\e( \qea+\e^2 \kappa_\e, \Vea + \e^2 W_\e)(t,z)  \exp \left( \e^2 \Xi_\e(t,z)  \right)	
	\end{multline}
	Thanks to \cref{control M,prop lin Ieps} where we linearized $\I_\e$ and the term in $\Xi_\e$, one can expand the right hand side : 
	\begin{multline}
	M(t,z)  \I_\e( \qea+\e^2 \kappa_\e, \Vea + \e^2 W_\e)(t,z) \exp \left( \e^2 \Xi_\e(t,z)  \right)  =  \\ M(t,z) \Big( 1 + O^\ast(\e^2)  + O(\e^3)\norm{(\kappa_\e,W_\e)} \Big) \Big(1 + \e^2 \Xi_\e(t,z) + O(\e^4) \norm{(\kappa_\e,W_\e)}  \Big) \\
	\label{eq linWint} = M(t,z) + \e^2 M(t,z)  \Xi_\e(t,z) +  M(t,z)  \Big( O^\ast(\e^2)  + O(\e^3)\norm{(\kappa_\e,W_\e)} \Big).
	\end{multline}
	The left hand side  of \cref{eq linWint1} is a little bit more involved. We will use our previous work on $(p_\e,\kappa_\e)$. First, thanks to \cref{def pe} that states the relationship verified by $p_\e$, we have 
	\begin{align*}
	-\e^2\Big( \pde(t) + m'(\zs) \qa(t) \Big) - \e^4 \kappa_\e m'(\zs) = 1 - \I_\e(\qa+\e^2 \kappa_\e,\Va+\e^2 W_\e)(t,\zs).
	\end{align*}
	We then use \cref{prop lin Ieps} about the linearization of $\I_\e$ to get that 
	\begin{align}\label{eq linWint3}
	-\e^2\Big( \pde(t) + m'(\zs) \qa(t) \Big) - \e^4 \kappa_\e m'(\zs) = O^\ast(\e^2) + O(\e^3)\norm{(\kappa_\e,W_\e)}.
	\end{align}
	From \cref{prop Vea bound}, we have the following uniform bound : 
	\begin{equation}\label{eq linWint4}
	\abs{ \p_t \Vea (t,z) } \leq K^\ast.
	\end{equation}
	Thanks to our preliminary work on $\kappa_\e$, and more precisely the \cref{eq kde int2} we know that  
	\begin{align*}
	\qda(t) + \e^2 \kde(t)=  O^\ast(1)  + O(\e)\norm{(\kappa_\e,W_\e)}.
	\end{align*} 
	Therefore, the affine terms are comparable to $M$, since $M$ is a superlinear function that admits a uniform lower bound by  hypothesis, see \cref{cond Gamma}:
	\begin{align}\label{eq linWint5}
	\abs{ \frac{ \Big(\qdea(t) +\e^2 \kde(t)\Big)(z-\zs) }{M(t,z) } }= O^\ast(1)  + O(\e)\norm{(\kappa_\e,W_\e)} .
	\end{align}
	When adding up the estimates of \cref{eq linWint4} and \cref{eq linWint5}, we have shown :
	\begin{multline}\label{eq linWint2}
	-\e^2\Big( \pde(t) + m'(\zs) \qa(t)+  \qdea(t)(z-\zs )  + \p_t \Vea(t,z) \Big)  \\- \e^4 \Big(  \kde(t)(z-\zs ) + m'(\zs(t)) \kappa_\e(t) + \p_t W_\e(t,z) \Big) 
	\\ = M(t,z) \Big( O^\ast(\e^2)+O(\e^3) \Big) \norm{(\kappa_\e,W_\e)} -\e^4 \p_t W_\e(t,z)  .
	\end{multline}
	We have divided by $M$ the relationships  \cref{eq linWint3} and \cref{eq linWint4}, which is possible  thanks to the uniform lower bound of $M$. 
	
	Finally, when putting together \cref{eq linWint5} and \cref{eq linWint} in \cref{eq linWint1}, the terms $M$ cancel each other, and we find \cref{PL Weps} factoring out $\e^2$.
\end{proof}
One can notice the similarity between what we just proved rigorously and the heuristics made in \cref{lin We heur}. From this result one can straightforwardly deduce a linear approximated equation verified by $\Xi_\e $.
\begin{cor}[Linearization in $\Xi_\e(t,z)$]\label{prop lin M}$ $\\
	For any ball $B$ of $\R \times \F$, there exists a constant $\e_B$ that depends only on $B$ such that for all $\e\leq \e_B$, any pair $(\kappa_\e,W_\e) \in B$ verifies the following estimate
	\begin{align}\label{PL M}
	\e^2 \p_t \Xi_\e(t,z) = M(t,z) \left( 2 \frac{ M(t, \bz(t))}{M(t,z)} \Xi_\e(t, \bar{z}) - \Xi_\e(t,z)  + O^\ast(1) + O(\e)\norm{(\kappa_\e,W_\e)} \right) 
	\end{align}
	where the $O(\e)$ depends only on $B$.
\end{cor}
\begin{remark}\label{remark linpzWe}\begin{itemize}[label=$\triangleright$]
		\item A careful reader may notice that the computation of $\p_t \Xi_\e$ yields a parasite term $ \e^2 \zds \p_z \Xi_\e(t,\bar z)$ not dealt by \cref{PL Weps}. However this is a lower order term since it verifies:
		\begin{align}
		\e^2 \zds(t) \p_z \Xi_\e(t,z) = O(\e^2)\norm{(\kappa_\e,W_\e)}
		\end{align}
\item 	Under the same assumption as \cref{prop lin M}, $W_\e$ also verifies the following linear equation:	\begin{align*}
	-\e^2 \p_t W_\e(t,z) = M(t,z) \Big( \Xi_\e(t,z) + O(1) \Big).
	\end{align*}
	However in \cref{sec lin stab}, we will study the stability of the solution of the linear problem. We will see that one needs precise estimates about the structure of the nonlinear negligible terms, which explains the more detailed \cref{PL Weps}, and is the purpose of all our previous sections.
		\end{itemize}
\end{remark}
\subsection{Linearization for $\p_z W_\e$}$ $\\
The computations for $\p_z W_\e$ are slightly more complex because of the differentiation of the triple product in the right-hand side \cref{perturb Weps}. However, the key point is that when we linearize $\I_\e(\qa+\e^2 \kappa_\e,\Va+\e^2 \kappa_\e)$ the derivatives of $\I_\e$ are negligible in $\e$. Therefore the intuitive linearized problem for $\p_z W_\e$, given by the derivation of the linearized equation for $W_\e$, actually holds true. This is the content of the following proposition :
\begin{prop}[Linearization in $\p_z W_\e$]\label{prop lin pzWeps}$ $\\
	For any ball $B$ of $\R \times \F$, there exists a constant $\e_B$ that depends only on $B$ such that for all $\e\leq \e_B$, any pair $(\kappa_\e,W_\e) \in B$  solution of \cref{perturb Weps}  verifies the following estimate :
	\begin{multline}\label{PL pzWeps+}
-	\e^2 \p_t \p_z W_\e(t,z)= M(t,z) \left(\p_z \Xi_\e(t,z) +  \frac{  O^\ast(1)  + O(\e)\norm{(\kappa_\e,W_\e)} }{\vphia(t,z)}\right) \\ + \p_z M(t,z) \Big(\Xi_\e(t,z)  +O^\ast(1) + O(\e) \norm{(\kappa_\e,W_\e)}\Big),
	\end{multline}
	where $O(\e)$ depends only on $B$. 
\end{prop}

\begin{proof}[Proof of \cref{prop lin pzWeps}]$ $\\
	One starts by differentiating \cref{perturb Weps} as in the proof of \cref{prop eq kde} to highlight $\kappa_\e$. This yields : 
	\begin{align*}
	\p_zM(t,z) -\e^2  \qdea(t) -\e^2 \p_z \p_t \Vea(t,z)  - & \e^4  \kde(t)-\e^4 \p_z \p_t W_\e(t,z)  = \\ 
	& M(t,z) 	\p_z  \I_\e(\qea + \e^2 \kappa_\e, \Vea + \e^2 W_\e)(t,z)  \exp( \e^2\Xi_\e(t,z))  \\
	+ & \, \p_z M(t,z)  \I_\e(\qea + \e^2 \kappa_\e, \Vea + \e^2 W_\e)(t,z)  \exp( \e^2\Xi_\e(t,z))   \\
	+ & \, \e^2  M(t,z) \I_\e(\qea  + \e^2 \kappa_\e, \Vea   +  \e^2 W_\e)(t,z)  \exp( \e^2\Xi_\e(t,z) ) \p_z \Xi_\e(t,z) .
	\end{align*}
	However contrary to the case where we were studying $\kde$, we will not evaluate in $\zs$. We introduce the notations $R_i$ corresponding to each of the three terms of the right hand side of the previous equation.
	We will linearize each $R_i$ starting with $R_1$ which we estimate  thanks to \cref{prop lindecay,control M}, paired with the estimate of \cref{decay pzIea} :
	\begin{align*}
	R_1  : = & \p_z  \I_\e(\qea + \e^2 \kappa_\e, \Vea + \e^2 W_\e)(t,z) M(t,z)  \exp( \e^2\Xi_\e(t,z)) \\
	= &   M(t,z)   \left(  \p_z \Iea(t,z)+ \frac{O(\e^3)\norm{(\kappa_\e,W_\e)}}{\vphia(t,z) }  \right) \Big(1 + \e^2\Xi_\e(t,z) + O(\e^4) \norm{(\kappa_\e,W_\e)} \Big) \\
	= &  M(t,z) \left( \frac{O^\ast(\e^2) + O(\e^3)\norm{(\kappa_\e,W_\e)}}{\vphia(t,z) } \right) \Big(1  + \e^2\Xi_\e(t,z) + O(\e^4)\norm{(\kappa_\e,W_\e)} \Big). 
	\end{align*}
	Therefore, the final contribution of $R_1$ is:
	\begin{align}\label{eq R1}
	R_1 =   M(t,z) \left( \frac{O^\ast(\e^2) + O(\e^3)\norm{(\kappa_\e,W_\e)}}{\vphia(t,z) } \right).
	\end{align}
	Next, one looks at $R_2$. Thanks to \cref{prop lin Ieps}, 
	\begin{align}
	R_2 :=  &  \nonumber \p_z M(t,z)  \I_\e(\qea + \e^2 \kappa_\e, \Vea + \e^2 W_\e)(t,z)    \exp( \e^2\Xi_\e(t,z)) \\
	=& \nonumber\p_z M(t,z) \Big( 1 + O^\ast(\e^2)+ O(\e^3)\norm{(\kappa_\e,W_\e)} \Big) \Big(1 + \e^2\Xi_\e(t,z) + O(\e^4) \norm{(\kappa_\e,W_\e)} \Big), \\
	= & \label{eq R2} \p_z M(t,z)  + \e^2 \p_zM(t,z)\Xi_\e(t,z) + \p_z M(t,z) \Big( O^\ast(\e^2)+ O(\e^3)\norm{(\kappa_\e,W_\e)} \Big). 
	\end{align}
	We finally tackle $R_3$ with the same techniques, using \cref{prop lin Ieps,control M} :
	\begin{align}
	\nonumber R_3 :=& \e^2 M(t,z) \I_\e(\qea + \e^2 \kappa_\e, \Vea + \e^2 W_\e)(t,z) \exp( \e^2\Xi_\e(t,z)) \p_z \Xi_\e(t,z), \\
	=& \nonumber  \e^2 M(t,z) \p_z \Xi_\e(t,z)\Big( 1 + O^\ast(\e^2)+ O(\e^3)\norm{(\kappa_\e,W_\e)} \Big)   \Big(1+\e^2 \Xi_\e(t,z) + O(\e^4)\norm{(\kappa_\e,W_\e)} \Big) ,\\
	\label{eq R3} = &  \e^2 M(t,z) +  M(t,z) \frac{O(\e^4) \norm{(\kappa_\e,W_\e)} }{\vphia(t,z)} .
	\end{align}
	In that last estimate, we chose to write $O^\ast(\e^4)$ as a regular $O(\e^4)$. 
	If we come back to our initial problem, when we assemble \cref{eq R1,eq R2,eq R3}, we obtain:
	\begin{multline}\label{eq linpzWe int}
	\p_zM(t,z) -\e^2  \qdea(t) -\e^2 \p_z \p_t \Vea(t,z)  -\e^4  \kde(t)-\e^4 \p_z \p_t W_\e(t,z)   \\ 
	=  \p_z M(t,z) +  \e^2  \p_z M(t,z)  \Big( \Xi_\e(t,z) +  O^\ast(1)+ O(\e)\norm{(\kappa_\e,W_\e)} \Big) + \\ 
	\e^2  M(t,z)  \left(\p_z \Xi_\e(t,z) 
	+ \frac{ O^\ast(1) + O(\e)\norm{(\kappa_\e,W_\e)} }{\vphia(t,z)}   \right).
	\end{multline}
We now deal with the left hand side of  \cref{eq linpzWe int}. First, the terms  $\p_z M(t,z) $ on each side cancel. Next, using the ODE that defines $\qa$ in \cref{def qda}, our linearized equation on $\kde$ stated in \cref{EDO kappa} and finally our bound of $\p_t \Va$ made in \cref{prop Vea bound}, we find: 
	\begin{align}\label{eq linpzWe int2}
	-\e^2  \Big( \qdea(t)  + \p_z \p_t \Vea(t,z)   + \e^2   \kde(t) \Big)=  O^\ast(\e^2)+ O(\e^3)\norm{(\kappa_\e,W_\e)}. 
	\end{align}
	Finally, if we divide by  $M$, the following estimate holds true since $\alpha<1$:
	\begin{align*}
	\abs{\frac{ O^\ast(\e^2) + O(\e^3)\norm{(\kappa_\e,W_\e)}}{M(t,z)}} \leq \frac{ O^\ast(\e^2) + O(\e^3)\norm{(\kappa_\e,W_\e)}}{\vphia(t,z)}.
	\end{align*}
	Plugging this into \cref{eq linpzWe int}, and dividing each side by $\e
	^2$, we therefore recover the relationship we wanted to prove :
	\begin{multline*}
	-\e^2 \p_t \p_z W_\e(t,z) = M(t,z)  \left( \p_z \Xi_\e(t,z) + \frac{ O^\ast(1) + O(\e)\norm{(\kappa_\e,W_\e)}}{\vphia(t,z)} \right)   \\ +\p_z M(t,z) \Big(\Xi_\e(t,z)  +O^\ast(1) + O(\e) \norm{(\kappa_\e,W_\e)}\Big).
	\end{multline*}
\end{proof}
We deduce straightforwardly a linearization result upon the quantity $ \p_z \Xi_\e$.
\begin{cor}[Linearization for $\p_z \Xi_\e(t,z)$]\label{prop lin N}$ $\\
	For any ball $B$ of $\R \times \F$, there exists a constant $\e_B$ that depends only on $B$ such that for all $\e\leq \e_B$, any pair $(\kappa_\e,W_\e) \in B$  solution of \cref{perturb Weps}  verifies the following estimate :
	\begin{multline*}
	\e^2 \p_t  \p_z \Xi_\e(t,z) =  M(t,z)  \left[  \frac{M(t,\bz)}{M(t,z)}  \p_z \Xi_\e(t,\bz) - \p_z \Xi_\e(t,z) +  \frac{O^\ast(1) + O(\e)\norm{(\kappa_\e,W_\e)} }{\vphia(t,z)} \right] \\
	+  \p_z M(t,z)  \left[ \frac{\p_z M(t,\bz)}{\p_z M(t,z)} \Xi_\e(t,\bz) -\Xi_\e(t,z) + O^\ast(1)+O(\e)\norm{(\kappa_\e,W_\e)} \right].
	\end{multline*}
	where the $O(\e)$ depends only on $B$.
\end{cor}
\subsection{Linearization for $\p_z^2 W_\e(t,z)$}$ $\\
We now tackle the linearized equation for $\p_z^2 W_\e$. 
\begin{prop}[Linearization for $\p_z^2 W_\e$]\label{prop lin pzzWeps}$ $\\
	For any ball $B$ of $\R \times \F$, there exists a constant $\e_B$ that depends only on $B$ such that for all $\e\leq \e_B$, any pair $(\kappa_\e,W_\e) \in B$  solution of \cref{perturb Weps} verifies the following estimate:
	\begin{multline}\label{PL pzzWeps+}
	-\e^2 \p_z^2 \p_t W_\e(t,t) =    \p_z^2 M(t,z) \Big( \Xi_\e(t,z)  + O^\ast(1)+ O(\e)\norm{(\kappa_\e,W_\e)} \Big)  \\ + 2  \p_zM(t,z) \left(  \p_z \Xi_\e(t,z) + \frac{O^\ast(1)+O(\e)\norm{(\kappa_\e,W_\e)}}{\vphia(t,z)} \right)    +  M(t,z) \left(  \p_z^2 \Xi_\e(t,z) + \frac{O^\ast(1)+O(\e)\norm{(\kappa_\e,W_\e)}}{\vphia(t,z)} \right).
	\end{multline}
	where the $O(\e)$ depends only on $B$. 
\end{prop}
We will choose later to write the second derivative $\p_z^2 \Xi_\e(t,z)$ in full : $\p_z^2 W_\e(t,z)- \frac12 \p_z^2 W_\e(t,\bz)$ in the next sections as  the factor $\frac12$ will be the key to ensure the uniform boundedness of $\p_z^2 W_\e$, see \cref{sec lin stab}.
\begin{proof}[Proof of \cref{prop lin pzzWeps}]$ $\\
	We start by differentiating twice \cref{perturb Weps}. This yields : 	\begin{align*}
	\p_z^2 M(t,z) -\e^2 \p_z^2 \p_t \Vea(t,z)  -\e^4 \p_z^2 \p_t W_\e(t,z)  = R_1 + R_2 + R_3 + R_4 + R_5+R_6,
	\end{align*} 
	with the following notations :
	\begin{align*}
	& R_1  : = \p_z^2  \I_\e(\qea + \e^2 \kappa_\e, \Vea + \e^2 W_\e)(t,z) M(t,z) \exp( \e^2\Xi_\e(t,z)) ,  \\ 
	& R_2 := 2 \p_z M(t,z)\p_z  \I_\e(\qea + \e^2 \kappa_\e, \Vea + \e^2 W_\e)(t,z)    \exp( \e^2\Xi_\e(t,z)), \\
	& R_3 := 2 M(t,z) \e^2 \p_z  \I_\e(\qea + \e^2 \kappa_\e, \Vea + \e^2 W_\e)(t,z)  \exp( \e^2\Xi_\e(t,z))  \p_z \Xi_\e(t,z) , \\
	& R_4 :=  \I_\e(\qea + \e^2 \kappa_\e, \Vea + \e^2 W_\e)(t,z)  \p_z^2 M(t,z)   \exp( \e^2\Xi_\e(t,z)),\\
	& R_5 := 2 \e^2 \I_\e(\qea + \e^2 \kappa_\e, \Vea + \e^2 W_\e)(t,z) \p_z M(t,z) \exp( \e^2\Xi_\e(t,z)) \p_z \Xi_\e(t,z),
	\end{align*}
	and finally :
	\begin{align*}
	R_6 := \e^2M (t,z)  \I_\e(\qea + \e^2 \kappa_\e, \Vea + \e^2 W_\e)(t,z)   \exp( \e^2\Xi_\e(t,z))  \Big( \e^2 \p_z \Xi_\e(t,z)^2 + \p_z^2 \Xi_\e(t,z) \Big). 
	\end{align*}	
	We will estimate each term separately, starting with $R_1$, for which we apply the \cref{prop lindecay,control M,decay pzIea} : 
	\begin{align*}
	R_1 = &   M(t,z)   \left(  \p_z^2 \Iea(t,z)+ \frac{O(\e^3)\norm{(\kappa_\e,W_\e)}}{\vphia(t,z) }  \right) \Big(1 + \e^2\Xi_\e(t,z) + O(\e^4) \norm{(\kappa_\e,W_\e)} \Big) \\
	= &  M(t,z) \left( \frac{O^\ast(\e^2) + O(\e^3)\norm{(\kappa_\e,W_\e)}}{\vphia(t,z)} \right) \Big(1  + \e^2\Xi_\e(t,z) + O(\e^4)\norm{(\kappa_\e,W_\e)} \Big). \end{align*}
	Therefore, the final estimate of $R_1$ is :
	\begin{align}\label{eq R1b}
	R_1 =   M(t,z) \left( \frac{O^\ast(\e^2) + O(\e^3)\norm{(\kappa_\e,W_\e)}}{\vphia(t,z) } \right).
	\end{align}
	Next, for the other term $R_2$ we use \cref{prop lindecay,decay pzIea}  :
	\begin{align*}
	R_2 = & 2 \left( \p_z \Iea(t,z) + \frac{  O(\e^3)\norm{(\kappa_\e,W_\e)}}{\vphia(t,z)} \right) \p_z M(t,z)\Big( 1 + \e^2\Xi_\e(t,z) + O(\e^4) \norm{(\kappa_\e,W_\e)} \Big), \\
	= & 2 \p_z M(t,z)  \left( \frac{O^\ast(\e^2) + O(\e^3) \norm{(\kappa_\e,W_\e)}}{\vphia(t,z)} \right) \Big( 1 + \e^2\Xi_\e(t,z)  + O(\e^4) \norm{(\kappa_\e,W_\e)} \Big).
	\end{align*}
	We can simplify this expression :
	\begin{align}\label{eq R2b}
	R_2 & =  \p_z M(t,z)  \left( \frac{O^\ast(\e^2) + O(\e^3) \norm{(\kappa_\e,W_\e)}}{\vphia(t,z)} \right).
	\end{align}
	The term, $R_3$ will not contribute at the order $\e^2$, because of  \cref{decay pzIea}, and $\abs{\p_z \Xi_\e(t,z) } \leq \norm{W_\e}_\F $:
	\begin{align}
	\nonumber R_3 = &  2\e^2 M(t,z) \p_z \Xi_\e(t,z )\left( \frac{ O^\ast(\e^2) + O(\e^3) \norm{( \kappa_\e ,W_\e) }} {\vphia(t,z)} \right)  \Big( 1+\e^2\Xi_\e(t,z) + O(\e^4) \norm{(\kappa_\e,W_\e)} \Big) \\
	\label{eq R3b}=&   \frac{O(\e^3) \norm{( \kappa_\e ,W_\e) }}{\vphia(t,z)} M(t,z).
	\end{align} 
	For $R_4$, zeroth order terms are more entangled. With  \cref{prop lin Ieps,control M} :
	\begin{align}
	\nonumber R_4 & =  \p_z^2 M(t,z)   \Big(1 + O^\ast(\e^2) + O(\e^3)  \norm{( \kappa_\e ,W_\e) }\Big)\Big(1 + \e^2\Xi_\e(t,z) + O(\e^4) \norm{( \kappa_\e ,W_\e) }\Big), \\
	&\label{eq R4b} = \p_z^2 M(t,z) + \e^2 \p_z^2 M(t,z) \Big( \Xi_\e(t,z) + O^\ast(1) + O(\e) \norm{(\kappa_\e,W_\e)} \Big).
	\end{align}
	We see in $R_4$ the appearance of the term $ \e^2 \p_z^2 M(t,z) \Xi_\e(t,z)$ that is also in \cref{PL pzzWeps+}, and so it is a good opportunity to do a first a summary of the computations when adding \cref{eq R1b,eq R2b,eq R3b,eq R4b} :
	\begin{multline}\label{eq Rsum14}
	R_1 + R_2 + R_3  + R_4 =   \p_z^2 M(t,z) + \e^2 \p_z^2 M(t,z) \Big( \Xi_\e(t,z) +   O^\ast(1)+O(\e)\norm{(\kappa_\e,W_\e)} \Big) \\ + \e^2 M(t,z)\frac{O^\ast(1)+O(\e)\norm{(\kappa_\e,W_\e)}}{\vphia(t,z)} + \e^2 \p_z M(t,z)\frac{O^\ast(1)+O(\e)\norm{(\kappa_\e,W_\e)}}{\vphia(t,z)}  .%
	\end{multline}
	We continue the estimations by looking at $R_5$, thanks to \cref{prop lin Ieps} :
	\begin{align}
	\nonumber R_5 = &  2 \e^2 \p_z M(t,z) \p_z \Xi_\e(t,z)\Big(1 +O^\ast(\e^2)+ O(\e^3)\norm{(\kappa_\e,W_\e)} \Big) \Big[ 1 + \e^2\Xi_\e(t,z) + O(\e^4) \norm{(\kappa_\e,W_\e)} \Big], \\
	\label{eq R5}=&  2 \e^2  \p_z M(t,z) \p_z \Xi_\e(t,z) + \e^2  \p_z M(t,z)\frac{ O^\ast(\e) +O(\e^2) \norm{(\kappa_\e,W_\e)} }{\vphia(t,z)}.
	\end{align}
	Finally, we tackle the last term, $R_6$, with \cref{prop lin Ieps}
	\begin{align}
\nonumber	R_6   = \e^2 M (t,z)  \Big( 1 + O^\ast(\e^2)+ O(\e^3)\norm{(\kappa_\e,W_\e)} \Big)  \Big(1 +  \e^2\Xi_\e(t,z) +O(\e^4) \norm{(\kappa_\e,W_\e)} \Big) \\ \nonumber \times
	\left( \frac{O(\e^2) \norm{(\kappa_\e,W_\e)}}{\vphia(t,z)}+ \p_z^2\Xi_\e(t,z) \right),\\
\label{eq R6}  =   \e^2 M(t,z) \p_z^2 \Xi_\e(t,z) + \e^2 M(t,z) \frac{O(\e^2) \norm{(\kappa_\e,W_\e)}}{\vphia(t,z)} .
	\end{align}
	Thanks to those last two estimates \cref{eq R5} and \cref{eq R6}, that we add with the previous result of \cref{eq Rsum14}, we obtain for the full equation : 
	\begin{multline*}
	\p_z^2 M(t,z) -\e^2 \p_z^2 \p_t \Vea(t,z)  -\e^4 \p_z^2 \p_t W_\e(t,z) = \\   \p_z^2 M(t,z)   + \e^2 \p_z^2 M(t,z) \Big( \Xi_\e(t,z)  + O^\ast(1) + O(\e)\norm{(\kappa_\e,W_\e)} \Big)   +  	\\ 2 \e^2  \p_zM(t,z) \left(  \p_z \Xi_\e(t,z)  + \frac{ O^\ast(1) `+  O(\e)\norm{(\kappa_\e,W_\e)} }{\vphia(t,z)}\right) \\  + \e^2 M(t,z)\left(  \p_z^2 \Xi_\e(t,z)+ \frac{O^\ast(1)+O(\e)\norm{(\kappa_\e,W_\e)}}{\vphia(t,z)} \right).
	\end{multline*}
	Thanks to \cref{prop Vea bound} we know that $\norm{\e^2 \p_z^2 \p_t \Vea(t,z)}_\infty \leq O^\ast(\e^2)$. Then,
	\begin{multline*}
	-\e^4 \p_z^2 \p_t W_\e(t,t) =   \e^2 \p_z^2 M(t,z) \Big( \Xi_\e(t,z)  + O^\ast(1)+ O(\e)\norm{(\kappa_\e,W_\e)} \Big)  \\ + 2 \e^2  \p_zM(t,z) \left(  \p_z \Xi_\e(t,z) + \frac{O^\ast(1)+O(\e)\norm{(\kappa_\e,W_\e)}}{\vphia(t,z)} \right)    \\+ \e^2 M(t,z) \left(  \p_z^2 \Xi_\e(t,z) + \frac{O^\ast(1)+O(\e)\norm{(\kappa_\e,W_\e)}}{\vphia(t,z)} \right).
	\end{multline*}
	which proves \cref{PL pzzWeps+} after dividing by $\e^2$.
\end{proof}
\subsection{Linearization of $\p_z^3 W_\e(t,z)$}$ $\\
Our last linearized equation is the one for $\p_z^3 W_\e$ and we proceed with the same technique, with slightly more complex formulas.
\begin{prop}[Linearization in $\p_z^3 W_\e$]\label{prop lin pzzzWeps}$ $\\
	For any ball $B$ of $\R \times \F$, there exists a constant $\e_B$ that depends only on $B$ such that for all $\e\leq \e_B$, any pair $(\kappa_\e,W_\e) \in B$  solution of \cref{perturb Weps} verifies the following estimate:
	\begin{multline}\label{PL pzzzWeps+}
	-\e^2 \p_t  \p_z^3 W_\e(t,z) =  \p_z^3 M(t,z) \Big( \Xi_\e(t,z) + O^\ast(1) + O(\e) \norm{(\kappa_\e,W_\e)} \Big)  \\ + 3  \p_z^2 M (t,z)  \left( \p_z \Xi_\e(t,z) +\frac{O^\ast(1)  +  O(\e) \norm{(\kappa_\e,W_\e)}}{\vphia(t,z)} \right)   \\ +  	 3  \p_z M(t,z) \times \left(  \p_z^2 \Xi_\e(t,z) + \frac{O^\ast(1) + O(\e) \norm{(\kappa_\e,W_\e)}}{\vphia(t,z)} \right)  \\
	+ 	  M(t,z)    \left(  \p_z^3 \Xi_\e(t,z)   +  \frac{\norm{\vphia \p_z^3 W_\e}_\infty}{2^{1-\alpha} \vphia(t,z)} + \right.   \left. \frac{O^\ast(1) + O(\e^\alpha)\norm{(\kappa_\e,W_\e)}}{\vphia(t,z) }  \right).
	\end{multline}
	where the $O(\e)$ depend only on $B$. 
\end{prop}
\begin{proof}[Proof of \cref{prop lin pzzWeps}]$ $\\
	We start, as ever, by differentiating \cref{perturb Weps}, but now three times. This yields for the right hand side ten terms : 	\begin{align}\label{eq pzzzW}
	\p_z^3 M(t,z) \! - \! \e^2 \p_z^3 \p_t \Vea(t,z)  -\e^4 \p_z^3 \p_t W_\e(t,t) \! = \! R_1 + R_2 + R_3 + R_4 + R_5+R_6+R_7+R_8+R_9+R_{10},
	\end{align} 
	with the following notations :
	\begin{align*}
	& R_1  : = \p_z^3  \I_\e(\qea + \e^2 \kappa_\e, \Vea + \e^2 W_\e)(t,z) M(t,z) \exp( \e^2\Xi_\e(t,z)) ,  \\ 
	& R_2 := 3 \p_z^2  \I_\e(\qea + \e^2 \kappa_\e, \Vea + \e^2 W_\e)(t,z)  \p_z M(t,z)  \exp( \e^2\Xi_\e(t,z)), \\
	& R_3 := 3 \e^2  \p_z^2  \I_\e(\qea + \e^2 \kappa_\e, \Vea + \e^2 W_\e)(t,z) M(t,z)  \exp( \e^2\Xi_\e(t,z))  \p_z \Xi_\e(t,z) , \\
	& R_4 :=  6 \e^2 \p_z \I_\e(\qea + \e^2 \kappa_\e, \Vea + \e^2 W_\e)(t,z)  \p_z M(t,z) \exp( \e^2\Xi_\e(t,z))\p_z \Xi_\e(t,z),\\
	& R_5 := 3 \p_z \I_\e(\qea + \e^2 \kappa_\e, \Vea + \e^2 W_\e)(t,z) \p_z^2 M(t,z) \exp( \e^2\Xi_\e(t,z)),	\end{align*}
	and moreover :
	\begin{multline*}
	R_6 := 3 \e^2   \p_z \I_\e(\qea + \e^2 \kappa_\e, \Vea + \e^2 W_\e)(t,z)  M (t,z) \exp( \e^2\Xi_\e(t,z))  \Big( \e^2 \p_z \Xi_\e(t,z)^2 + \p_z^2 \Xi_\e(t,z) \Big) ,
	\end{multline*}
	\begin{multline*}
	R_7 := 3 \e^2 \I_\e(\qea + \e^2 \kappa_\e, \Vea + \e^2 W_\e)(t,z)  \p_z M (t,z) \exp( \e^2\Xi_\e(t,z))   \Big( \e^2  \p_z \Xi_\e(t,z)^2 + \p_z^2 \Xi_\e(t,z)\Big) ,
	\end{multline*}
	\begin{align*}
	& 	R_8 :=  3 \e^2   \I_\e(\qea + \e^2 \kappa_\e, \Vea + \e^2 W_\e)(t,z) \p_z^2 M(t,z)  \exp( \e^2\Xi_\e(t,z))  \p_z \Xi_\e(t,z) , \\
	&	R_9 :=   \I_\e(\qea + \e^2 \kappa_\e, \Vea + \e^2 W_\e)(t,z) \p_z^3 M(t,z)  \exp( \e^2\Xi_\e(t,z)).
	\end{align*}
	The last term corresponds to the third derivative of the exponential term $\exp(\e^2\Xi_\e)	$. 
	\begin{multline*}
	R_{10}:=  \e^2 \I_\e(\qea + \e^2 \kappa_\e, \Vea + \e^2 W_\e)(t,z)  M (t,z) \exp( \e^2\Xi_\e(t,z))  \\ \times \Big( \e^4\p_z \Xi_\e(t,z)^3 + 3\e^2  \p_z \Xi_\e(t,z)  \p_z^2 \Xi_\e(t,z) +  \p_z^3 \Xi_\e(t,z) \Big) .
	\end{multline*}
	We first tackle $R_1$. We use the linearization of the third derivative of $\I_\e$ in \cref{prop lindecay}.
	\begin{align*}
	R_1  = &   M(t,z)   \left(  \p_z^3 \Iea(t,z)+ \frac{\e^2 \norm{\vphia \p_z^3 W_\e}_\infty }{2^{1-\alpha}\vphia(t,z)}  +   \frac{ O(\e^{2+\alpha})\norm{(\kappa_\e,W_\e)}}{\vphia(t,z) } \right) \Big(1 + \e^2\Xi_\e(t,z) + O(\e^4) \norm{(\kappa_\e,W_\e)} \Big) \\
	= &  \e^2M(t,z) \left(  \frac{\norm{\vphia \p_z^3 W_\e}_\infty}{2^{1-\alpha} \vphia(t,z)} + \frac{O^\ast(1) + O(\e^\alpha)\norm{(\kappa_\e,W_\e)}}{\vphia(t,z) } \right) \Big(1  + \e^2\Xi_\e(t,z) + O(\e^4)\norm{(\kappa_\e,W_\e)} \Big). 
	\end{align*}
	We end up with the following estimate
	\begin{align}\label{eq R1c}
	R_1 =  \e^2 M(t,z) \left( \frac{ \norm{\p_z^3 W_\e}_\infty}{2^{1-\alpha}\vphia(t,z)} +  \frac{O^\ast(1) + O(\e^\alpha)\norm{(\kappa_\e,W_\e)}}{\vphia(t,z) }\right).
	\end{align}
	For $R_2$, with \cref{prop lindecay} we have
	\begin{align*}
	R_2 = & 3\p_z M(t,z) \left( \p_z^2 \Iea(t,z) +  \frac{ O(\e^3)\norm{(\kappa_\e,W_\e)}}{\vphia(t,z)} \right) \Big( 1 + \e^2\Xi_\e(t,z) + O(\e^4) \norm{(\kappa_\e,W_\e)} \Big), \\
	= & 3 \p_z M(t,z)  \left( \frac{O^\ast(\e^2) + O(\e^3) \norm{(\kappa_\e,W_\e)}}{\vphia(t,z)} \right) \Big( 1 + \e^2\Xi_\e(t,z)  + O(\e^4) \norm{(\kappa_\e,W_\e)} \Big).
	\end{align*}
	We can simplify this expression to
	\begin{align}\label{eq R2c}
	R_2 = \e^2 \p_z M(t,z)\left( \frac{O^\ast(1) + O(\e)\norm{(\kappa_\e,W_\e)}}{\vphia(t,z) } \right).
	\end{align}
	For $R_3$ we get 
	\begin{align*}
	R_3 = & 3 \e^2 M(t,z) \p_z \Xi_\e(t,z) \left( \p_z^2 \Iea(t,z) +  \frac{ O(\e^3)\norm{(\kappa_\e,W_\e)}}{\vphia(t,z)} \right) \Big( 1 + \e^2\Xi_\e(t,z) + O(\e^4) \norm{(\kappa_\e,W_\e)} \Big), \\
	= & 3 \e^2 M(t,z) \p_z \Xi_\e(t,z) \left( \frac{O^\ast(\e^2) + O(\e^3) \norm{(\kappa_\e,W_\e)}}{\vphia(t,z)} \right) \Big( 1 + \e^2\Xi_\e(t,z)  + O(\e^4) \norm{(\kappa_\e,W_\e)} \Big).
	\end{align*}
	We can simplify roughly this expression to
	\begin{align}\label{eq R3c}
	R_ 3 =   \frac{O(\e^3) \norm{( \kappa_\e ,W_\e) }}{\vphia(t,z)} M(t,z).
	\end{align}
	For $R_4$ one has very similarly 
	\begin{align*}
	R_4 = & 6 \e^2 \p_z M(t,z) \p_z \Xi_\e(t,z) \left( \p_z \Iea(t,z) +  \frac{ O(\e^3)\norm{(\kappa_\e,W_\e)}}{\vphia(t,z)} \right) \Big( 1 + \e^2\Xi_\e(t,z) + O(\e^4) \norm{(\kappa_\e,W_\e)} \Big), \\
	= & 6 \e^2 \p_z M(t,z) \p_z \Xi_\e(t,z) \left( \frac{O^\ast(\e^2) + O(\e^3) \norm{(\kappa_\e,W_\e)}}{\vphia(t,z)} \right) \Big( 1 + \e^2\Xi_\e(t,z)  + O(\e^4) \norm{(\kappa_\e,W_\e)} \Big).
	\end{align*}
	We can simplify this expression to
	\begin{align}\label{eq R4c}
	R_ 4 =   \frac{O(\e^3) \norm{( \kappa_\e ,W_\e) }}{\vphia(t,z)} \p_z M(t,z).
	\end{align}
	The expression for $R_5$ still follows the same road
	\begin{align*}
	R_5 = & 3 \p_z^2  M(t,z) \left( \p_z \Iea(t,z) +  \frac{ O(\e^3)\norm{(\kappa_\e,W_\e)}}{\vphia(t,z)} \right) \Big( 1 + \e^2\Xi_\e(t,z) + O(\e^4) \norm{(\kappa_\e,W_\e)} \Big), \\
	= &  3 \p_z^2 M(t,z) \left( \frac{O^\ast(\e^2) + O(\e^3) \norm{(\kappa_\e,W_\e)}}{\vphia(t,z)} \right) \Big( 1 + \e^2\Xi_\e(t,z)  + O(\e^4) \norm{(\kappa_\e,W_\e)} \Big).
	\end{align*}
	The last expression can be shortened in
	\begin{align}\label{eq R5c}
	R_ 5 =  3 \e^2 \p_z^2 M(t,z) \frac{O^\ast(1) + O(\e) \norm{( \kappa_\e ,W_\e) }}{\vphia(t,z)}.
	\end{align}
	For $R_6$, the expression is a little more involved due to the second derivative of the exponential
	\begin{multline*}
	R_6  =  \e^2 M(t,z)\left( \frac{O^\ast(\e^2)+ O(\e^3)\norm{(\kappa_\e,W_\e)}}{\vphia(t,z)} \right)   \Big(1 +  \e^2\Xi_\e(t,z) +O(\e^4) \norm{(\kappa_\e,W_\e)} \Big) \\ \times
	\left( \frac{O(\e^2) \norm{(\kappa_\e,W_\e)}}{\vphia(t,z)}+ \p_z^2 \Xi_\e(t,z) \right) .
	\end{multline*}
	We eventually shorten $R_6$ as
	\begin{align}\label{eq R6c}
	R_6 = 3 M(t,z)  \frac{ O(\e^3 )\norm{(\kappa_\e,W_\e)}}{\vphia(t,z) } .
	\end{align}
	If we bridge together all of our previous estimates in \cref{eq R1c}, \cref{eq R2c}, \cref{eq R3c}, \cref{eq R4c} and \cref{eq R5c}, \cref{eq R6c} we obtain that 
	\begin{multline}\label{eq sumR16}
	R_1 + R_2 + R_3 + R_4 +R_5+R_6 = \e^2   M(t,z) \left( \frac{O^\ast(1) + O(\e^\alpha )\norm{(\kappa_\e,W_\e)}}{\vphia(t,z) } \right) \\ + \e^2   \p_z M(t,z) \left( \frac{O^\ast(1) + O(\e)\norm{(\kappa_\e,W_\e)}}{\vphia(t,z) } \right) + \e^2   \p_z^2 M(t,z) \left( \frac{O^\ast(1) + O(\e)\norm{(\kappa_\e,W_\e)}}{\vphia(t,z) } \right)  \\ + \frac{\e^2 \norm{\vphia \p_z^3 W_\e}_\infty}{2^{1-\alpha} \vphia(t,z)}M(t,z). 
	\end{multline}
	In that first round of estimates, we have shown that all the contributions of the terms with the derivatives of  $\I_\e$ do not appear when linearizing because they are of high order in $\epsilon$. Therefore, the most meaningful contribution will now appear, because $\I_\e$ now contributes mainly as $1$ and no longer vanishes. 
	
	We start with $R_7$  :
	\begin{multline*}
	R_7 = 3 \e^2 \p_z M (t,z)  \Big( 1 + O^\ast(\e^2)+ O(\e^3)\norm{(\kappa_\e,W_\e)} \Big)  \Big(1 +  \e^2\Xi_\e(t,z) +O(\e^4) \norm{(\kappa_\e,W_\e)} \Big) \\ \times
	\left( \frac{O(\e^2) \norm{(\kappa_\e,W_\e)}}{\vphia(t,z)}+ \p_z^2 \Xi_\e(t,z) \right) ,
	\end{multline*}
	which can be rewritten as
	\begin{align*}
	R_7= 3 \e^2 \p_z M(t,z) \Big( 1 + O^\ast(\e) + O(\e^2) \norm{(\kappa_\e,W_\e)} \Big) \left(  \p_z^2 \Xi_\e(t,z)   + \frac{O(\e^2) \norm{(\kappa_\e,W_\e)}}{\vphia(t,z)} \right).
	\end{align*}
	Finally, for $R_7$ : 
	\begin{align}\label{eq R7c}
	R_7= 3 \e^2 \p_z M(t,z)   \p_z^2 \Xi_\e(t,z) + \p_z M(t,z)\left( \frac{O^\ast(\e^3) + O(\e^4)\norm{(\kappa_\e,W_\e)} }{\vphia(t,z)}\right).
	\end{align}
	For $R_8$, the following estimates hold true, 
	\begin{align*}
	\nonumber R_8 =  3 \e^2 \p_z^2 M(t,z) \p_z \Xi_\e(t,z)\Big( 1 + O^\ast(\e) + O(\e^2) \norm{(\kappa_\e,W_\e)} \Big)\Big(1 +  \e^2\Xi_\e(t,z) +O(\e^4) \norm{(\kappa_\e,W_\e)} \Big) .
	\end{align*}
	Therefore
	\begin{align}\label{eq R8c}
	R_8 = 3 \e^2 \p_z^2 M (t,z) \p_z \Xi_\e(t,z) +  \p_z^2 M(t,z) \left( \frac{O^\ast(\e^3)+  O(\e^4 )\norm{(\kappa_\e,W_\e)}}{\vphia(t,z)}\right) .
	\end{align}
	For the last two terms, the  derivatives up to the third order appear. The simplest is given by  $R_9$ :
	\begin{align}
	\nonumber R_9 & =  \p_z^3 M(t,z)   \Big(1 + O^\ast(\e^2) + O(\e^3)  \norm{( \kappa_\e ,W_\e) }\Big)\Big(1 + \e^2\Xi_\e(t,z) + O(\e^4) \norm{( \kappa_\e ,W_\e) }\Big), \\
	&\label{eq R9c} = \p_z^3 M(t,z) + \e^2 \p_z^3 M(t,z) \Big( \Xi_\e(t,z) + O^\ast(1) + O(\e) \norm{(\kappa_\e,W_\e)} \Big).
	\end{align}
	At last, for the term $R_{10}$, 
	\begin{multline}
	R_{10} = \e^2  M (t,z)  \Big( 1 + O^\ast(\e^2)+ O(\e^3)\norm{(\kappa_\e,W_\e)} \Big)  \Big(1 +  \e^2\Xi_\e(t,z) +O(\e^4) \norm{(\kappa_\e,W_\e)} \Big) \\ \times
	\left( \frac{O(\e^2) \norm{(\kappa_\e,W_\e)}}{\vphia(t,z)}+ \p_z^3 \Xi_\e(t,z)\right) .
	\end{multline}
	It is shortened to 
	\begin{align}\label{eq R10c}
	R_{10}= \e^2 M(t,z) \p_z^3 \Xi_\e(t,z)+ \e^2 M(t,z)  \frac{O(\e^2) \norm{(\kappa_\e,W_\e)}}{\vphia(t,z)} .
	\end{align}
	We now add every estimate, starting from \cref{eq sumR16} and with \cref{eq R7c}, \cref{eq R8c}, \cref{eq R9c} and \cref{eq R10c} to obtain
	\begin{multline}
	\sum_{j=1}^{10}R_j =  \p_z^3 M(t,z)  + \e^2 \p_z^3 M(t,z) \Big( \Xi_\e(t,z) + O^\ast(\e^2) + O(\e^3) \norm{(\kappa_\e,W_\e)} \Big) \\  
	+ 3 \e^2 \p_z^2 M (t,z) \left( \p_z \Xi_\e(t,z)  + \frac{O^\ast(1)  +  O(\e) \norm{(\kappa_\e,W_\e)}}{\vphia(t,z)} \right)  \\ +  3 \e^2 \p_z M(t,z) \left(  \p_z^2 \Xi_\e(t,z)+ \frac{ O^\ast(1) + O(\e) \norm{(\kappa_\e,W_\e)} }{\vphia(t,z)}\right)  \\
	 + \e^2 M(t,z)  \left(  \p_z^3 \Xi_\e(t,z)  + \frac{\norm{\vphia \p_z^3 W_\e}_\infty}{2^{1-\alpha}\vphia(t,z)}    +\frac{O^\ast(1) + O(\e^\alpha )\norm{(\kappa_\e,W_\e)}}{\vphia(t,z) }\right)    .
	\end{multline}
	To conclude the proof, we deal with the left hand side of \cref{eq pzzzW} as in the linearization of the second derivative, noticing that the terms $\p_z^3 M$ cancel on each side.
\end{proof}

\section{Stability of  the linearized equations}\label{sec lin stab}
Building upon the series of linear approximations, we can study the stability of $W_\e$ in the space $\F$. The first result is to control the different terms of $\F$ in the norm $\normf{\cdot}$, see \Cref{def F}. The weight function introduced in the definition of $\E$ is meant  to enable  controlling the behavior at infinity.  
\begin{thm}[Stability analysis]\label{stab We}$ $\\
	For any ball $B$ of $\R \times \F$, there exists a constant $\e_B$ that depends only on $B$ such that for all $\e\leq \e_B$, any pair $(\kappa_\e,W_\e) \in B$  solution of \cref{perturb Weps} verifies the following bounds : 
	%
	\begin{align*}
	\norm{\Xi_\e}_{\infty} \leq &   O^\ast_0(1) + O(\e) \norm{(\kappa_\e,W_\e)},\\
	\norm{\p_z W_\e}_{\infty} \leq &   O^\ast_0(1) + O(\e) \norm{(\kappa_\e,W_\e)},\\
	\norm{\vphia \p_z \Xi_\e}_{\infty} \leq &    O^\ast_0(1) + O(\e) \norm{(\kappa_\e,W_\e)},  \\
	\norm {\vphia \p_z^2 W_\e}_{\infty} \leq &  O_0^\ast(1) + O(\e) \norm{(\kappa_\e,W_\e)},\\
	\norm{\vphia\p_z^3 W_\e}_{\infty} \leq &  O_0^\ast(1) + O(\e^\alpha) \norm{(\kappa_\e,W_\e)}+ k(\alpha) \norm{W_\e}_\F.
	\end{align*}
where $O^\ast_0 (1) = \max\Big( O^\ast(1), O(1)\norm{  W_\e(0,\cdot)}_\F \Big)$, and $k(\alpha)<1$ is a uniform constant.
\end{thm}
The proof of this theorem is quite intricate and will be divided in several subsections. The plan is a follows :
\begin{itemize}
	\item First, we focus on a small ball around $\zs(t)$. The first step is to get bounds only on a small time interval on this ball, and the second step is to propagate this bound uniformly in time, locally in space.
	\item Next, we propagate this bound on the whole space by dividing it in successive dyadic rings $D_n$ centered around $\zs$, see \cref{def Dn}.
\end{itemize}
The main arguments are the maximum principle coupled with a suitable  division of the space that accounts for the non local nature of  the infinitesimal operator. The purpose of this dyadic decomposition in rings is to obtain a decay of the norm with  respect to the radius of the ring.  
\subsection{Division of the space in  a ball surrounded by dyadic rings}$ $\\

Let us first consider a time $T_\ast$. Then for all times such that $ 0 \leq t,s \leq T_\ast$, the inequality 
\begin{align*}
\abs{\zs(t)-\zs(s)} \leq \sup_{s\geq 0} \abs{ m'(\zs(s))} T_{\ast}:=r_\ast
\end{align*}
holds true, and the supremum is finite because $\zs$ lives in a bounded domain uniquely determined by $m$ and $\zs(0)$, see \cref{eqdef zs}.  

We slightly expand this ball by a constant $r_0$ to be defined later and define $$B_0 := \left\{ z \text{ such that}   \abs{z-\zs(0)}  \leq r_0 +r_\ast \right\}.$$ 
Our intention behind this choice is that the ball $B_0$ verifies the following property :
\begin{align}\label{eq R0T0}
\forall t \leq T_\ast, \forall z \in B_0, \quad  \abs{ z- \bar{z}(t)} = \frac{\abs{z-\zs(t)}}{2} = \frac{\abs{z-\zs(0) + \zs(0) - \zs(t)}}{2}\leq \frac{r_0}{2}+r_\ast.
\end{align}
We recall that $\bz(t) := \dfrac{z+\zs(t)}{2}$. We will split the rest of the space around $B_0$ in successive dyadic rings. The first ring is defined as $D_1 = \left\{ z:  r_0+r_\ast \leq \abs{z-\zs(0)} \leq 2 r_0+ r_\ast \right\}$.
It verifies for every $t \leq T_\ast $ the following identity on the middle point :
\begin{align*}
\abs{\bar z -\zs(0)}= \abs{\frac{z+\zs(t)}{2}-\zs(0)} & \leq \abs{\frac{z-\zs(0)}{2}} + \abs{\frac{\zs(0)-\zs(t)}{2}}, \\
& \leq r_0 + r_\ast. 
\end{align*}
This shows that for any $z \in D_1$, and time $t \leq T_{\ast}$, any  middle point $\bz(t)$  lies in $B_0$. More generally, the following lemma holds true if we define for $n\geq 2$ : 
\begin{align}\label{def Dn}
D_n := \left\{  2^{n-1}r_0+r_\ast \leq \abs{z-\zs(0)} \leq 2 ^n r_0+r_\ast \right\},
\end{align} 
\begin{lem}[Middle point property]\label{middle point}$ $\\
	For every time $ 0 \leq t \leq T_\ast$ : 
	\begin{align*}
	\forall n \geq 1  \ \forall z \in D_n, \quad  \bz(t) \in D_{n-1},
	\end{align*}
with the convention $D_0=B_0$.
\end{lem}
\begin{figure}
	\begin{center}
		\begin{tikzpicture}[scale=0.45]
		\draw[red!80!black, thick] plot [smooth] coordinates {(-1.5,0) (0,-1.3) (0.2,0.8)(1.5,0)};
		
		\draw[] (0,0) circle (1.5);
		\draw[thick, orange] (0,0) circle (2);
		\draw[] (0,0) circle (3.5);
		\draw[] (0,0) circle (7.5);
		\draw[orange] (0,2) node[above] {$B_0$};
		\draw[red] (0,0) node[left ] {$\zs$};
		\draw[] (0,3.5) node[above] {$D_1$};
		\draw[] (0,7.5) node[above] {$D_2$};
		\end{tikzpicture}
	\end{center}
	\caption{Illustration of the division of the space in successive dyadic rings.}\label{fig rings}
\end{figure}
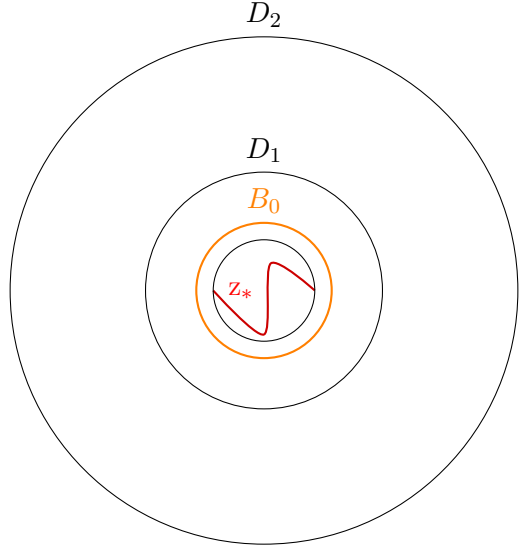
Moreover, the following inequalities are a direct consequence of the definition of $D_n$ and $T^\ast$ :
\begin{align}\label{center zs}
\forall t \leq T_\ast, \forall z \in D_n, \quad  2^{n-1} r_0 \leq \abs{z-\zs(t)} \leq 2^n r_0 + 2r_\ast. 
\end{align}
\paragraph*{\textbf{Notations for this section :}} We will denote $\norm{\cdot}_\infty^n$ the $L^\infty$ norm on $\R_+\times D_n$.
\subsection{Local bounds on $B_0$}$ $\\
Our first step consists in getting  bounds on the ball $B_0$, uniformly in time. 
\begin{prop}[Local bounds]\label{bound B0}$ $\\
For a convenient choice of $T^\ast$ and $r_0$ introduced above, and  made explicit in \cref{choice R0T0}, there exists a constant $\e_B$ that depends only on $B$, such that upon the conditions of \cref{stab We}, $W_\e$ verifies for $\e\leq \e_B$
	\begin{align*}
\norm{\Xi_\e}_{\infty}^0 \leq &   O^\ast_0(1) + O(\e) \norm{(\kappa_\e,W_\e)},\\
\norm{\p_z W_\e}_{\infty}^0 \leq &   O^\ast_0(1) + O(\e) \norm{(\kappa_\e,W_\e)},\\
\norm{\vphia \p_z \Xi_\e}_{\infty}^0 \leq &    O^\ast_0(1) + O(\e) \norm{(\kappa_\e,W_\e)},  \\
\norm {\vphia \p_z^2 W_\e}_{\infty}^0 \leq &  O_0^\ast(1) + O(\e) \norm{(\kappa_\e,W_\e)},
\end{align*}
	where $\ds \norm{W}_{\infty}^0:= \sup_{(t,z)\in \R_+\times B_0} \abs{W(t,z)}$ and $O^\ast_0 (1) = \max\Big( O^\ast(1), O^\ast(1)\norm{  W_\e(0,\cdot)}_\F \Big)$.
\end{prop}
To prove this "local" bound, \textit{i.e.} in the ball $B_0$, one must start with the higher order derivative to build a contraction argument. Estimates of the lower order derivatives are then successively deduced by integration. Clearly, our argument for the third derivative is  the more technical because it involves a lot of terms through the linearized approximation made in \cref{prop lin pzzzWeps}. Therefore, for clarity reason, third derivatives are left out from \cref{bound B0}, we will deal with them, locally and on the rings, in \cref{rings pzzzWe}. We present here our argument on the simpler derivatives up to order two, and we refer to \cref{sec pz3 stab} for the generalization of the method to the third derivative.

 Interestingly, to prove the non local estimates on the rings, we will proceed in the reverse way by first dealing with the lower order derivatives.
\begin{proof}[Proof of \cref{bound B0}]$ $\\
By the derivation of the linearized equation in \cref{prop lin pzzWeps}, $W_\e$ verifies, see \cref{def M}:
	\begin{multline*}
	\e^2 \p_t \p_z^2  W_\e(t,z) =   - \p_z^2 M(t,z) \Big(  W_\e(t,z)- 2 W_\e \left(t, \bz \right)  + O^\ast(1)+ O(\e)\norm{(\kappa_\e,W_\e)}\Big)  \\ - 2   \p_z M(t,z) \Big[  \p_zW_\e(t,z)-\p_z W_\e(t,\bz )  + O^\ast(1)+ O(\e)\norm{(\kappa_\e,W_\e)} \Big] \\ + M(t,z)  \left(  \frac12 \p_z^2 W_\e(t,\bz )  -\p_z^2 W_\e(t,z)+ O^\ast(1)+ O(\e)\norm{(\kappa_\e,W_\e)} \right).
	\end{multline*}
	We will  use the maximum principle on the ball $B_0$.  The key point is that on this ball, all other factors are controlled by $ \norm{\p_z^2 W_\e}_\infty$. To compare all those terms with $\p_z^2 W_\e$, we perform Taylor expansions with respect to the space variable. First, we write that for any $z \in B_0$, thanks to \cref{eq R0T0},
	\begin{align*}
	\p_z W_\e(t,\bz) -\p_z W_\e(t,z) \leq \left( \dfrac{r_0}{2}+r_\ast \right) \norm{ \p_z^2 W_\e(t,\cdot) }_{L^\infty\left( B_0\right)}.
	\end{align*}
	Similarly, there exists $\xi\in \left( z,\bz \right)$ and $\xi' \in \left(\zs , \bz \right)$  such that 
	\begin{multline}\label{eq expMe}
\Xi_\e(t,z) =	W_\e(t,z) -2 W_\e(t,\bz)  + W_\e(t,\zs)  = \left( \frac{z-\zs}{2} \right)\p_z W(t,\bar{z}) + \frac12 \dfrac{(z-\zs)^2}{4} \p_z^2 W(t,\xi)  \\  - \left( \frac{z-\zs}{2} \right)\p_z W(\bar{z}) +  \frac12 \dfrac{(z-\zs)^2}{4} \p_z^2 W(\xi')  \leq \frac14 \left( \dfrac{r_0}{2} + r_\ast
	\right)^2    \norm{ \p_z^2 W_\e(t,\cdot)}_{L^\infty\left(  B_0\right)}. \
	\end{multline}
	Moreover by the hypothesis made in \cref{decay Gamma} on $M$, for $j=1,2$
	\begin{align*}
	\sup_{(t,z)\in \R_+ \times B_0}\abs{\frac{\p_z^{(j)} M(t,z)}{M(t,z)}} \leq O^\ast(1).
	\end{align*}
	Thanks to those a priori bounds, when we evaluate \cref{PL pzzWeps+} at the  point of maximum of $\p_z^2 W_\e$ on $B_0$ we get 
	\begin{multline*}
	\e^2 \p_t 	\left(\norm{ \p_z^2 W_\e(t,\cdot)}_{L^\infty \left( B_0\right)} \right) \leq  M(t,z) \left[  \frac{1}{2} \norm{\p_z^2 W_\e(t,\cdot)}_{L^\infty\left( B_0\right)}  -  \norm{ \p_z^2 W_\e(t,\cdot)}_{L^\infty\left( B_0\right)} \right.  \\    + O^\ast(1)  \left(   \frac14\left( \dfrac{r_0}{2}+r_\ast \right)^2 +  \frac{r_0}{2}  +  r_\ast  \right)  \norm{\p_z^2 W_\e(t,\cdot)}_{L^\infty\left( B_0\right)}    
	+  O^\ast(1) +  O(\e )  \norm{(\kappa_\e,W_\e)} \Big].
	\end{multline*} 
	The \emph{crucial} step is that we  choose  $T^\ast$ and $r^\ast$ so  small so that
	\begin{align}\label{choice R0T0}
	O^\ast(1) \left(   \frac14\left( \dfrac{r_0}{2}+r_\ast \right)^2 + \frac{r_0}{2}+ r_\ast \right) \leq \frac14.
	\end{align} 
	The consequence is that
	\begin{align*}
	\e^2 \p_t 	\left( \norm{  \p_z^2 W_\e(t,\cdot)}_{L^\infty\left( B_0\right)} \right) \leq M(t,z) \left[ -\frac{1}{4}  \norm{ \p_z^2 W_\e(t,\cdot)}_{L^\infty\left( B_0\right)} + O^\ast(1) + O(\e) \norm{(\kappa_\e,W_\e)} \right] .
	\end{align*}
The function $M(t,z)$ admits a lower bound. Therefore, we can apply the maximum principle, on the ball $B_0$ :
	\begin{align*}
	\norm{  \p_z^2 W_\e}_{L^\infty\left(\left[0,T^\ast\right] \times B_0\right)} \leq \max \Big(  O^\ast(1) + O(\e) \norm{(\kappa_\e,W_\e)}, \norm{ \p_z^2 W_\e(0,\cdot)}_{L^\infty\left(	B_0\right)} \Big).	\end{align*}
	We now detail how to propagate this bound uniformly in time. One can  renew every previous estimate on each interval  $I_k:=\left[kT_\ast,(k+1)T_\ast\right]$. By going over the same steps, we notice  that the only argument that changes for different $k$ is the center of the ball $B_0$ around $\zs$, but interestingly not its radius see \cref{choice R0T0}. Every other estimate is the same and is independent of $k$. Therefore, since the condition \cref{choice R0T0} is uniform in time ($O^\ast(1)$ does not depend on time), once the radius is chosen small enough depending only on $K^\ast$, see \cref{choice R0T0}, we can repeat recursively the estimates on each interval $I_k$. Considering all $k \in \N$, we have therefore proven that 
	\begin{align}
\nonumber 	\norm{\p_z^2 W_\e}_\infty^0 & \leq  \max \Big(  O^\ast(1) + O(\e) \norm{(\kappa_\e,W_\e)}, \norm{ \p_z^2 W_\e(0,\cdot)}_{L^\infty(B_0)}\Big),\\
& \label{eqd2W} \leq 	 O_0^\ast(1) + O(\e) \norm{(\kappa_\e,W_\e)}.
	\end{align}
We will use this estimate  as the starting point in order to prove the rest of \cref{bound B0}. First, notice that adding the weight function $\vphia$ is straightforward, since it is uniformly bounded on $B_0$ :
	\begin{align*}
\norm{ \vphia \p_z^2 W_\e}_{\infty}^0 \leq  O_0^\ast(1) + O(\e) \norm{(\kappa_\e,W_\e)}.
\end{align*}
Next, taking advantage that both $W_\e$ and $\p_z W_\e$ vanish at $z^\ast$, we write 
	\begin{align*}
	\abs{\p_z W_\e(t,z)} =\abs{\int_{\zs(t)}^z \p_z^2 W_\e(t,z') dz'} \leq  (r_0+2r_\ast)\norm{  \p_z^2 W_\e}_\infty^0.
	\end{align*}
	As a consequence, using again the expansion of \cref{eq expMe}, 
	\begin{align*}
	\abs{\Xi_\e(t,z)} = \abs{2 W_\e(t,\bz(t)) -W_\e(t,z)} \leq 
	 \frac14 \left(	\frac{r_0}{2}+r_\ast \right)^2 \norm{  \p_z^2 W_\e}_\infty^0.
	\end{align*}
Similarly, we get a uniform bound on $\p_z \Xi_\e$.  
Combining those estimates with the first one in \cref{eqd2W}, that comes from the maximum principle, the proof of \cref{bound B0} is concluded.
\end{proof}
\subsection{Bound in the rings, $\Xi_\e$}\label{sec xie}$ $\\
We will now propagate those bounds beyond the small ball. It is very important to keep the level of precision of $O^\ast(1) + O(\e)\norm{(\kappa_\e,W_\e)}$, to which we will add some decay property due to the increasing size of the rings. 
\begin{prop}[In the rings, $\Xi_\e$]\label{rings Me}$ $\\
	There exists a constant $\e_B$ that depends only on $B$ such that upon the conditions of \cref{stab We}, $W_\e$ verifies for $\e\leq \e_B$
	\begin{align}\label{bd_W1_W0}
	\norm{\Xi_\e }_{\infty}^{n} & \leq  O_0^\ast(1) +  O(\e) \norm{(\kappa_\e,W_\e)},
	\end{align}
	for all $n \geq 1$.
\end{prop}
\begin{proof}[Proof of \cref{rings Me}]$ $\\
	For any $n \geq 1$, take  $z$ in the ring $D_{n}$ defined previously. Then, $\bz \in D_{n-1}$ by \cref{middle point}.  Next, we use the linearized equation given by \cref{prop lin M}. For $t \in \R_+$ and $z \in D_{n}$ the following inequality holds true
	\begin{align*}
	\e^2 \p_t \Xi_\e(t,z) \leq  M(t,z) \left( 2 \frac{ M(t, \bz)}{M(t,z)} \norm{\Xi_\e}_\infty^{n-1} - \Xi_\e(t,z)  + O^\ast(1) + O(\e) \norm{(\kappa_\e,W_\e)} \right) 
	\end{align*}
We define $a_n$ such that the quotient of $M$ verifies :
	\begin{align*}
	\sup_{(t,z) \in \R_+\times D_{n} } \abs{\frac{ M(t, \bz)}{M(t,z)}} := a_n,
	\end{align*}
	where the sequence $a_n$ is bounded and verifies $a_n \to a<\frac12$ as $ n \to \infty$ by the hypothesis made in \cref{cond m infty}.

	Moreover since $M$ admits a uniform lower bound by \cref{cond Gamma}, we can apply the maximum principle:
	\begin{align}\label{bd_Wn1_Wn}
	\norm{\Xi_\e }_{\infty}^{n} \leq \max \Big( 2 a_n \norm{ \Xi_\e }_{\infty}^{n-1}  + O^\ast(1) +  O(\e) \norm{(\kappa_\e,W_\e)}, \norm{\Xi_\e(0,\cdot)}_{L^\infty(D_n)} \Big).
	\end{align}
The interplay between the recursion and the $\max$ in the formula above requires a careful argument.  We first notice that for all $n \in \N$:
\begin{align*}
\norm{\Xi_\e(0,\cdot)}_{L^\infty(D_n)} \leq O_0^\ast(1).
\end{align*}
Therefore, from \cref{bd_Wn1_Wn},
\begin{align}\label{bd Wn}
\norm{\Xi_\e }_{\infty}^{n}   \leq 2 a_n \norm{ \Xi_\e }_{\infty}^{n-1}  + O_0^\ast(1) +  O(\e) \norm{(\kappa_\e,W_\e)}.
\end{align}
Here lies the motivation behind the introduction of the notation $O_0^\ast(1)$. It allows to take into account the initial data and to make recursive estimates that were \textit{a priori} not possible with \cref{bd_Wn1_Wn}.  

Since $2 a_n \to 2a<1$ when $n \to \infty$, we know from \cref{bd Wn} that the sequence $\left(\norm{\Xi_\e }_{\infty}^{n}\right)_n$ is a contraction, with, for instance, a factor $\theta =a + \frac12$, such that $2a<\theta<1$. Since $2 a_n \leq \theta$ but for a finite number of terms, we deduce
\begin{align*}
\norm{ \Xi_\e }_{\infty}^{n} & \leq \max \left( \frac{ O_0^\ast(1) +  O(\e) \norm{(\kappa_\e,W_\e)} }{1-\theta}, \norm{\Xi_\e }_{\infty}^{0} \right), \\
& \leq O_0^\ast(1) +  O(\e) \norm{(\kappa_\e,W_\e)}.
\end{align*}
\end{proof}

\subsection{Bound on the rings : $\p_z \Xi_\e$}$ $\\
We now state a similar result  for $\p_z \Xi_\e$. We see the appearance of the weight function $\vphia$ in the estimates. It slightly worsen the expressions but the methodology is the same than the one deployed to prove \cref{rings Me}. 
\begin{prop}[In the rings, $\p_z \Xi_\e$]\label{rings Ne}$ $\\
	There exists  a constant $\e_B$ that depends only on $B$ such that upon the condition of \cref{stab We}, $W_\e$ verifies for $\e\leq \e_B$
	\begin{align*}
	\norm{\vphia \p_z \Xi_\e }_{\infty}^{n} & \leq O_0^\ast(1) +  O(\e) \norm{(\kappa_\e,W_\e)},
	\end{align*}
	for $n \geq 1$.
\end{prop}

\begin{proof}[Proof of \cref{rings Ne}]$ $\\
	The proof is similar to the bound on $\Xi_\e$, but we have to take the weight function into account. We first make the following computation:
	\begin{align*}
\p_t (\vphia \p_z \Xi_\e)(t,z) = \vphia(t,z) \p_t \p_z \Xi_\e(t,z) + \p_z \Xi_\e(t,z) \p_t \vphia(t,z). 
	\end{align*}  
First,
\begin{align*}
\p_z \Xi_\e(t,z) \p_t \vphia(t,z) =  \alpha \p_z \Xi_\e (t,z) \frac{m'(\zs) \sign(z-\zs)}{\Big(1 + \abs{z-\zs}\Big)^{1-\alpha} } = O^\ast(1) \p_z \Xi_\e(t,z),
\end{align*}
and therefore, 
\begin{align}\label{eq:pzXie3}
\e^2 \p_z \Xi_\e(t,z) \p_t \vphia(t,z) =  O^\ast (\e^2) \norm{(\kappa_\e,W_\e)}.
\end{align}
Second, we gave an linear equation verified by $\p_t \p_z \Xi_\e$ in  the \Cref{prop lin N}. With those two ingredients, we find that  for $z\in D_{n}$ and $t \in \R_+ $:  
\begin{multline}\label{eq:pzXi}
\e^2 \p_t  \p_z \Xi_\e(t,z) \leq M(t,z) \left[ a_n \p_z \Xi_\e(t,\bar z) -\p_z \Xi_\e(t,z) +   \frac{O^\ast(1)+O(\e)\norm{(\kappa_\e,W_\e)}}{\vphia(t,z)} \right. \\ +  \left. \frac{O^\ast(1)}{\vphia(t,z)} \left(b_n \norm{\Xi_\e}_\infty^{n-1} + \norm{\Xi_\e}_\infty^{n} + O^\ast(1) + O(\e)\norm{(\kappa_\e,W_\e)} \right) \right],
\end{multline} 
with the following notations: 
	\begin{align*}
	\sup_{(t,z) \in \R_+ \times D_{n} } \abs{\frac{M(t,\bz)}{ M(t,z) }} := a_n, \quad \quad \sup_{(t,z) \in \R_+ \times D_{n} } \abs{\frac{\p_z M(t,\bz)}{\p_z M(t,z)} } := b_n.
	\end{align*}
We  used that, thanks to \cref{decay Gamma}:
	\begin{align*}
	\sup_{(t,z) \in \R_+ \times \R } \left( \vphia(t,z) \abs{\frac{\p_z M(t,z)}{M(t,z)} } \right) \leq O^\ast(1).
	\end{align*}
Coming back to \cref{eq:pzXi}, we first know thanks to our assumption made in  \cref{cond m infty}, that the sequence $b_n$ is uniformly bounded. Moreover, thanks to \cref{rings Me} we can estimate the terms involving $\Xi_\e$ on the rings. Therefore, by multiplying \cref{eq:pzXi} by $\vphia$, one gets, with \cref{eq:pzXie3}:
	 \begin{multline}\label{eq:pzXie2}
	 \e^2 \p_t \Big[ \vphia \p_z \Xi_\e \Big](t,z) \leq  M(t,z)  \left[a_n \abs{\frac{\vphia(t,z )}{\vphia(t,\bar z)}} \norm{\p_z \Xi_\e}_\infty^{n-1} + O_0^\ast(1)+O(\e)\norm{(\kappa_\e,W_\e)} \right. \\   \left.  -  \vphia(t,z) \p_z \Xi_\e(t,z)  + \frac{O^\ast(\e^2) \norm{(\kappa_\e,W_\e)}}{M(t,z)} \right] .
	 \end{multline}
The weight function was chosen  precisely to satisfy the following scaling estimate:
\begin{align}\label{scaling estimate}
\sup_{\R_+ \times \R}\abs{ \frac{\vphia(t,z)}{\vphia(t,\bz)}} \leq 2^\alpha.
\end{align}
The function $1/M$ has a uniform upper bound. Therefore, thanks again to the maximum principle on the equation \cref{eq:pzXie2} we get 
	\begin{align*}
	\norm{\vphia \p_z \Xi_\e}_\infty^{n} \leq \max  \left( 2^\alpha a_n \norm{\p_z \Xi_\e}_\infty^{n-1} +   O_0^\ast(1)+O(\e)\norm{(\kappa_\e,W_\e)}
	, \norm{ \vphia \p_z \Xi_\e(0,\cdot)}_\infty^n \right).
	\end{align*}
To deduce any result by recursion, we proceed  as in the previous proof. Notice that for all $n\in \N$, 
\begin{align*}
\norm{\vphia \p_z \Xi_\e(0,\cdot)}_{L^\infty(D_n)} \leq \norm{W_\e(0,\cdot)}_\F \leq O_0^\ast(1).
\end{align*}
 Therefore,
 	\begin{align*}
 \norm{\vphia \p_z \Xi_\e}_\infty^{n} \leq  2^\alpha a_n \norm{\p_z \Xi_\e}_\infty^{n-1} +   O_0^\ast(1)+O(\e)\norm{(\kappa_\e,W_\e)}. \end{align*}
As before, by hypothesis, $2^\alpha a_n \to 2^\alpha a<1$ when $n \to \infty$,  and therefore, but for a finite number of terms, $2^\alpha a_n \leq 2^\alpha a<1$. We deduce that the sequence $\left(\norm{\vphia \p_z \Xi_\e }_{\infty}^{n}\right)_n$ is a contraction, with, for instance, a factor $\theta =a + \frac12<1$. 
Therefore, using the initialization on the small ball $B_0$ made in \cref{bound B0}:
\begin{align*}
\norm{\vphia \p_z \Xi_\e }_{\infty}^{n} & \leq \max \left( \frac{ O_0^\ast(1) + O(\e) \norm{(\kappa_\e,W_\e)} }{1-\theta}, \norm{\vphia \p_z \Xi_\e }_{\infty}^{0} \right), \\
& \leq O_0^\ast(1) +O(\e) \norm{(\kappa_\e,W_\e)}.
\end{align*}

\end{proof}
\subsection{Bound on the rings : $\p_z^2 W_\e$}$ $\\
We now make a similar statement upon the second derivative. 
\begin{prop}[In the rings, $\p_z^2 W_\e$]\label{rings pzzWe}$ $\\
	There exists a constant $\e_B$ that depends only on $B$ such that upon the condition of \cref{stab We}, $W_\e$ verifies for $\e\leq \e_B$
	\begin{align*}
	\norm{\vphia \p_z^2 W_\e}_{\infty}^{n} & \leq  O_0^\ast(1) +  O(\e) \norm{(\kappa_\e,W_\e)},	\end{align*}
	for $n \geq 1$.
\end{prop}

\begin{proof}[Proof of \cref{rings pzzWe}]$ $\\
We proceed as in the proof of \cref{rings Ne}. 	We already know a linearized approximation for $\p_z^2 W_\e$, thanks to \cref{PL pzzWeps+}.  Taking this into account, one finds that  $ \vphia \p_z^2 W_\e$ solves :
	\begin{multline}\label{eq pzzrings}
	\e^2 \p_t\Big[ \vphia(t,z) \p_z^2  W_\e \Big](t,z) =    -\p_z^2 M(t,z) \vphia(t,z) \Big( \Xi_\e \left(t, z \right) + O^\ast(1)+ O(\e)\norm{(\kappa_\e,W_\e)} \Big)  \\ - 2   \p_z M(t,z)  \Big[\vphia(t,z) \p_z \Xi_\e(t,z) + O^\ast(1) + O(\e)\norm{(\kappa_\e,W_\e)} \Big]  \\+  M(t,z)    \left( \frac{\vphia(t,z)}2 \p_z^2 W_\e(t,\bz ) - \vphia(t,z) \p_z^2 W_\e(t,z) + O^\ast(1)+ O(\e)\norm{(\kappa_\e,W_\e)}  \right) \\ +  O(\e^2)\norm{(\kappa_\e,W_\e)}.
	\end{multline}
	The last term comes from the same computation of $\p_t \vphia$ as the one made in \cref{eq:pzXie3}. We can estimate on the rings most of the terms involved in \cref{eq pzzrings}. First, 	we dispose of the following uniform controls on the ring by \cref{decay Gamma}:
	\begin{align*}
\sup_{(t,z) \in \R_+ \times \R } \left( \vphia(t,z) \abs{\frac{\p_z^2 M(t,z)}{M(t,z)} } \right) \leq O^\ast(1), \quad \sup_{(t,z) \in \R_+ \times \R } \abs{\frac{\p_z M(t,z)}{M(t,z)} }  \leq O^\ast(1).
\end{align*}
We also need the scaling estimate of the weight function, stated in \cref{scaling estimate}. Then, we can bound the right hand side of \cref{eq pzzrings} after factorizing by $M$, for $t \in \R_+$ and $z \in D_{n}$: 
	\begin{multline*}
	\e^2 \p_t\Big[  \vphia(t,z) \p_z^2 W_\e\Big](t,z)  \leq M(t,z)  \left[   -\vphia(t,z) \p_z^2 W_\e(t,z) +\dfrac{1}{2^{1-\alpha}} \norm{ \vphia \p_z^2 W_\e   }_{\infty}^{n-1}  \right. \\ +  O^\ast(1)+ O(\e)\norm{(\kappa_\e,W_\e)}      + O^\ast(1) \Big( \norm{\Xi_\e}_{\infty}^{n} + O^\ast(1)+ O(\e)\norm{(\kappa_\e,W_\e)} \Big)  \\ 
	+ \left.  O^\ast(1) \Big(\norm{ \vphia \p_z \Xi_\e}_{\infty}^{n}  +  O^\ast(1)+ O(\e)\norm{(\kappa_\e,W_\e)} \Big)  \right]  .
	\end{multline*}

	We also control $\Xi_\e$ and $\p_z \Xi_\e$ on the rings thanks to \cref{rings Me,rings Ne}. We therefore can write our last bound as
	\begin{multline*}
	\e^2 \p_t\Big[  \vphia(t,z) \p_z^2 W_\e \Big](t,z)  \leq M(t,z) \Big[   -\vphia(t,z) \p_z^2 W_\e(t,z) \\ \left.+\dfrac{1}{2^{1-\alpha}} \norm{ \vphia \p_z^2 W_\e   }_{\infty}^{n-1}  +  O_0^\ast(1)+ O(\e)\norm{(\kappa_\e,W_\e)} \right].
	\end{multline*}
The function $M(t,z)$ admits a positive lower bound by \cref{cond Gamma}. We can apply the maximum principle:
	\begin{align*}
	\norm{ \vphia \p_z^2W_\e}_{\infty}^{n}  \leq  \max \left( \dfrac{1}{2^{1-\alpha}} \norm{ \vphia \p_z^2W_\e }_{\infty}^{n-1} + O_0^\ast(1)+O(\e) \norm{(\kappa_\e,W_\e)} ,\norm{\vphia \p_z^2 W_\e(0,\cdot)}_\infty^n \right) .
	\end{align*}
The recursive arguments are somehow a little easier in that case compared to the proofs of \cref{rings Me,rings Ne} since the geometric term, $2^{\alpha-1}$, does not depend on $n$. However, first, as earlier, we get rid of the maximum before any recursion, by stating that for all $n\in \N$, 
\begin{align}
	\norm{\vphia \p_z^2 W_\e(0,\cdot)}_\infty^n \leq \norm{W_\e(0,\cdot)}_\F\leq O_0^\ast(1).
\end{align}
Then, 
	\begin{align*}
\norm{ \vphia \p_z^2W_\e}_{\infty}^{n}  \leq  \dfrac{1}{2^{1-\alpha}} \norm{ \vphia \p_z^2W_\e }_{\infty}^{n-1} + O_0^\ast(1)+O(\e) \norm{(\kappa_\e,W_\e)}.
\end{align*}
 Therefore, straightforwardly, we get, because $2^{\alpha-1}<1$ :
\begin{align*}
\norm{\vphia \p_z^2 W_\e }_{\infty}^{n} & \leq \max \left( \frac{ O_0^\ast(1) +  O(\e) \norm{(\kappa_\e,W_\e)} }{1-2^{\alpha-1}}, \norm{\vphia \p_z^2 W_\e }_{\infty}^{0} \right), \\
& \leq O_0^\ast(1) +  O(\e) \norm{(\kappa_\e,W_\e)}.
\end{align*}
\end{proof}
\subsection{Local and on the rings bound for  $\p_z^3 W_\e$}\label{sec pz3 stab}$ $\\
We dedicate this section to the study of $\p_z^3 W_\e$ since it does not exactly fits the mold of the previous estimates due to the additional  factor $\frac1{2^{1-\alpha}}\norm{\vphia \p_z^3 W_\e }_\infty$ in the linearized equation in \cref{prop lin pzzzWeps}.

$\triangleright$ We highlight the difference by first proving the initial bound on the small $B_0$. We  write the  linear equation solved by $\vphia \p_z^3 W_\e$:
\begin{multline}\label{eq lin3}
	-\e^2 \p_t \Big[\vphia \p_z^3 W_\e\Big](t,z) =  \vphia(t,z)\p_z^3 M(t,z) \Big( \Xi_\e(t,z) + O^\ast(1) + O(\e) \norm{(\kappa_\e,W_\e)} \Big)  \\ + 3  \p_z^2 M (t,z)   \Big( \vphia(t,z)\p_z \Xi_\e(t,z) + O^\ast(1)  +  O(\e) \norm{(\kappa_\e,W_\e)} \Big) \\  +  	 3  \p_z M(t,z)   \Big( \vphia(t,z) \p_z^2 \Xi_\e(t,z) +  O^\ast(1) + O(\e) \norm{(\kappa_\e,W_\e)} \Big)  \\
+ 	  M(t,z)   \left(  \vphia(t,z) \p_z^3 \Xi_\e(t,z)   +  \frac{\norm{\vphia \p_z^3 W_\e}_\infty}{2^{1-\alpha}} + O^\ast(1) + O(\e^\alpha)\norm{(\kappa_\e,W_\e)}  \right)\\ - \e^2 \p_z^3W_\e(t,z) \p_t\vphia(t,z).
\end{multline}
Straightforwradly, one finds
\begin{align*}
\e^2 \p_z^3 \Xi_\e(t,z) \p_t \vphia(t,z) =  O^\ast (\e^2) \norm{(\kappa_\e,W_\e)}.
\end{align*}
 We recall that  $\Xi_\e, \p_z \Xi_\e $ and $\p_z^2 \Xi_\e$ are all uniformly bounded on $B_0$, with the weight, by \cref{bound B0}. 
Moreover, from \cref{cond Gamma}, for $j=1,2$
	\begin{align}\label{eq pz3lin}
\sup_{(t,z) \in \R_+ \times \R } \abs{\frac{\p_z^{(j)} M(t,z)}{M(t,z)} }  \leq O^\ast(1), \quad \sup_{(t,z) \in \R_+ \times \R } \left( \vphia(t,z) \abs{\frac{\p_z^3 M(t,z)}{M(t,z)} } \right) \leq O^\ast(1).
\end{align}
Finally, 
\begin{align*}
	 \vphia(t,z) \p_z^3 W_\e(t,\bz) \leq  \frac{2^{\alpha}}{4} \abs{ \vphia(t,\bar z) \p_z^3 W_\e(t,\bar z)}.
\end{align*}
When plugging all of this into \cref{eq lin3}, we obtain, by evaluating at the point of maximum on $B_0$,
\begin{multline*}
\e^2 \p_t  \norm{\vphia(t,\cdot)\p_z^3 W_\e(t,\cdot)}_{L^\infty(B_0)} \leq   M(t,z)    \left[ -\norm{\vphia(t,\cdot)\p_z^3 W_\e(t,\cdot)}_{L^\infty(B_0)} + \frac1{2^{2-\alpha}} \norm{\vphia(t,\cdot)\p_z^3 W_\e(t,\cdot)}_{L^\infty(B_0)} \right.\\
\left. +  \frac{\norm{\vphia \p_z^3 W_\e}_\infty}{2^{1-\alpha}} +  O_0^\ast(1) + O(\e^\alpha)\norm{(\kappa_\e,W_\e)} \right].
\end{multline*}
Since there is a positive lower bound of $M$,  we recognize a contraction argument on the ball $B_0$, and for bounded times $0 < t \leq T^\ast$:
\begin{multline*}
\norm{  \vphia \p_z^3 W_\e}_{L^\infty( \left[0, T^\ast \right] \times B_0)}  \leq  \\
\max \left( \left(\frac{1}{1-2^{\alpha-2}}\right)\left[  O_0^\ast(1) + O(\e^\alpha) \norm{(\kappa_\e,W_\e)}  + \frac1{2^{1-\alpha}} \norm{\vphia \p_z^3 W_\e}_\infty \right] , \norm{\vphia(0,\cdot)\p_z^3 W_\e(0,\cdot)}_{L^\infty(B_0)} \right).
\end{multline*}
Therefore, since the initial data is conctrolled by $O_0^\ast(1)$, we may write :
\begin{align*}
\norm{  \vphia \p_z^3 W_\e}_{L^\infty( \left[0, T^\ast \right] \times B_0)}  \leq  O_0^\ast(1) + O(\e^\alpha) \norm{(\kappa_\e,W_\e)}  + \frac{2^{\alpha-1}}{1-2^{\alpha-2}} \norm{\vphia \p_z^3 W_\e}_\infty.
\end{align*}

As explained before, we can now repeat the procedure on each interval of time $I_k:=\left[kT_\ast,(k+1)T_\ast\right]$ and end up with a bound uniform in time on the ball $B_0$:
\begin{align}\label{bound pzzzW B0}
\norm{\vphia  \p_z^3 W_\e}_\infty^0 \leq  O_0^\ast(1) + O(\e^\alpha) \norm{(\kappa_\e,W_\e)}  + \frac{2^{\alpha-1}}{1-2^{\alpha-2}} \norm{\vphia \p_z^3 W_\e}_\infty.
\end{align} 

$\triangleright$ We now proceed to propagate this bound on the rings, starting again from \cref{eq lin3} and using the maximum principle. For any $t \in \R_+$ and $z \in D_{n}$, we have:
\begin{multline*}
\e^2 \p_t \Big[ \vphia \p_z^3 W_\e\Big] (t,z) \leq  M(t,z)   \left[   -\vphia(t,z)\p_z^3 W_\e(t,z) + \frac 1{2^{2-\alpha}} \norm{\p_z^3 W_\e }_\infty^{n-1} +  \frac{1}{2^{1-\alpha}} \norm{\vphia \p_z^3 W_\e}_\infty  \right.     \\ + O_0^\ast(1) + O(\e^\alpha)\norm{(\kappa_\e,W_\e)}   +   \abs{\vphia(t,z) \frac{\p_z^3 M(t,z)}{M(t,z)}} \Big( \norm{\Xi_\e}^{n}_\infty + O^\ast(1) + O(\e) \norm{(\kappa_\e,W_\e)} \Big)  \\  + \left.
3 \abs{\frac{\p_z^2 M (t,z)}{M(t,z)}} \Big( \norm{ \vphia \p_z \Xi_\e}_\infty^{n} +O^\ast(1)  +  O(\e) \norm{(\kappa_\e,W_\e)} \Big)  \right. \\ \left.  +  	3 \abs{ \frac{\p_z M(t,z) }{M(t,z)} }  \Big( \norm{\vphia\p_z^2 \Xi_\e}_\infty^{n} + O^\ast(1) + O(\e) \norm{(\kappa_\e,W_\e)} \Big) \right] .
\end{multline*}
We will use once more our hypothesis \cref{decay Gamma}, under the form stated in \cref{eq pz3lin}.  We also need all our previous estimates on the rings, \cref{rings Me,rings Ne,rings pzzWe}. We then  obtain 
\begin{multline*}
\e^2 \p_t \Big[ \vphia \p_z^3 W_\e\Big] (t,z) \leq  M(t,z)  \left[   -\vphia(t,z)\p_z^3 W_\e(t,z) + \frac 1{2^{2-\alpha}} \norm{\p_z^3 W_\e }_\infty^{n-1} \right.     \\ \left. +  \frac{1}{2^{1-\alpha}} \norm{\vphia \p_z^3 W_\e}_\infty  + O_0^\ast(1) + O(\e^\alpha)\norm{(\kappa_\e,W_\e)}   \right] .
\end{multline*}
We recall that the term $\norm{\vphia \p_z^3 W_\e}_\infty$ is a control on the whole space $\R$ and not only on the ball $B_0$. By applying the maximum principle, one gets
\begin{multline*}
\norm{\vphia \p_z^3W_\e}_{\infty}^{n}  \leq \\  \max \left( \frac 1{2^{2-\alpha}} \norm{\p_z^3 W_\e }_\infty^{n-1}  +  \frac{1}{2^{1-\alpha}} \norm{\vphia \p_z^3 W_\e}_\infty  + O_0^\ast(1) + O(\e^\alpha )\norm{(\kappa_\e,W_\e)} ,\norm{\vphia(0,\cdot) \p_z^3 W_\e(0,\cdot)}_\infty^n \right) .
\end{multline*}
We can absorb the initial data in the $O_0^\ast(1)$ to deduce:
\begin{align*}
\norm{\vphia \p_z^3W_\e}_{\infty}^{n}  \leq \frac 1{2^{2-\alpha}} \norm{\p_z^3 W_\e }_\infty^{n-1}  +  \frac{1}{2^{1-\alpha}} \norm{\vphia \p_z^3 W_\e}_\infty  + O_0^\ast(1) + O(\e^\alpha)\norm{(\kappa_\e,W_\e)}.
\end{align*}
This sequence is bounded, because its ratio verifies: $2^{\alpha-2} <1$.
\begin{align}\label{eq pz3 stab}
\norm{\vphia \p_z^3 \Xi_\e }_{\infty}^{n} & \leq \max \left( \frac{ O_0^\ast(1) +  O(\e^\alpha) \norm{(\kappa_\e,W_\e)} }{1-2^{\alpha-2}} + \frac{2^{\alpha-1}}{1-2^{\alpha-2}} \norm{\vphia \p_z^3 W_\e}_\infty , \norm{\vphia \p_z^3 W_\e }_{\infty}^{0} \right).
\end{align}
We define $k(\alpha)$ as follows:
\begin{align*}
\boxed{ k(\alpha) : = \frac{2^{\alpha-1}}{1-2^{\alpha-2}}  } 
\end{align*}
and from \cref{eq pz3 stab} we finally conclude, taking the initial data \cref{bound pzzzW B0} into account:
\begin{align*}
\norm{\vphia \p_z^3 \Xi_\e }_{\infty}^{n}  \leq O_0^\ast(1) +  O(\e^\alpha) \norm{(\kappa_\e,W_\e)} + k(\alpha) \norm{\vphia \p_z^3 W_\e}_\infty .
\end{align*}
We have therefore proven the following proposition :
\begin{prop}[In the rings, $\p_z^3 W_\e$]\label{rings pzzzWe}$ $\\
	There exists a constant $\e_B$ that depends only on $B$ such that upon the condition of \cref{stab We}, $W_\e$ verifies for $\e\leq \e_B$
	\begin{align*}
	\norm{\p_z^3 W_\e}_{\infty}^{n} \leq O_0^\ast(1) +  O(\e^\alpha) \norm{(\kappa_\e,W_\e)}+ k(\alpha) \norm{ \vphia \p_z^3 W_\e}_\infty,
	\end{align*}
	for $n \geq 1$, with \begin{align}
0< k(\alpha) : =  \frac{2^{\alpha-1}}{1-2^{\alpha-2}}  <1. 
	\end{align}
\end{prop}
The scalar $k(\alpha)$ is a contraction factor, only upon the condition
\begin{align}
\alpha < 2 - \frac{\ln 3}{\ln 2} \approx 0.415.
\end{align}
We make that assumption retrospectively when we introduce $\E$ in \cref{def F}. It appears to be the same threshold than in the stationary case, see \cite[Equation 5.11]{spectralsex}. It appeared in that case for seemingly very different reasons than here. Another reason for which $\alpha$ cannot be taken too large is that is worsens the contraction estimate $\vphia(t,z) \leq 2^\alpha \vphia(t, \bz)$.
\subsection{Conclusion : proof of \cref{stab We}}$ $\\
All our previous estimates of \cref{rings Me,rings Ne,rings pzzWe,rings pzzzWe}  are uniform in $n$, and therefore apply to the whole space. Thus, every bound of \cref{stab We} has been proved except for the one upon $\p_z W_\e$. Its proof can be straightforwardly adapted of the one of \cref{rings Ne}, starting from the linearized equation of \cref{prop lin pzWeps}. A more elegant argument is to notice that we dispose of the following uniform bound for all times $t>0$ and $z\in \R$:
\begin{align*}
\p_z \Xi_\e(t,z) = \p_z W_\e(t,z) - \p_z W_\e(t,\bz) \leq \frac{O_0^\ast(1) + O(\e) \norm{(\kappa_\e,W_\e)}}{\vphia(t,z)}.
\end{align*}
Therefore, since $\p_z W_\e(t,\zs)=0$,  we get that, by means of a  series, for all $h\in \R$ :
\begin{align*}
\p_z W_\e(t,\zs +h) \leq \Big(O_0^\ast(1) + O(\e)\norm{(\kappa_\e,W_\e)}\Big) \sum_{k\geq 0} \frac{1}{\vphia(t, \zs + 2^{-k}h)}.
\end{align*}
The series $\ds  \sum_{k \geq 0} 2^{\alpha k}$ converge, and therefore 
\begin{align*}
\norm{\p_z W_\e}_{\infty} \leq O_0^\ast(1) +  O(\e)\norm{(\kappa_\e,W_\e)} .
\end{align*}
One sees that if $\alpha=0$, the series above does not converge. This shows that the weight $\vphia$ is necessary to ensure uniform Lipschitz bounds of $W_\e$.
\section{Proof of \cref{main theo}}\label{sec proof}$ $\\
We now prove the main result of this paper, that is the boundedness of $(\kappa_\e,W_\e)$ in $\R \times \F$. 
We first suppose that there exists $K_0$ such that \begin{align}\label{sup K}
\abs{\kappa_\e(0)}\leq K_0 \quad \text{ and } \normf{W_\e(0,\cdot)} \leq K_0,
\end{align}
and we look to prove 
\begin{align*}
\abs{\kappa_\e} \leq K'_0 \quad \text{ and } \normf{W_\e} \leq K'_0,
\end{align*}
with $K$ to be determined by the analysis.

By \cref{stab We}, that we can apply with our assumption \cref{sup K}, we have precise bounds of   $W_\e$.  More precisely, there exists a constant $C_0^\ast$ that depends only on $K_0$ and $K^\ast$ and a constant $C'_K$ that depends only $K'_0$, such that :
\begin{align*}
\normf{W_\e}  \leq  C_0^\ast + C'_K\e^\alpha K'_0 + k(\alpha)  \norm{W_\e}_\F.
\end{align*}
Therefore, up to renaming the constants, 
\begin{align}\label{eq cont}
\normf{W_\e}  \leq  \frac{ C_0^\ast +C'_K \e^\alpha K'_0}{1-k(\alpha)}\leq C_0^\ast +C'_K \e^\alpha K'_0.
\end{align}
Now we work on $\kappa_\e$. We go back to \cref{prop eq kde} since we made suitable assumptions and we get that $\kappa_\e$ solves
\begin{align}\label{EDO kappa2}
-\kde(t)  =   R_\e^\ast(t)  \kappa_\e   + O^\ast(1)+O(\e)\norm{(\kappa_\e,W_\e)} + O^\ast(1) \norm{W_\e}_\F. \end{align}
Thanks to our previous contraction argument,  we have an estimate  of the term $\norm{W_\e(t,\zs)}_\F$.
Keeping in mind this estimate \cref{eq cont}, we can finally conclude the argument on $\kappa_\e$. 

Since $R_\e^\ast$ is a positive function  that admits for $t\geq t_0$ a uniform lower bound $R_0$,  see \cref{control pzgVIea}, it is straightforward from \cref{eq cont} and  \cref{EDO kappa2}, and our subsequent bounds, that there exists $C_0^\ast$ and $C'_K$  such that for all time $t$
\begin{align}\label{final kappa}
\abs{\kappa_\e(t)} \leq C_0^\ast + C'_K \e^\alpha K'_0. 
\end{align}
Coupled with \cref{eq cont}, those are the stability results we needed. 
Set a scalar  $K$ such that
\begin{align}\label{choice K}
K'_0\geq 2C_0^\ast. 
\end{align}
Then, choose $\e_0$ in the following way
\begin{align*}
\e_0  :=  \left( \frac{1}{2C'_K}\right)^{\frac 1\alpha},
\end{align*}
where $C_K$ is the constant corresponding to the choice made in \cref{choice K} of the size of the ball $K$.
Then for $\e \leq \e_0$, starting from an initial data that verifies \cref{sup K}, the bound is propagated in time and
\begin{align*}
\boxed{\normf{W_\e} \leq K'_0,\quad 
\abs{\kappa_\e}\leq K'_0.}
\end{align*}
Since $V_\e = \Va+\e^2 W_\e, \quad q_\e=\qa+\e^2 \kappa_\e, $ \cref{main theo} is proven.

\section{Numerical simulations and discussion}\label{sec num}
In this section we will display some numerical simulations showing the behavior of the solution of the Cauchy problem for positive $\e$, and we will provide an insight on the structural assumption we made in \cref{cond Gamma}. 

\begin{figure}[H]
	\setlength{\tabcolsep}{8pt}	\begin{tabular}{cc}
		\raisebox{2pt}	{\includegraphics[scale=0.5]{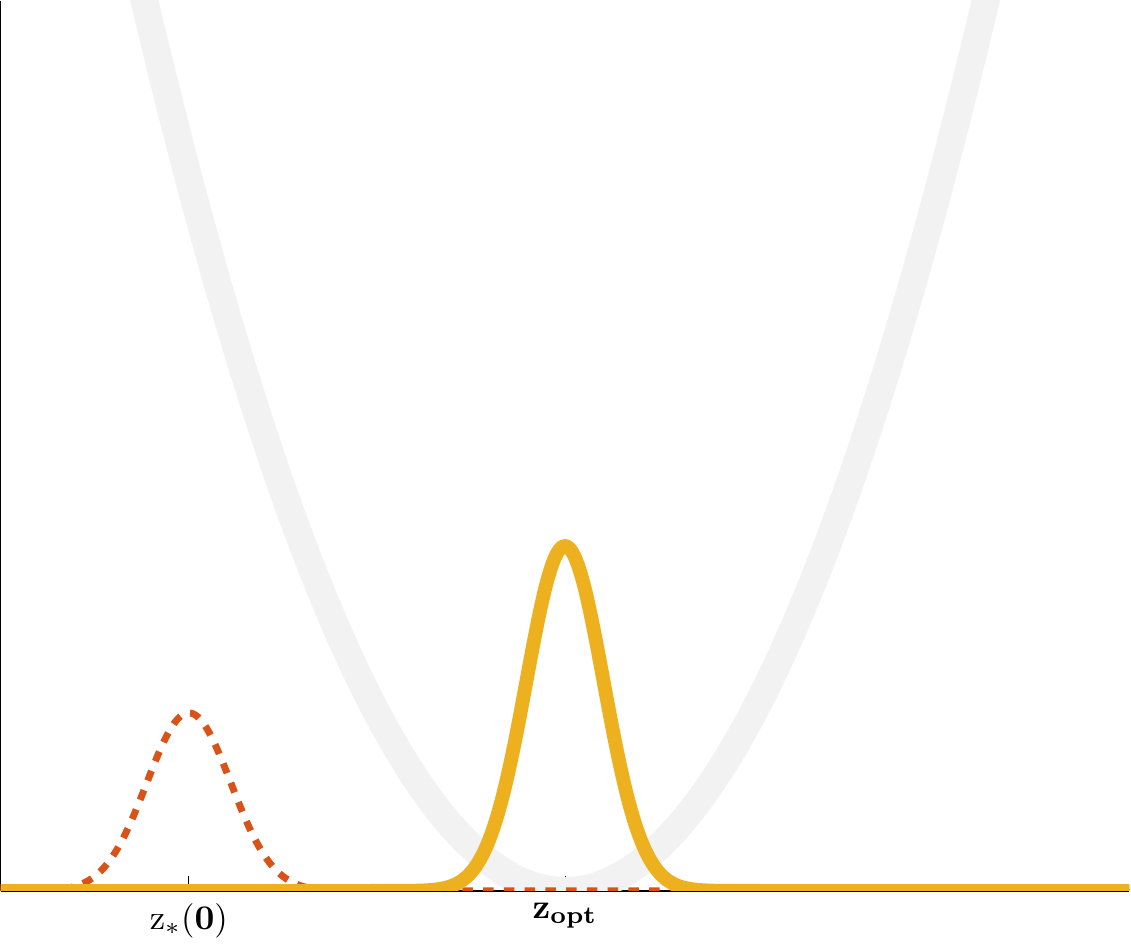}} &	\includegraphics[scale=0.42]{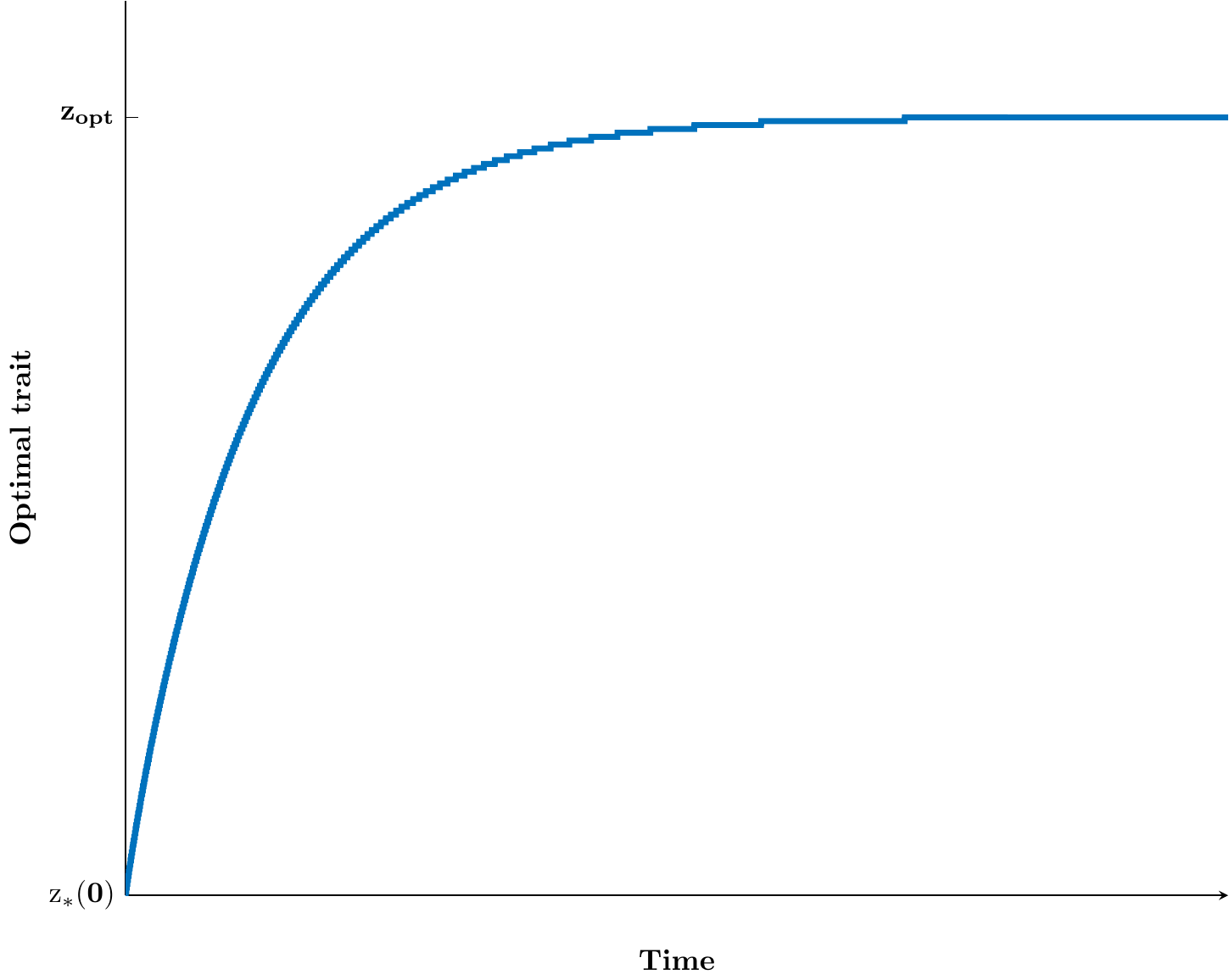}
	\end{tabular}
	\caption{On the left, in dotted red, the initial data $f_\e(0,\cdot)$, and in orange the distribution $f_\e$ after a long time. In the background the selection function $m$ with a global optimum $z_{opt}$. On the right, the trajectory of the dominant trait $\zs$.}\label{quadratic}
\end{figure}

\begin{figure}[H]
	\setlength{\tabcolsep}{8pt}	\begin{tabular}{cc}
		\raisebox{8pt}	{\includegraphics[scale=0.5]{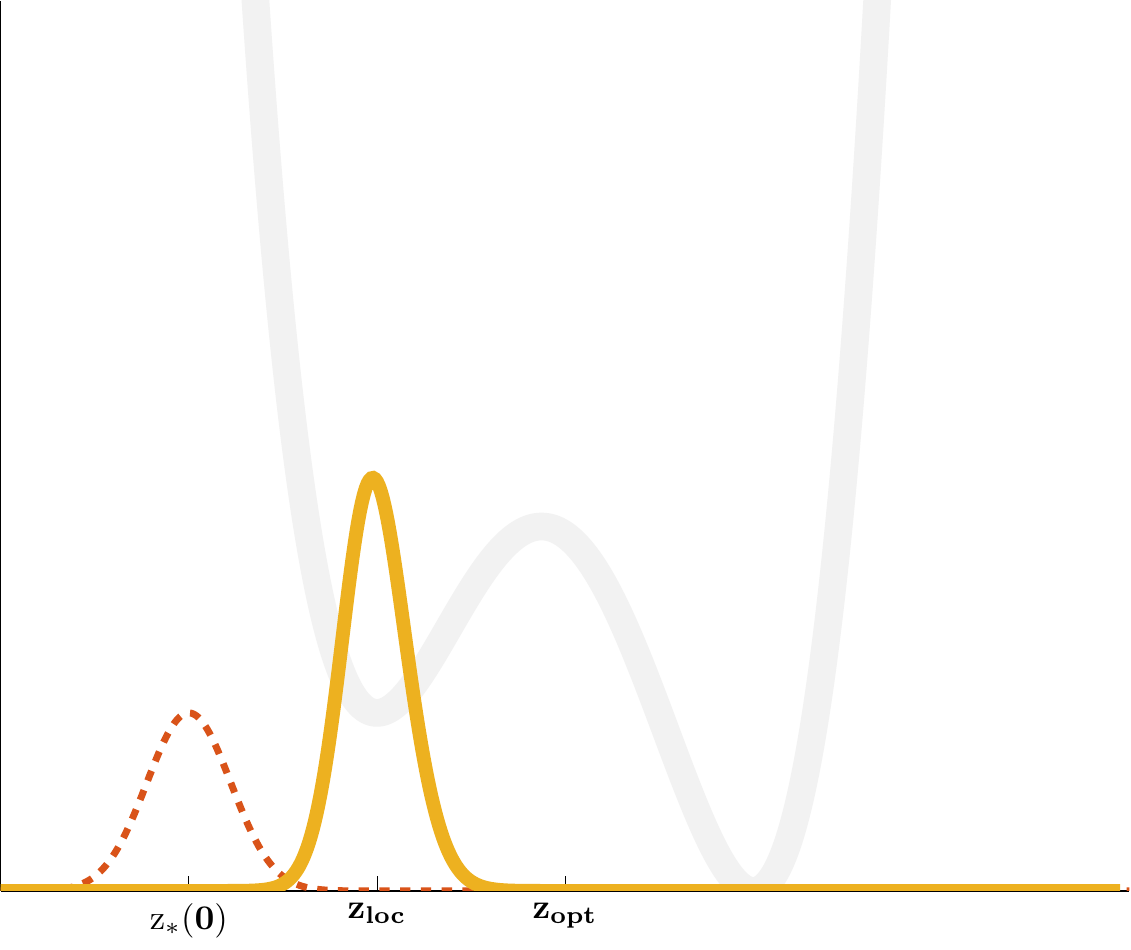}} &	\includegraphics[scale=0.5]{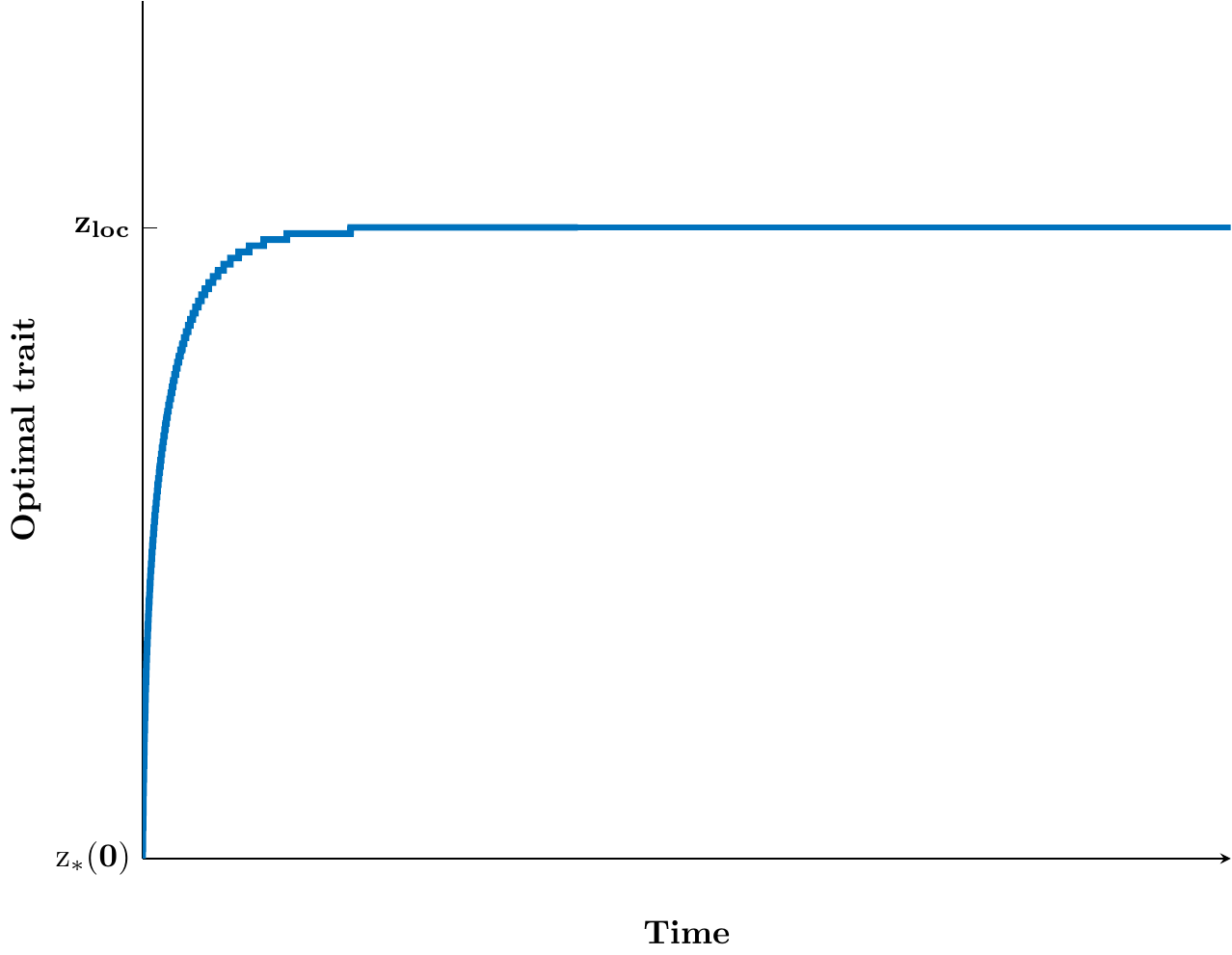}
	\end{tabular}
	\caption{On the left, in dotted red, the initial data $f_\e(0,\cdot)$, and in orange the distribution $f_\e$ after a long time. In the background the selection function $m$ with a global optimum $z_{opt}$ and a local optimum $z_{loc}$. On the right, the trajectory of the dominant trait $\zs$. The function $M$ is uniformly positive.}\label{normal}
\end{figure}

\textbf{Influence of the condition \cref{cond Gamma}.} A first example for our study is to consider quadratic selection function, as depicted in \Cref{quadratic}. In that case, according to \Cref{main theo},  starting from any initial data $\zs(0)$, the solution $f_\e$ stays  close to a Gaussian density with variance $\e^2$. In addition, its mean $\zs$ converges to the unique minimum of $m$ when the time is large.
Our framework encompasses more general selection functions with multiple local minima, as depicted in \Cref{normal}. The condition in \cref{cond Gamma} restricts somehow the position of those minima. If one assumes that $\zs$ starts from a local minimum, that is $m'(\zs(0))=0$, then this condition is that the selective difference between minima must be inferior to $1$ : $m(\zs(0))-m(z_{opt})<1$. We recover the structural condition under which the analysis for the stationary case was performed, see \cite{spectralsex}. 

The selection function depicted in \Cref{normal}, coupled with $\zs(0)$ verifies the condition \cref{cond Gamma}. Then as stated by \cref{main theo} the population density $f_\e$ concentrates around the local minimum, accordingly to the gradient flow dynamics of \cref{def m cauchy}.

\begin{figure}
\setlength{\tabcolsep}{8pt}	\begin{tabular}{cc}
\raisebox{8pt}	{\includegraphics[scale=0.35]{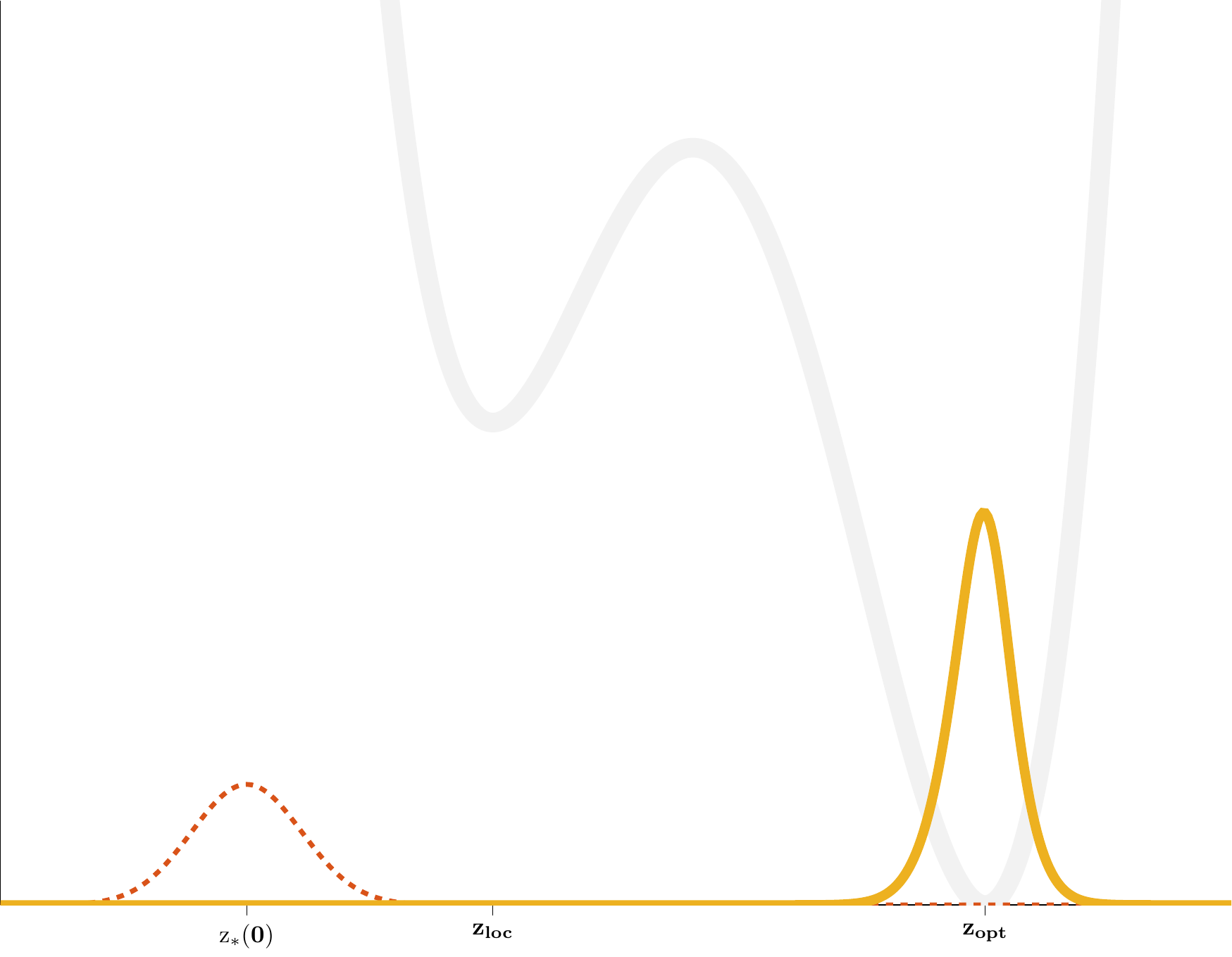}} &	\includegraphics[scale=0.51]{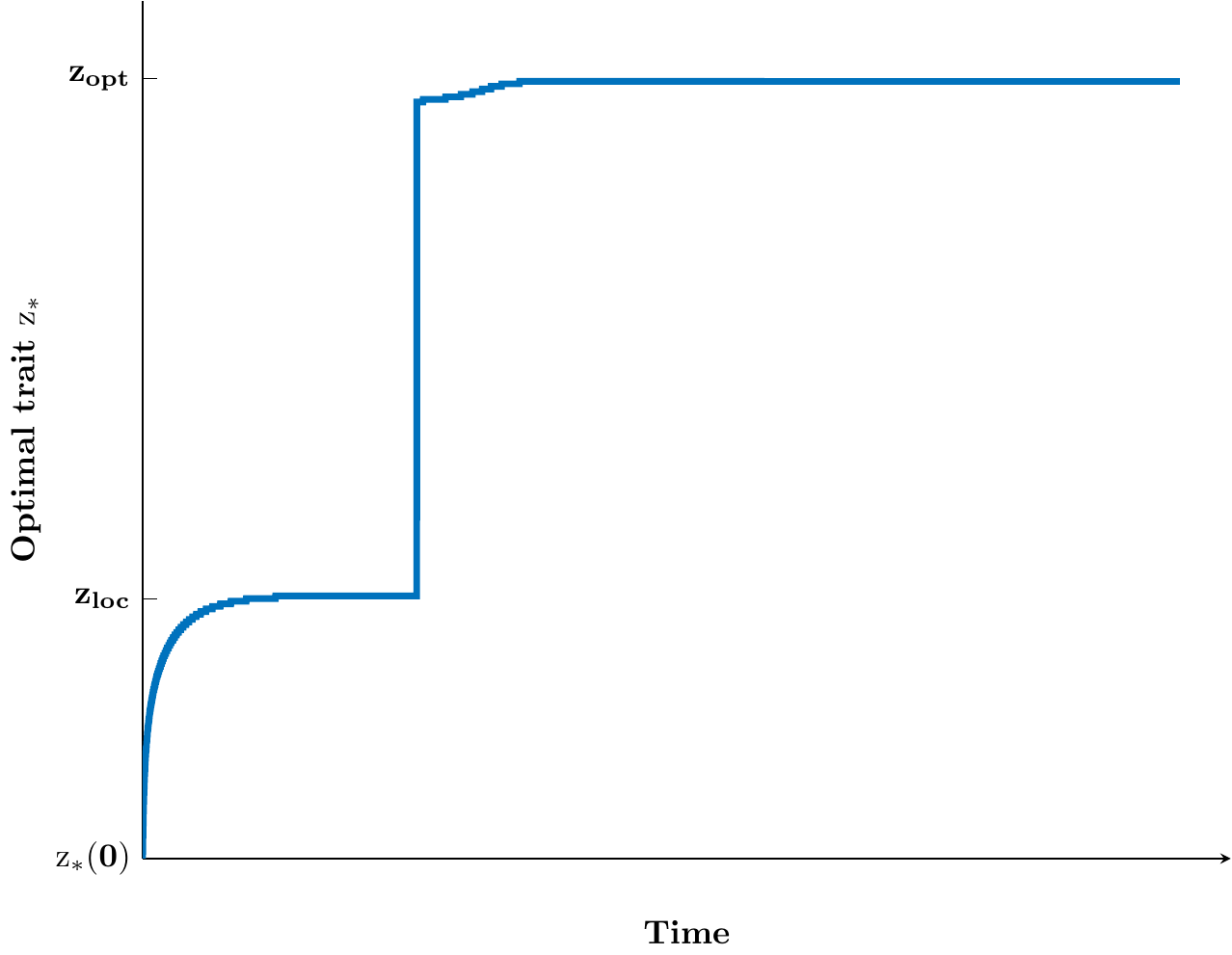}
\end{tabular}
	\caption{On the left, in dotted red, the initial data $f_\e(0,\cdot)$, and in orange the distribution $f_\e$ after a long time. In the background the selection function $m$ with a global optimum $z_{opt}$ and a local optimum $z_{loc}$. On the right, the trajectory of the dominant trait $\zs$. The function $M$ is not uniformly positive.}\label{critical}
\end{figure}

A case not taken into account by our methodology is when \cref{cond Gamma} is not verified at all times. This is the case if the slopes of the lines between local and global minima are too sharp. For instance, this is true in the case of  \Cref{critical}. Interestingly, what is observed is a critical behavior. The solution will first concentrate around the first local minimum before jumping sharply in the attraction basin of the global minimum see the right hand picture of \Cref{critical}.
Under this model it would seem that the population will concentrate around the global minimum of selection if it is much better than the other selective optima. Interestingly, the value of the local maximum in between the two minima, that could act as an obstacle between the two convex selection valleys, do not appear to play a role. On the other hand if the global minimum is not much better than a local minima, in the sense that each of them falls under the regime of \cref{cond Gamma}, the population can concentrate around this local minimum.	

\textbf{Influence of the sign of $q_\e$.} We introduced the scalar $q_\e$ in \cref{decomp Ue} as part of the decomposition of $U_\e$ between the affine parts and the rest of the function, which we later justified by heuristics on the linearized problem, see \Cref{tab}. We can propose a different interpretation of this scalar,related to the Gaussian distribution. 

The logarithmic transform \cref{hopfcole} coupled with the decomposition \cref{decomp Ue} can be rewritten as the following transform on the solution of \cref{Ptfeps} : 
\begin{align}\label{hopf skewed}
f_\e(t,z) = \frac{1}{\e \sqrt{ 2\pi}}\exp \left( \dfrac{\lambda(t) - \e^2 p_\e(t) + \e^4 q_\e(t)^2}{\e^2} - \dfrac{\Big(z - (\zs(t)-\e^2q_\e(t) )\Big)^2}{2\epsilon^2}  - V_\epsilon(t,z) \right). 
\end{align}
Therefore one can see that $q_\e$ is the correction to the mean of the Gaussian distribution at the next order in $\e$. Its sign corresponds to the sign of the error made on the mean of the Gaussian distribution. If $q_\e$ is positive, the correction of $\zs$ lies on its left. This is consistent  with the following reasoning on the limit value $\qa = \ds  \lim_{\e \to 0} q_\e$, defined in \cref{def qda}. For clarity, suppose that $\zs$ does not depend on time, that is the regime of the stationary case. Then from \cref{def qda}, we find an explicit value for $\qa$,  which coincides with  \cite[equation 3.2]{spectralsex} : 
\begin{align*}
\qa = \frac{m^{(3)}(\zs)}{2 m''(\zs)}.
\end{align*}
By local convexity of $m$ around $\zs$, see \cref{zs infty}, the sign of $\qa$ is the same than the sign of $m^{(3)}(\zs)$. Therefore, if this scalar is positive, selection leans the profile towards the left, which has better selective values than the right, since it is flatter. Therefore, we recover what we deduced from \cref{hopf skewed}, the sign of $q_\e$ is linked to the skewness of the selection function $m$ around $\zs$.

\bibliographystyle{apalike}
\bibliography{Biblio}

\end{document}